\documentclass[a4paper,12pt]{article}
\usepackage{amssymb,amsmath,amsfonts,amsthm}
\usepackage{a4wide}

\numberwithin{equation}{section}

\providecommand{\abs}[1]{\left\vert#1\right\vert}
\providecommand{\norm}[1]{\left\Vert#1\right\Vert}
\providecommand{\pnorm}[2]{\left\Vert#1\right\Vert_{L^{#2}}}
\providecommand{\pnormspace}[3]{\left\Vert#1\right\Vert_{L^{#2}(#3)}}
\providecommand{\pqnorm}[3]{\left\Vert#1\right\Vert_{L^{#2,#3}}}
\providecommand{\pqnormspace}[4]{\left\Vert#1\right\Vert_{L^{#2,#3}(#4)}}

\providecommand{\lzspace}[5]{\left\Vert#1\right\Vert_{L^{#2,#3}\log^{#4}{L}(#5)}}

\providecommand{\Rn}[1]{\mathbb{R}^{#1}}

\providecommand{\csubset}{\subset\subset}

\def\wstar{\overset{*}{\rightharpoonup}}

\def\({\left(}
\def\){\right)}
\def\l|{\left|}
\def\r|{\right|}
\def\ep{\varepsilon}
\def\mr{\mathbb{R}}

\def\mc{\mathbb{C}}

\def\si{\mathbb{S}^1}
\def\p{\partial}
\def\lep{|\mathrm{log }\ \ep|}
\def\llep{\log\ \lep}
\def\nab{\nabla}
\def\np{\nab^\perp}

\def\om{\Omega}
\def\io{\int_{\Omega}}
\def\bo{\partial \Omega}
\def\hal{\frac{1}{2}}

\def\lti{L^{2,\infty}}
\def\dist{\mathrm{dist}}
\def\he{h_\text{ex}}
\def\hci{{H_{c_1}}}

\def\vchi{\text{\large{$\chi$}}}

\DeclareMathOperator{\curl}{curl}
\DeclareMathOperator{\diverge}{div}
\DeclareMathOperator{\supp}{supp}

\newtheorem{lem}{Lemma}[section]
\newtheorem{cor}[lem]{Corollary}
\newtheorem{prop}[lem]{Proposition}
\newtheorem{thm}{Theorem}
\newtheorem{remark}{Remark}[section]

\title{Lorentz Space Estimates and Jacobian
  Convergence for the Ginzburg-Landau Energy with Applied Magnetic Field}
\author{Ian Tice\footnote{Supported by an NSF
Graduate Research Fellowship}\\
{\small Courant Institute of Mathematical Sciences}\\
{\small 251 Mercer St., New York, NY 10012}\\
{\small\tt tice@cims.nyu.edu} }

\begin{document}

\maketitle
\begin{abstract}
In this paper we continue the study of Lorentz space estimates
for the Ginzburg-Landau energy started in \cite{p1}. We focus on
getting estimates for the Ginzburg-Landau energy with 
external magnetic field $h_{ex}$ in certain interesting regimes of $h_{ex}$.
This  allows us to show that for configurations close to
minimizers or local minimizers of the energy, the vorticity mass
of the configuration $(u,A)$ is comparable to the $L^{2,\infty}$
Lorentz space norm of $\nabla_A u$.  We also establish 
convergence of the gauge-invariant Jacobians (vorticity measures)
in the dual of a function space defined in terms of Lorentz
spaces.

\end{abstract}

\section{Introduction}
This is the sequel to the paper \cite{p1}, where we proved Lorentz
space estimates for the Ginzburg-Landau free energy
\begin{equation}
\label{fep}  F_\ep(u,A)= \hal \io|\nab_A u|^2 + |\curl A|^2 +
\frac{1}{2\ep^2}(1-|u|^2)^2.
\end{equation}
In the present paper we consider  the full
Ginzburg-Landau energy with applied magnetic field
\begin{equation}
\label{gep} G_\ep (u,A) = \frac{1}{2} \int_{\Omega} \abs{\nabla_A
u}^2  + \abs{\curl{A} - \he}^2+
\frac{1}{2\varepsilon^2}(1-\abs{u}^2)^2,
\end{equation} which models a superconductor submitted to an
external magnetic field of intensity $\he$. In \eqref{fep} and
\eqref{gep} $\Omega\subset \Rn{2}$ is a bounded regular domain, and
$u$ is a complex-valued function called the ``order parameter," which indicates
the local state of the material (normal or superconducting): $|u|^2$ is the local density of
superconducting electrons.  The vector field  $A:\Omega \rightarrow \Rn{2}$ is the
vector-potential of the induced magnetic field,  $h:=\curl A=\p_1 A_2 - \p_2 A_1$.
The notation $\nab_A $ refers to the
covariant gradient $\nab_A u = (\nab - i A) u$. We are interested
in the regime of small $\ep$, corresponding to ``extreme type-II"
superconductors.

The Ginzburg-Landau energy with magnetic field admits a
gauge-invariance: for every smooth $\Phi$, $G_\ep(u,A) = G_\ep(u
e^{i \Phi} , A + \nab \Phi)$. The physically intrinsic quantities
are those that are gauge-invariant, such as $|u|$ and $|\nab_A
u|$.

   We refer to \cite{ss_book} for a more thorough presentation of this
functional.

\subsection{Results of \cite{p1}}
The objects of interest are the zeroes of the complex-valued
function $u$, which can have a nonzero topological degree. These
are called the {\it vortices} of the configuration $(u,A)$.

Starting with Bethuel-Brezis-H\'elein \cite{bbh}, several studies
have shown how to relate the value of the energy to the vortices
and their degrees. The method we focused on was
 the ``vortex-ball construction" introduced by Jerrard \cite{j} and
 Sandier \cite{sa}, which allows to construct disjoint  ``vortex balls" of
 small size and of degree $d$ containing at least a  $\pi |d|\lep$
  contribution to the energy. It was explained in \cite{p1} that the
  typical profile of a vortex  of degree $d$ is $ f(r)
  e^{id\theta} $ in polar coordinates, with $f(0)=0$ and $f(r)$
 very close to $1$ as soon as $r \gg \ep$; as a result $|\nab_A u|$ typically
 blows up like $ d/r$, leading to a logarithmic divergence of its
 $L^2$ norm,  hence of the energy. On the other hand,
  considering the Lorentz norm defined by
\begin{equation}
\label{ltinorm1} \|f\|_{\lti} = \sup_{|E|<\infty}
|E|^{-\hal}\int_E |f(x)|\, dx,
\end{equation}
where $|E|$ denotes the Lebesgue measure of $E$, one observes that
the $\lti$ norm of $d/r$ does not blow up, but is instead of order
$2\sqrt{\pi} |d|$. $\lti$ is critical in the sense that it is the
smallest Lorentz space to which the profile $1/r$ belongs. Based
on this observation, we searched for estimates on $\|\nab_A
u\|_{\lti}$ that would not blow up with $\ep\to 0$  but rather
would be of the order of the total vorticity mass $\sum |d_i|$ and
could thus serve to estimate the total number of vortices.

The method used in \cite{p1} consisted in giving an improvement of
the lower bounds of \cite{j,sa,ss_book},
 which allowed us to  gain an extra term that served to evaluate
 $\|\nab_A u\|_{\lti}$.  Writing
\begin{equation*}
F_\varepsilon(\abs{u},\Omega):= \frac{1}{2}\int_\Omega \abs{\nabla \abs{u}}^2 + \frac{1}{2\varepsilon^2}(1-\abs{u}^2)^2,
\end{equation*}
for the free energy of $\abs{u}$, our first main result was

\begin{thm}[Improved ball construction]
\label{energy_bound} Let $\alpha \in (0,1)$.  There exists
$\varepsilon_0 >0$ (depending on $\alpha$) such that for
$\varepsilon \le \varepsilon_0$ and $u, A$ both $C^1$ such that
$F_{\varepsilon}(\abs{u},\Omega) \le \varepsilon^{\alpha-1}$,  the
following hold.

For any $1 > r > C \varepsilon^{\alpha/2}$, where $C$ is a
universal constant, there exists a finite, disjoint collection of
closed balls, denoted by $\mathcal{B}$, with the following
properties.

1. The sum of the radii of the balls in the collection is $r$.

2.  Defining $\Omega_\ep= \{ x \in \Omega \, | \,\dist(x, \bo) >
\ep\}$, we have 
\begin{equation*}
\{x\in \Omega_{\varepsilon} \; | \; \abs{u(x)-1} \ge  \varepsilon^{\alpha/4} \}
\subset V := \Omega_{\varepsilon} \cap \left(\cup_{B\in
\mathcal{B}} B\right).
\end{equation*}

3. We have
\begin{multline}\label{e_b_15}
\frac{1}{2} \int_V \abs{\nabla_A u}^2 + \frac{1}{2\varepsilon^2}(1-\abs{u}^2)^2 + r^2(\curl{A})^2
 \\  \ge \pi n \left( \log{\frac{r}{\varepsilon n}} -C  \right) +
\frac{1}{18} \int_{V} \abs{\nabla_{A}u -i u Y}^2 +
\frac{1}{2\varepsilon^2}(1-\abs{u}^2)^2,
\end{multline}
where $Y$ is some explicitly constructed vector field, $d_B$ denotes
$\deg(u, \p B)$ if $B \subset \Omega$ and $0$ otherwise, and
\begin{equation*}
n = \sum_{\substack{B\in \mathcal{B} \\ B \subset
\Omega_\varepsilon}} \abs{d_B}
\end{equation*}
is assumed to be nonzero and $C>0$ is universal.
\end{thm}

\begin{remark}
 In the earlier paper \cite{p1} we denoted the vector field $Y$ by $G$.  We switch notation here to avoid confusion with $G_\varepsilon$, the energy functional.  In what follows we will construct such a field $Y$ for each configuration $\{(u_\varepsilon,A_\varepsilon)\}$ in a sequence; when we do so we will write $Y_\varepsilon$ to denote this dependence.
\end{remark}

Then we could bound from below $\int|\nab_{A}u  -iuY|^2 $ by
$\|\nab_{A} u -iuY\|_{\lti}^2$, and by controlling $\|Y\|_{\lti}$, we
obtained:

\begin{thm}[Lorentz norm bound]
\label{norm_switch} Assume the hypotheses of Theorem
\ref{energy_bound}.  Then there exists a universal constant $C>0$
such that
\begin{multline}\label{n_s_1}
\frac{1}{2} \int_V \abs{\nabla_A u}^2 +
\frac{(1-\abs{u}^2)^2}{2\varepsilon^2}
 + r^2(\curl{A})^2 + \pi \sum  |d_{B}|^2\\
\ge C \pqnormspace{\nabla_A u}{2}{\infty}{V}^2 + \pi \sum |d_B|
\left( \log{\frac{r}{\varepsilon \sum |d_B|}} -C  \right).
\end{multline}
where the sums are taken over all the balls $B$ in the final
collection $\mathcal{B}$ which are included in $\Omega_\ep$.
\end{thm}
This result allows us to control $\|\nab_A u\|_{\lti}$ by the
``energy-excess," the difference between the total energy and the
vortex energy. In \cite{p1} we presented several corollaries, and in particular, direct applications to minimizers $u_\varepsilon$ of the Ginzburg-Landau functional without magnetic field: we bounded $\|\nab u_\ep\|_{\lti(\om)}$  in terms of the
total degree and then deduced  $\lti(\om)$ weak-$*$ convergence results for $\nab
u_\ep$.

\subsection{The Ginzburg-Landau energy with applied magnetic field}

Applying Theorem \ref{norm_switch} to get useful estimates for the
full functional (\ref{gep}) with magnetic field is a more
complicated task. Writing $n =\sum |d_B|$ for the total vorticity, which implicitly depends on $\varepsilon$, we may search for estimates of the type $\|\nab_A u\|_{\lti} \le Cn
$. Such estimates follow from Theorem \ref{norm_switch} if an
upper bound on the free energy like $F_\ep(u,A)\le \pi n \lep
+O(n^2)$ holds, but this is in general not true for arbitrary
configurations, and not even true for energy-minimizers. The
reason is that when there is an applied magnetic field, the vortices
tend to be confined near the center of the sample by the magnetic
field; this is so because the energy of interaction is not only proportional to
$n^2$, but also contains the cost of the interaction of confined vortices. We must compensate by extracting
a new term, to be used in the Lorentz space norms, in each step
used in the proofs of \cite{ss_book}. Because this term will be
inserted between matching lower and upper bounds, the result will
only  be true for configurations whose energy is close to optimal.
In particular they will apply to various minimizers and locally
minimizing solutions found in \cite{ss_book}. \medskip

Let us now look more closely at the way of describing the vortices. In
\cite{ss_book} as well as in previous papers, the vortices of a
configuration were described through its ``vorticity" $\mu(u,A)$,
a gauge-invariant version of the Jacobian determinant of $u$:
\begin{equation}
\label{jacdef}\mu(u,A) = \curl (iu,\nab_A u)+ \curl A,
\end{equation} 
where the vector field $(iu,\nab_A u)$ is called the current, and $(\cdot,\cdot)$ denotes the scalar product in
$\mc$ as identified with $\mr^2$, i.e. $(iu,\nab_A u)= \Im
(\overline{u} \nab_A u)$. This is an intrinsic and gauge-invariant
quantity that is analogous to the vorticity in fluid mechanics.
It can be related to ``vorticity measures" $\sum_i d_i \delta_{a_i} $ obtained via the
ball-construction (like the result of Theorem \ref{energy_bound}),
where $a_i$'s are the centers of the balls and $d_i$'s their
degrees, via the following ``Jacobian estimates" (see the work of Jerrard-Soner \cite{js},
or Theorem 6.2 of \cite{ss_book} ):
\begin{equation} \label{jacest}
\|\mu(u,A) - 2\pi \sum_i d_i
\delta_{a_i}\|_{(C^{0,\gamma}(\Omega))^*} \le  C r^\gamma
(F_\ep(u,A) +1).
\end{equation}
That is, the vorticity measures
constructed via the ball construction -- nonunique and nonintrinsic -- are
very close to the intrinsic vorticity $\mu(u,A)$ when the total
radii of the balls is small. Then, after normalizing by
 the possibly divergent $n=\sum_i |d_i|$,
  these measures are weakly compact in the sense
of measures, and this yields that $\mu(u,A)$, similarly normalized
by $n$, is compact in $(C^{0,\gamma})^*$, and after
extraction, converges to a {\it measure}.

The problem with this normalizing factor $n$ is that it depends on
the ball-construction, and thus is not intrinsic. This is where
the introduction of  $\|\nab_A u\|_{\lti}$ may help since it is
an intrinsic quantity expected  to behave like $n$; we
will make this rigorous.

\subsection{Regimes of applied field}

We will employ the notation $a \ll b$ to mean that $a/b \rightarrow 0$ as $\varepsilon \rightarrow 0$.  We write $o_z(1)$ (resp. $O_z(1)$) for a quantity, depending on $z$, that vanishes (resp. is bounded) as either $z \rightarrow 0$ or $z\rightarrow \infty$ (it will always be clear in context what the limit of $z$ is).  We write $o_z(a)$ and $O_z(a)$ for quantities such that $o_z(a)/a = o_z(1)$ and $O_z(a)/a = O_z(1)$ respectively.  The symbols $o(1)$ and $O(1)$ always mean $o_\varepsilon(1)$ and $O_\varepsilon(1)$.  Let us now recall some of the results summarized in  \cite{ss_book}: Sandier
and Serfaty, studying minimizers of the energy (\ref{gep}) for all
applied magnetic fields satisfying $ \he \ll \frac{1}{\ep^2}$ as
$\ep \to 0$, showed that there are essentially four regimes, which
follow some physical ``phase-transitions" (all constants below are
positive).
\begin{enumerate}
\item  For $\he\le \hci,$ minimizers have no vortices, and $\hci$
is the first critical field, which has an asymptotic expansion
$\hci \sim C(\Omega) \lep+ O(1) $ as $\ep\to 0$, $C(\Omega)$ being
 a constant determined by the domain.  
\item  For $\hci \le \he
\le \hci + O(\llep)$ minimizers have a bounded number  $n$ of vortices,
determined by the value of $\he$. 
\item For $\llep \ll \he - \hci
\ll \lep$ there are roughly $n$ vortices, with $n$ determined by $\he$ and
 $1\ll n \ll \he$ as $\ep \to 0$. 
\item For $ C \lep \le \he -
\hci \ll \frac{1}{\ep^2}$ there are roughly $n$ vortices with $
C_1\he \le n \le C_2 \he $.\end{enumerate}

For each regime  the asymptotic value of the minimal energy was
given, and  the optimal limiting vorticities were identified,
through explicit limiting problems obtained via
$\Gamma$-convergence. In regime 2, the vortices tend to minimize a
function of their $n$ locations. In regimes 3 and 4, the vorticity
$\mu (u,A)$, suitably blown-up and normalized by $n$, converges to
some identified probability measure with constant density.

 We will focus here on the regimes 2 and 3 where  $n \ll \he$.  The reason is that regime
4 is easy to treat. Indeed in that regime, $n$ and $\he$ are of
the same order, thus the natural normalization of  $\mu(u,A)$ is by 
$\he$, a quantity that does not depend on the vortex-ball construction.
This is what was done in \cite{ss_book}, Chapter 7.  Moreover, the a priori upper bound
$G_\ep(u,A) \le C \he^2$ is always satisfied for minimizers
(comparing with $u=1, A=0$), and thus $\|\nab_A u\|_{L^2} \le C \he$
always holds. Therefore in the regime 4,  the desired control
$\|\nab_A u\|_{\lti} \le C n$ comes trivially. A lower bound for
it by $n$ times the norm of a weak limit  also  follows. As for
the Jacobian vorticity, its $H^{-1}$ compactness is proved in
\cite{ss_book}, which is stronger than what one can prove with
norms involving Lorentz spaces. 

For these reasons, we now focus on regimes 2 and 3.

\subsection{The result of \cite{ss_book}, Chapter 9}

Let us first recall some notation and results from \cite{ss_book}
in these regimes. We introduce $\xi_0$ to be the solution to
\begin{equation}\label{xi0}
 \begin{cases}
  -\Delta \xi_0 + \xi_0 +1 = 0 & \text{ in } \Omega \\
  \xi_0 =0& \text{ on } \partial \Omega.
 \end{cases} 
\end{equation}
The significance of $\xi_0$ is that it conveys the geometry of the
domain to the vortex confinement by the applied magnetic field: when $n \ll
h_{ex}$, the $n$ vortices tend to nucleate in a neighborhood of
the set on which $\xi_0$ achieves its minimum, a set determined by
the geometry of the domain. It can be shown (see \cite{ss_glmin})
that this set is composed of a finite number of points; for
simplicity we assume that the set is a single point $p$, i.e. that
$\xi_0$ has a unique minimum at the point $p\in \Omega$. We shall further assume that $D^2 \xi_0(p)$ is positive definite.  We also
define, for any vector-field $A$,  
\begin{equation*}
A':=A- \he \np\xi_0.
\end{equation*}

In the regime $\he -\hci \ll \lep$, all the vortices concentrate
around the point $p$, at the scale $\ell= \sqrt{\frac{n}{\he}}$
(where $n$ again is the number of vortices).  Thus the vorticity
measure $\mu(u,A)$ normalized by $2\pi n$ will converge to $
\delta_p$, the Dirac mass at $p$. In order to obtain more
interesting information on the vortex-locations, we need to blow
up the vorticity measure $\mu(u,A)$: its push forward under the
rescaling $x\mapsto \sqrt{\frac{\he}{n}}(x-p)$ is denoted
$\tilde{\mu}(u,A)$. We will also use the  function $G_p$,  defined
as the solution of
\begin{equation}\label{Gp}
 \begin{cases}
  -\Delta G_p + G_p  = 2\pi \delta_p & \text{ in } \Omega \\
  G_p =0& \text{ on } \partial \Omega.
 \end{cases}
\end{equation}
It was proved that the energy in that regime is equivalent to
$f_\ep(n)+O(n^2)$ where $f_\ep$ is an explicit function of $n$,
depending only on $\he$, $\ep$,  and the domain $ \om$. When $n\gg
1$, i.e. in the regime 3, we have  the following
$\Gamma$-convergence result to the function $I$ defined over the
set of probability measures on $\mr^2$ by
\begin{equation}\label{I}
I(\mu)= - \pi \iint_{\mr^2\times\mr^2} \log |x-y|\, d\mu(x) \,
d\mu(y) + \pi \int_{\mr^2} Q(x) \, d\mu(x),\end{equation}
 where $Q$ is the  quadratic form of the Hessian of $\xi_0$ at the point $p$.

\begin{thm}[$\Gamma$-convergence in the intermediate regime -
\cite{ss_book}, Theorem 1.5]\label{ss} \indent\par
Let  $\{(u_\ep,A_\ep)\}_\ep$
be a family of configurations such that
$G_\ep(u_\ep,A_\ep)<\ep^{-1/4} $ with $\he <C \lep$. Defining
$n =\sum_i|d_i|$, where the $d_i$'s are the degrees of some
collection of  vortex-balls of total radius
$r=\frac{1}{\sqrt{\he}}$ constructed by Theorem
\ref{energy_bound},  assume that
$$1\ll n \ll \he$$  and $G_\ep(u_\ep, A_\ep)\le f_\ep(n) + Cn^2$,
 as $\ep \to 0$. Then
there exists a probability measure $\mu_*$  such that, after
extraction of a subsequence, $\frac{\tilde{\mu}(u_\ep,A_\ep)}{2\pi
n} \to \mu_*$ in $(C^{0,\gamma}_c(\mr^2))^*$ for some
$1\ge \gamma>0$ and, as $\varepsilon \rightarrow 0$, 
$$G_\ep(u_\ep, A_\ep)-f_\ep(n)\ge n^2 I(\mu_*) +o(n^2).$$
Conversely, for each probability measure $\mu$  with compact
support in $\mr^2$ and each $1\ll n \ll \he\le C \lep$, there
exists $\{(u_\ep, A_\ep)\}_\ep$ such that
$\frac{\tilde{\mu}(u_\ep,A_\ep)}{2\pi n} \to \mu_*$ in
$(C^{0,\gamma}_c(\mr^2))^*$ for each $1\ge \gamma>0$ and such that, as $\varepsilon \rightarrow 0$, 
$$G_\ep(u_\ep,A_\ep) - f_\ep(n) \le n^2 I(\mu) +o(n^2).$$
\end{thm}

An analogous result was proved in \cite{ss_book} for the case
$n = O(1)$ (regime 2), which we do not quote here for the sake of
brevity.

\subsection{Main result}

Here we obtain several improvements of this result. We quote here
some of our results under simpler assumptions; more results, with a weaker set of assumptions, can be found in the theorems below.  A first improvement is obtained under this same assumption
 \begin{equation}
 \label{hyp1}
G_\varepsilon(u_\varepsilon,A_\varepsilon) \le
f_\varepsilon(n) + C_0 n^2\end{equation} for some constant
$C_0\ge 0$, where $n$ is defined as in the theorem above.

 \begin{thm} \label{3}
 Suppose configurations $\{(u_\varepsilon,A_\varepsilon)\}$
satisfy
\begin{equation*}
 F_\varepsilon(u_\varepsilon,A_\varepsilon') \le \varepsilon^{\alpha -1}, \;\; 10 \le \he \le C \lep,
\end{equation*}
for some $\alpha \in (2/3,1)$, and $1\ll n \ll h_{ex}$. Assume that (\ref{hyp1})
holds. Then there exists an explicitly constructed vector field
$X_\ep$ such that, as $\ep \to 0$,
\begin{equation} \label{asser1}
G_\varepsilon(u_\varepsilon,A_\varepsilon)
  - f_\varepsilon(n) \ge n^2 I(\mu_*) + \frac{n^2}{36} \int_{\Omega}
  \abs{\frac{1}{n}\nabla_{A_\varepsilon'}
   u_\varepsilon - i u_\varepsilon X_\ep}^2
+o(n^2).
\end{equation}
Through estimates on $X_\ep$, it follows that for $\ep$ small
enough,
 $$\frac{1}{n} \|\nab_{A_\ep'} u_\ep\|_{\lti(\om)} \le C_0 +
C, $$ where $C$ is a constant depending only on $\om$
\end{thm}

In this result, we have inserted an extra term in the lower
bounds, measuring the $L^2$ distance between $\nab_{A'}u $ and
some known vector field, just as in Theorem \ref{energy_bound}
with the vector field $Y$. Just like in passing from Theorem
\ref{energy_bound} to Theorem \ref{norm_switch}, an estimate of
the $\lti$ norm of $X_\ep$ allowed us to deduce an upper bound for
$\|\nab_{A'} u\|_{\lti}$  by an order $n$, as desired.

\subsection{Application to convergence of the vorticity}
The estimates above have a direct application for vorticity
measures (Jacobians); indeed, $\|\nab_{A'} u\|_{\lti} \le C n$
implies that $\frac{1}{n}\curl(iu, \nab_{A'} u) $  and
$\frac{1}{n} \mu(u,A)$ are bounded as the derivative of an $\lti$
function.  More precisely, we need to use the Lorentz space $L^{2,1}$, whose dual space is $\lti$  (see Section \ref{cpt_jac} for definitions). Introducing the space $ \mathcal{X}(\Omega) = \{ f \in H^1_0
(\om)|\, \nab f \in L^{2,1}(\om)\}$, we obtain that
$\frac{\mu(u,A)}{n}$ is bounded, hence weakly-$*$ compact, in $\mathcal{X}^*(\Omega)$,
the dual of $\mathcal{X}(\Omega)$.

We deduce
\begin{thm}[see Proposition \ref{app_jac}]
Under the same assumptions as in Theorem \ref{3}, we have that 
\begin{equation*}
 \begin{alignedat}{2}
    & \frac{\mu(u_\varepsilon,A_\varepsilon)}{2\pi n} \wstar \delta_p   
	 && \text{ weakly-$*$ in }  \mathcal{X}^*(\Omega) \text{ and} \\
    & \frac{\tilde{\mu}(u_\varepsilon,A_\varepsilon)}{2 \pi n} \wstar  \mu_*    
	 &&\text{ locally weakly-$*$ in }  \mathcal{X}^*_{loc}(\Rn{2}),
 \end{alignedat}
\end{equation*}
where $\mu_*$ is the probability measure
given by Theorem \ref{ss}, and locally-weak-$*$ convergence in $\mathcal{X}^*_{loc}$ means weak-$*$ convergence in $\mathcal{X}^*(V)$ for every $V \csubset \Rn{2}$.
\end{thm}

This is a  slight improvement or alternative  to the known
compactness result in $(C^{0,\gamma})^*$. Since the space  $\mathcal{X}$
embeds into continuous functions, $\mathcal{X}^*$ is slightly larger
than measures, but neither $\mathcal{X}$ nor $C^{0,\gamma}$ embeds into the other.  Again in the regime 4, when $n$ is of the same order as $\he$, $\mu(u,A)/\he $ was shown to be compact in
$H^{-1}(\om)$, which is better than $\mathcal{X}^*$.

A second set of more precise results is obtained when one makes
the stronger assumption \begin{equation}\label{hyp2}
G_\ep(u_\ep,A_\ep) \le f_\ep(n) + n^2
I(\mu_*)+o(n^2),\end{equation} where $\mu_*$ is the weak limit of
$\frac{1}{2\pi n}\tilde{\mu}(u_\ep, A_\ep)$ as given by Theorem
\ref{ss}.
 Among them are a lower bound for $ \|\nab_{A'}u\|_{\lti(\om)}$ (see Theorem
 \ref{application_1}): as $\varepsilon \rightarrow 0$, 
$$C_0 \le \frac{1}{n} \|\nab_{A'}
u\|_{\lti(\om)} \le C_1,$$ where $0<C_0< C_1$ depend on $\Omega$, and the
extra convergence (see Corollary \ref{london_convergence}):
\begin{equation}
  \frac{1}{n}(iu,\nab_{A'}u) \wstar -\nabla^\bot G_p
   \text{ weakly-$*$} \text{ in }  L^{2,\infty}(\Omega) \text{ as } \varepsilon \rightarrow 0.
\end{equation}

Of course, the above results immediately apply to
energy-minimizing solutions in that regime. We also show that the
same results apply to the locally minimizing solutions found in
Chapter 11 of \cite{ss_book} (see our Theorem \ref{stable_solns}), and we get:
\begin{thm}\label{4}
Let $(u_\ep,A_\ep)$ be either global minimizers of the energy for $\he \le \lep$, or the locally minimizing solutions constructed in Theorem 11.1 of \cite{ss_book}, with the assumption that $h_{ex}$ is sufficiently large (see Theorem \ref{stable_solns} for the precise condition). Then we have that for $\varepsilon$ sufficiently small, 
$$C_0\le \frac{1}{n}\|\nab_{A_\ep'} u_\ep\|_{\lti(\Omega)}
 \le C_1 ,$$
where $C_0$ and $C_1$ are positive constants that depend only on
$\Omega$. 
\end{thm}
 Thus in both cases, $\|\nab_{A'} u\|_{\lti}$
can indeed serve as a normalizing factor to replace the
nonintrinsic $n$. Convergence results for $\mu(u,A')/n$ and $(iu,\nabla_{A'}u)/n$ such as the above are also stated.

\subsection{Strong convergence in Lorentz-Zygmund spaces}\label{oc}
The convergence of the vorticity measures  $\frac{1}{2\pi
n}\mu(u,A)$ and $\frac{1}{2\pi n} \tilde{\mu}(u,A)$  (as well as
those of the currents) to their limits is  weak-$*$ in $\mathcal{X}^*$. It
{\it does not hold strongly} in $\mathcal{X}^*$, nor is it true that $(iu,\nabla_{A'}u)/n
\to \np G_p$ strongly in $\lti$. The reason is (as pointed out
also in \cite{p1}) that the $\mathcal{X}^*$ norm acts a bit like the strong
norm on measures: for Dirac masses, we do not have
$\delta_{p_n}\to \delta_p$ strongly in $\mathcal{X}^*$ when the points
$p_n\to p$, but rather
$$2\sqrt{\pi} \le \|\delta_{p_n} - \delta_p\|_{X^*} \le 4\sqrt{\pi},$$ while we {\it do}
have $\delta_{p_n} \to \delta_p$ weakly-$*$ in $\mathcal{X}^*$. This explains
why similarly the strong convergences above do not hold in general.

One may then wonder if there is a weaker space (but still stronger
than $W^{-1,p}$ for $p<2$), in which strong convergence results
hold. We find that spaces based on the Lorentz-Zygmund spaces
$\lti\log^\gamma L(\om)$  with $\gamma<0$, which are just slightly
bigger than $\lti(\om)$, provide such a setting. We then obtain the strong
convergence analogues of the above results.  This is the object of
Section \ref{results_in_lz}.

\subsection{Plan}
The paper is organized as follows.
 In Section \ref{alg_tricks} we present the  ``completion of the square'' algebraic
 trick that serves as a general basis  for extracting new terms in energy lower
 bounds. We give a general statement
  for such lower bounds, which can be of independent interest, as well as applications
 in our setting.

  In Section
\ref{f_ep_apriori}, we present a first application of our results
of \cite{p1},  showing that under certain a priori energy upper
bounds on $F_\varepsilon$, the $L^{2,\infty}$ norm of $\nabla_{A}
u$ is comparable to the vorticity  mass $n$. However these a
priori bounds are rarely satisfied except for local minimizers for
low applied fields.

Section \ref{more_lower_bounds} refines the lower bounds of
 the energy $G_\varepsilon$ of \cite{ss_book} Chapter 9
to extract the terms used in the $L^{2,\infty}$ estimates.

Section \ref{G_ep_apriori} gives results similar to those of
Section \ref{f_ep_apriori} in the more useful case  of a priori
energy bounds on $G_\varepsilon$ that are satisfied by a large
class of minimizing solutions to the Ginzburg-Landau equations.

 Section \ref{convergence_results} establishes a compactness
result for the gauge-invariant Jacobians of configurations
satisfying the previously used  a priori energy upper bounds.  We
also establish $L^{2,\infty}$ weak-$*$ compactness results for the
gauge-invariant current.

 Section \ref{results_in_lz} improves the weak-$*$
compactness results to strong compactness in slightly larger
Lorentz-Zygmund spaces.

Section \ref{res_bnded} deals with the case of solutions with bounded vorticity.

Section \ref{solns_apps} applies all of these results to minimizing and locally minimizing 
solutions.

\vskip .5cm

{\it Acknowledgments}: This work was completed as a component of my Ph.D. dissertation at the Courant Institute.  I owe an enormous debt of gratitude to my advisor, Sylvia Serfaty, whose guidance and insight were integral to the writing of this paper; many thanks.  I would also like to extend my thanks to Etienne Sandier for suggesting the use of $\pqnorm{\nabla_A u}{2}{\infty}$ as a normalization for vorticity.

\section{Square completion and lower bounds }\label{alg_tricks}

In this section we review the main algebraic trick that is at
the core of our lower bounds, and we show how to generalize it to
obtain more lower bounds.

The main technical tool of \cite{p1}
 was the introduction of an auxiliary vector field, $Y$,
  defined on the collection of vortex balls, $\mathcal{B}$.
   The idea behind the introduction of this function is most
   easily seen when setting   $A=0$.  Write $u = \rho v$, with $\rho = \abs{u}$ and $v = e^{i\varphi}$.  Previous lower bounds on
$\int \abs{\nabla v}^2$ were found via the Cauchy-Schwarz inequality:
\begin{equation*}
\frac{1}{2}\int_{\partial B(a,r)}
 \abs{\nabla v}^2 = \frac{1}{2}\int_{\partial B(a,r)} \abs{\nabla \varphi}^2  \ge \frac{1}{4\pi r} \left(\int_{\partial B(a,r)} \nabla \varphi \cdot \tau  \right)^2  = \frac{4\pi^2 d^2}{4\pi r}
=  \pi \frac{d^2}{r},
\end{equation*}
where $d = \deg(u,\partial B(a,r))$.  Thus the inequality is sharp
if $u$ is radial and $\abs{\nabla v} \approx d/r$ on $\partial
B(a,r)$. We thus took the vector field $Y$  to be $\tau
d/r$, and rather than using Cauchy-Schwarz we ``completed the
square":
\begin{equation*}
\begin{split}
\frac{1}{2}\int_{\partial B(a,r)} \abs{\nabla v}^2 &= \frac{1}{2}\int_{\partial B(a,r)} \abs{\nabla v - \frac{d}{r} \tau}^2  + \frac{d}{r}\int_{\partial B(a,r)} \nabla \varphi \cdot \tau  - \frac{2\pi r d^2}{2r^2} \\
& =  \frac{1}{2}\int_{\partial B(a,r)} \abs{\nabla v - \frac{d}{r} \tau}^2 + \pi \frac{d^2}{r}.
\end{split}
\end{equation*}
This extracts a new term in the lower bound that measures the $L^2$ difference between $\nabla u$ and the optimal annular vortex profile, given by $Y$.  The implementation of this idea requires certain technical complications to handle the magnetic field and vorticity cancellation, but the main idea is as above: ``complete the square" with a function that $\nabla u$ ``should look like."

It is easy to extend this idea to domains. We begin with two
lemmas relying on the same algebraic manipulation.

\begin{lem}\label{j_square_completion}
Writing $j = (iu,\nabla_A u)$, we have that for any vector field $X:\Omega \rightarrow \Rn{2}$,
\begin{equation}
 \abs{\nabla_{A} u}^2
= \abs{\nabla_{A} u -i u X}^2 + 2 X \cdot j - \abs{X}^2\abs{u}^2.
\end{equation}
\end{lem}
\begin{proof}
 We calculate
 \begin{equation*}
 \begin{split}
  \abs{\nabla_{A} u}^2 &= \abs{\nabla_{A} u -iuX +iuX}^2 = \abs{\nabla_{A} u -i u X }^2 + 2 \Re((\nabla_{A} u -iuX)\cdot iuX) + \abs{uX}^2 \\
  & = \abs{\nabla_{A} u -i u X}^2 + 2 X \cdot \Re(\overline{iu} \nabla_{A} u) -2 \abs{uX}^2 + \abs{uX}^2 \\
  & = \abs{\nabla_{A} u -i u X }^2 + 2 X \cdot j - \abs{X}^2 \abs{u}^2.
 \end{split}
 \end{equation*}
\end{proof}

A simple modification of this lemma allows us to ignore the $\rho$ part of $u$ in the bound.

\begin{lem}\label{comp_sq_int_by_pts}
 Let $W \subseteq \Omega$ be a set on which $\abs{u} >0$, and hence on which $\nabla \varphi$ is well defined, where we write $u = \rho e^{i\varphi}$.  Let $H$ be a $C^1$ real-valued function on $W$.  Then
\begin{multline}
\frac{1}{2} \int_W \abs{\nabla \varphi -A}^2 + \abs{\curl{A}}^2
=\frac{1}{2} \int_W \abs{\nabla \varphi -A + \nabla^\bot H}^2 +
\frac{1}{2}\int_W \abs{\curl{A} -H}^2\\
- \frac{1}{2}\int_W \abs{\nabla H}^2 + \abs{H}^2   -
\int_{\partial W} H(\nabla \varphi -A)\cdot \tau
\end{multline}
where $\tau$ is counter-clockwise unit tangent vector field.
\end{lem}
\begin{proof}
Simple calculations show that $\abs{\nabla_A u}^2 = \abs{\nabla \rho}^2 + \rho^2 \abs{\nabla \varphi - A}^2$, $\abs{\nabla_A u - iuX}^2 = \abs{\nabla \rho}^2 + \rho^2 \abs{\nabla \varphi - A -X}^2$, and $j =(iu,\nabla_A u) = \rho^2(\nabla \varphi -A)$.  Use these equalities in Lemma \ref{j_square_completion} with $X = -\nabla^\bot H$, subtract $\abs{\nabla \rho}^2$ from both sides, and divide by $\rho^2$ to find the equality
\begin{equation}
\abs{\nabla \varphi -A + \nabla^\bot H}^2 - 2 \nabla^\bot H \cdot (\nabla \varphi -A)= \abs{\nabla \varphi -A}^2  + \abs{\nabla H}^2.
\end{equation}
Divide by $2$, integrate over $W$, and integrate the second term on the left by parts to get
\begin{multline}
 \frac{1}{2} \int_W \abs{\nabla \varphi -A + \nabla^\bot H}^2 -\int_W \curl{A}\cdot H - \int_{\partial W} H(\nabla \varphi -A)\cdot \tau \\= \frac{1}{2} \int_W \abs{\nabla \varphi -A}^2 +  \frac{1}{2}\int_W \abs{\nabla H}^2.
\end{multline}
The result follows by adding $\frac{1}{2}\int_W \abs{\curl{A}}^2 + \frac{1}{2} \int_W \abs{H}^2$ to both sides.

\end{proof}

These lemmas will be  put to crucial use in Section
\ref{more_lower_bounds}, where they are used to extract new terms
in the energy lower bounds in different parts of the exterior of
the vortex balls.

The identity obtained in Lemma \ref{comp_sq_int_by_pts} yields
convenient lower bounds when applied to well-chosen functions $H$.
More specifically, following the framework of \cite{bbh} Chapter
1, let $\{\omega_i\}$ be any finite family of disjoint closed ``holes"
with smooth boundary in $\om$ (for example balls), such that
$|u|>0$ in $\om\backslash \cup_i \omega_i$, with $d_i = \deg(u,
\partial \omega_i)$; we consider the function $H$ to be the solution to
\begin{equation}\label{external_h_def}
 \begin{cases}
  -\Delta H + H = 0  & \text{ in } \Omega \backslash \cup_i \omega_i \\
  H = c_i & \text{ on } \partial{\omega_i} \\
  H = 0 & \text{ on } \partial \Omega \\
  \int_{\partial \omega_i} \frac{\partial H}{\partial \nu} = 2\pi d_i .
 \end{cases}
\end{equation}
Here $c_i$ is an unknown constant, which is part of the problem,
 and $\nu$ is the outward pointing normal. The solution to this problem
is the minimizer of the variational problem
\begin{equation*}
 \inf_{Y} \frac{1}{2} \int_{\Omega \backslash \cup_i \omega_i} \abs{\nabla h}^2 + h^2 + 2\pi \sum_i d_i h\rvert_{\partial B_i},
\end{equation*}
where the space $Y$ is given by
\begin{equation*}
 Y = \{f\in H^1(\Omega \backslash \cup_i \omega_i) \;\vert\;   h\rvert_{\partial \omega_i} =\text{constant}, h\rvert_{\partial \Omega}=0   \}.
\end{equation*}
It should be noted that a function very similar to this one was
used in Chapter 1 of \cite{bbh} to obtain, through the same method,
lower bounds for $\si$-valued maps in punctured domains.

Such a function is useful in conjunction with  Lemma
\ref{comp_sq_int_by_pts} because it is constant on the boundary of
each $B_i$ and because  of the following simple identity, obtained
by integrating by parts and using \eqref{external_h_def}:
\begin{equation}\label{h_identity}
 \int_{\Omega \backslash \cup_i \omega_i } \abs{\nabla H}^2 + H^2 = \sum_i 2\pi d_i
 c_i.
\end{equation}
We thus obtain 

\begin{prop} \label{lbsimple}
Let $(u,A)$  be a  $C^1$ configuration defined on $ \om
\backslash \cup_i \omega_i$, where $\{\omega_i\}_i$ is a finite
collection of closed ``holes" with smooth boundaries in $\om$, with
$|u|>0$ in $\om \backslash \cup_i \omega_i$.  Let $v= u/|u|$,
$d_i = \deg(u,\partial \omega_i)$, and let $H$ be defined as in \eqref{external_h_def}. Then
\begin{multline}\label{res}
\hal \int_{\om\backslash \cup_i \omega_i} |\nab_A v|^2 + |\curl
A|^2  = \hal \int_{\om\backslash \cup_i \omega_i} |\nab H|^2 + H^2
\\+ \hal \int_{\om\backslash \cup_i \omega_i}
|\nab_{A}v +iv\np H|^2 + |\curl A - H|^2 - \sum_i c_i \int_{\omega_i} \curl A
.\end{multline}\end{prop}
\begin{proof}
We apply the result of Lemma \ref{comp_sq_int_by_pts} in $W =
\om\backslash \cup_i \omega_i $ with this $H$. Using the fact that
$H=c_i$ on each $\partial\omega_i$ and $0$ on $\bo$, and changing
the orientation to counterclockwise, we find
\begin{multline*}\hal\int_{\om\backslash \cup_i \omega_i} |\nab_A v|^2 +
|\curl A|^2 =\hal \int_{\om\backslash \cup_i \omega_i}  |\nab_{A}v +iv\np H|^2 + |\curl A - H|^2 \\ - \hal\int_{\om\backslash \cup_i
\omega_i} |\nab H|^2 + H^2
 + \sum_i c_i (2\pi d_i - \int_{\omega_i} \curl A)
\end{multline*}
Using \eqref{h_identity}, we conclude that \eqref{res} holds.
\end{proof}

Now this proposition provides, as in Chapter 1 of \cite{bbh}, lower
bounds on the energy in punctured domains by $\hal \int_{\om
\backslash \cup_i \omega_i} |\nab H|^2+H^2$, but in addition it
keeps track of the excess in the lower bound through the positive
term $\hal \int_{\om\backslash \cup_i \omega_i}  |\nab_{A}v +iv\np H|^2 + |\curl A - H|^2$ (the term $\sum_i c_i \int_{\omega_i}
\curl A$ can be shown to be small when the holes are small
enough). It then remains to bound from below $\hal \int|\nab H|^2
+H^2 $.

One application would be taking the holes $\omega_i$ to be the
smallest possible disjoint balls $B_i$ which cover the set where
$|u|< 1- \ep^{\alpha/4}$, such as the {\it  initial balls} in the
ball construction.  Then we obtain a lower bound on $F_\ep(u,A)$ by
$\hal\int_{\om\backslash \cup_i B_i} |\nab H|^2 + H^2.$  This term
can, in turn, be bounded below by the ball growth method (using equation 
\eqref{external_h_def} to estimate $\int\frac{\p H}{\p
\nu}$ on circles and easily readjusting the ball construction).
This would provide an alternate to Theorem \ref{energy_bound},
where this time the extra ``excess term" is
 $\int_{\omega \backslash \cup_i B_i}|\nab_{A}v +iv\np H
|^2 + |\curl A- H|^2$.
 This has the advantage that $H$ is well-described; for example
 $-\Delta H +H \approx  2\pi \sum_i d_i \delta_{a_i},$ where the $a_i$'s
 are the centers of the (small) initial balls. This can serve
 to control the difference between the Jacobian vorticity measure
  $\mu(u,A) $  and the quantity $2\pi \sum_i d_i \delta_{a_i}$, in $H^{-1} $ norm, by the energy-excess, as done by Jerrard-Spirn \cite{jspirn} in a different \medskip metric.

\section{The case of a priori upper bounds on $F_\varepsilon$}\label{f_ep_apriori}
We now focus on our initial question of obtaining upper and lower
bounds for $\|\nab_A u\|_{\lti(\om)}$, in terms of the number of
vortices. We start with a simple case where there is a strong
upper bound on the energy.  We use Theorem \ref{energy_bound} with final radius $r=1/2$
to produce a collection of balls, $\mathcal{B}$, and we let $n$ be the vorticity mass of these balls.  Again note that $n$ implicitly depends on $\varepsilon$, though we do not write the dependence explicitly. We also
heavily employ the convention that $C$ denotes a generic,
positive, universal constant that can change from line to line and
can stand for different constants even in the same expression.
When constants explicitly depend on other parameters it is noted.


We begin with a general argument that shows that if a
configuration has free energy $F_\varepsilon$ not too different
from $\pi n\abs{\log{\varepsilon}}$, then the $L^{2,\infty}$ norm
of the covariant derivative is of order $n$.  The next proposition
establishes both  upper and lower  bounds in terms of $n$.

\begin{prop}\label{degree_control}
 Suppose that $\{(u_\varepsilon,A_\varepsilon)\}$ are configurations satisfying the upper bound 
 $F_{\varepsilon}(\abs{u_\varepsilon}) \le \varepsilon^{\alpha-1}$ for some $\alpha \in (0,1)$.  The following hold.

 1.  Supposing  that $n\ge 1$ and
 \begin{equation}\label{deg_c_1}
  F_{\varepsilon}(u_\varepsilon,A_\varepsilon)  \le  \pi n \abs{\log{\varepsilon}}  + Mn^2,
 \end{equation}
we have that 
\begin{equation}
 \pqnormspace{\nabla_{A_\varepsilon} u_\varepsilon}{2}{\infty}{\Omega} \le C n,
\end{equation}
where $C$ depends only on $M$.

2. Supposing that $\pnorm{\nabla_{A_\varepsilon}
u_\varepsilon}{\infty} \le C/\varepsilon$ and that $n \ll
\abs{\log{\varepsilon}}$, we have that for $\varepsilon$ sufficiently
small,
 \begin{equation}
  \frac{\pi}{2} n \le \pqnormspace{\nabla_{A_\varepsilon}
u_\varepsilon}{2}{\infty}{\Omega}^2 +
\frac{1}{2\abs{\log{\varepsilon}}} \int_{\Omega}
(\curl{A_\varepsilon})^2  + o(1).
 \end{equation}
\end{prop}

\begin{proof}
 We neglect to write the subscript $\varepsilon$.  For the first assertion,
 Corollary 5.2 of \cite{p1}, applied with $r=1/2$, provides the bound
 \begin{equation}
  F_{\varepsilon}(u,A,\Omega) \ge C \pqnormspace{\nabla_A u}{2}{\infty}{\Omega}^2 + \pi n \left( \log\frac{1}{2\varepsilon n} -C\right) - \pi \sum d_i^2,
 \end{equation} for some universal constant $C$.
Noting that $\sum d_i^2 \le \left( \sum \abs{d_i} \right)^2 =n^2,$ we deduce
\begin{equation}
  F_{\varepsilon}(u,A,\Omega) - \pi n\abs{\log{\varepsilon}} \ge C \pqnorm{\nabla_A u}{2}{\infty}^2 - 3\pi n^2 -C n.
\end{equation}
Utilizing the upper bound of the hypothesis in conjunction with
this bound yields
\begin{equation}
  \pqnormspace{\nabla_A u}{2}{\infty}{\Omega}^2 \le  Cn^2,
\end{equation}
where $C$ depends only on $M$.

 For the second assertion, Proposition 1.4 of \cite{p1} applied to $|\nab_A u|$ gives us  $$\hal\io |\nab_A u|^2 \le \hal C^2 |\Omega| + \lep\|\nab_A u\|_{\lti(\om)}^2
 $$ where $|\nab_A u|\le C/\ep$.  Combining this with  Theorem \ref{energy_bound}
  applied with $r = 1/2$, we find
 \begin{equation}\label{l2inf_lb_c_0}
  \pi n \left( \log{\frac{1}{2n\varepsilon}} - C \right) \le \abs{\Omega} \frac{C^2}{2} + \abs{\log{\varepsilon}} \pqnormspace{\nabla_{A_\varepsilon} u_\varepsilon}{2}{\infty}{\Omega}^2 + \frac{1}{2} \int_{\Omega} (\curl{A_\varepsilon})^2.
 \end{equation}
The assumption $n \ll \abs{\log{\varepsilon}}$ proves that
\begin{equation}\label{l2inf_lb_c_1}
 \abs{\Omega} \frac{C^2}{2\abs{\log{\varepsilon}}} + C \pi \frac{n}{\abs{\log{\varepsilon}}} = o(1)
\text{ and that } \pi \frac{n\log{2n}}{\abs{\log{\varepsilon}}} \le \frac{\pi n}{2}
\end{equation}
for $\varepsilon$ small enough.  Inserting \eqref{l2inf_lb_c_1} into \eqref{l2inf_lb_c_0} yields the result.
\end{proof}

Unfortunately, in practice the above assumptions on the energy are
really only useful in the case of $n$ and $h_{ex}$ bounded
independently of $\varepsilon$.  Moreover, the upper and lower
bounds do not quite match, with the lower bound being of
order$\sqrt{n}$ and the upper bound of order $n$.  However, a
little extra work in what follows allows us to use an a priori
upper bound on the full energy $G_\varepsilon$ in conjunction with
the assumption that $1 \ll n \ll h_{ex}$ to prove that
$\pqnormspace{\nabla_A u}{2}{\infty}{\Omega}$ is bounded above and
below by terms of order $n$.  This improvement is accomplished by
examining the energy contained in a large annulus but outside the
balls produced by the ball construction.  This strategy follows
that employed in Chapter 9 of \cite{ss_book}.


\section{Improving the lower bounds}\label{more_lower_bounds}

\subsection{Definitions and notation}\label{defs_and_notation}
We are now in the setting of \cite{ss_book} Chapter 9, for regime 3 in the introduction.  The goal of this section is to prove an improved version of Theorem \ref{3}. The proof follows all the steps of Chapter 9 of \cite{ss_book}, adding an extra term in each lower bound via a square completion trick.  We recall  that we assume that $\xi_0$, defined by \eqref{xi0}, achieves its minimum at a single point $p$. This is satisfied, for instance, if we assume that $\Omega$ is convex. Indeed, if
$\Omega$ is convex, then the sub-level sets $\{ \xi_0 \le t
\}$ are convex (see \cite{caf_fried}); this, combined with the
fact that the set where $\xi_0$ achieves its minimum is finite,
proves that there is exactly one minimum point.  We shall further assume that $D^2 \xi_0(p)$ is positive definite.  We write
$\underline{\xi_0}=\xi_0(p)$ and define the constant $J_0$, which
depends only on the domain $\Omega$, by $J_0 = \frac{1}{2}
\norm{\xi_0}^2_{H^1(\Omega)}$.

As in \cite{ss_book} Chapter 9, we consider two sizes of balls: small and large.
We initially construct, via Theorem \ref{energy_bound}, a
collection of small balls $\mathcal{B}'$ such that $r':=
r(\mathcal{B}') = C\varepsilon^{\alpha/2}$.  Write $n' = \sum_i
\abs{d_i'}$ for the vorticity mass of the small balls.  We assume that the inequality $1/\sqrt{h_{ex}} > 2r'$ holds; below we will state an assumption on the size of $h_{ex}$ sufficient to give this property.  An application of the ball growth
lemma then allows us to grow $\mathcal{B}'$ into a collection of
large balls, $\mathcal{B}$, such that $r:= r(\mathcal{B}) =
1/\sqrt{h_{ex}}$.  Write $n= \sum_i \abs{d_i}$ for the vorticity
mass of the large balls.

In the remainder of the paper, we will work under the following set of hypotheses, borrowed from
\cite{ss_book}, Chapter 9: $\{(u_\varepsilon,A_\varepsilon)\}$ are
configurations which satisfy
\begin{enumerate}
\item[(H1)]\begin{equation}\label{H1}
 F_\varepsilon(u_\varepsilon,A_\varepsilon') \le \varepsilon^{\alpha -1}, \;\; 10 \le  h_{ex} \le C\varepsilon^{-\beta},
\end{equation}
for some $\alpha \in (2/3,1)$ and $\beta \in (0,3\alpha/2 -1)$.
\item[(H2)] Letting $n$ denote $\sum_i |d_i|$, the sum of the degrees of
the balls of total final radius $r= 1/\sqrt{\he}$, we have $1 \le n \ll h_{ex}$. 
\item[(H3)] One of the following holds:
\begin{equation}\label{H3}
 h_{ex} \le C \abs{\log{\varepsilon}} \;\; \text{ or }\;\; n' = n.
\end{equation}
\item[(H4)] The upper bound
\begin{equation}\label{H4}
G_\varepsilon(u_\varepsilon,A_\varepsilon) \le f_\varepsilon(n) + C_0n^2
\end{equation}
holds for some constant $C_0\ge 0$, where $f_\ep$ is defined by
 \begin{equation}\label{fep_def}
  f_\varepsilon(n) := h_{ex}^2 J_0 + \pi n \abs{\log{\varepsilon}} + 2 \pi n h_{ex} \underline{\xi_0} + \pi n^2 S_\Omega(p,p) + \pi(n^2-n)\log{\frac{1}{\ell}}.
\end{equation}
Here we have written $S_\Omega(\cdot,p)$ for the function defined by 
\begin{equation}
S_\Omega(x,p) = G_p(x) + \log{\abs{x-p}},
\end{equation}
where $G_p$ is defined by \eqref{Gp}.
\end{enumerate}

 In most of what follows we consider the case $1 \ll n \ll h_{ex}$ (regime 3), but we will
always explicitly state the assumption $1 \ll n$ in the hypotheses of the
results when it is needed.  As mentioned above, when $n \ll h_{ex}$,
vortices tend to form near the point $p$, and the typical
inter-vortex distance, and by extension, the typical distance
between a vortex and the point $p$, is of the order $\ell = \sqrt{n/h_{ex}}$ (see Section 9.1.1 of \cite{ss_book} for a more thorough discussion).  When $n/h_{ex}$ is not small, the vortices are dispersed throughout the domain and our method fails to capture the lower
bound in terms of $n$, as seen in Proposition \ref{degree_control}.

Besides repeated application of the square completion trick, the primary technical tool of this section, borrowed from Chapter 9 of \cite{ss_book}, is the introduction of the annulus $B(p,\delta)\backslash B(p,K\ell)$, where $K$ and $\delta$ are constants independent of $\varepsilon$ that will eventually be sent to $\infty$ and $0$ respectively.  For $t \in
(K\ell,\delta)$ we define the degree function
\begin{equation}\label{app_d_def}
 D(t) = \sum_{\abs{b_i-p} \le t} d_i,
\end{equation}
where the $\{b_i\}$ are points in the balls $\{B_i\} =
\mathcal{B}$ chosen later in Proposition \ref{app_balls}. Note
that $|D(t)|\le n$.  Finally, since there can be some vortices
(i.e. balls $B \in \mathcal{B}$) contained in the annulus
$B(p,\delta)\backslash B(p,K\ell)$, we must track their location by
defining the set
\begin{equation}\label{app_t_def}
T = \{ t \in (K\ell,\delta) \;|\; \partial B(p,t) \cap \mathcal{B} \neq
\varnothing \}.
\end{equation}
Note that for $t \notin T$,
\begin{equation*}
 D(t) = \deg{(u,\partial B(p,t))},
\end{equation*}
and that $\abs{T} \le 2r = 2/\sqrt{h_{ex}}$, where $\abs{T}$ denotes the measure of
$T$.



\subsection{Lower bounds in the balls and energy splitting}\label{lb_split}

Here we adapt the results of \cite{p1} to deal with the full
 energy $G_\varepsilon$.  This entails showing how energy
lower bounds hold on the two-phase ball construction (the balls in
$\mathcal{B}'$ and $\mathcal{B}$) and also proving an ``energy
splitting lemma" that allows us to pass from the full energy
$G_\varepsilon$ to a sum of the free energy $F_\varepsilon$ and
other terms.

Our first lemma provides lower bounds on the free energy in the
balls.  It is a modification of Lemma 9.1 of \cite{ss_book} that
incorporates a term involving the vector field $Y_\varepsilon$ of Theorem \ref{energy_bound} into lower bounds
in $\mathcal{B}$. This is different from the result of Theorem \ref{energy_bound} only in that the estimates are constructed in two stages:
first in $\mathcal{B}'$ and then in $\mathcal{B} \backslash
\mathcal{B}'$. In all that follows, we  abuse notation by writing
$\mathcal{B}$ in place of $\cup_{B \in \mathcal{B}} B$.


\begin{lem}{(\cite{ss_book} Lemma 9.1 Redux)}\label{redux_92}
Suppose that configurations $\{(u_\varepsilon,A_\varepsilon)\}$
satisfy assumption (H1).  Let $Y_\varepsilon$ be the vector field of Theorem \ref{energy_bound}, applied to the large balls $\mathcal{B}$ with $r= 1/\sqrt{h_{ex}}$.  Then
 \begin{multline*}
  \frac{1}{2} \int_{\mathcal{B}} \abs{\nabla_{A_\varepsilon'} u_\varepsilon}^2 + \frac{1}{2\varepsilon^2} (1- \abs{u_\varepsilon}^2)^2 + r^2(\curl{A_\varepsilon'})^2 \\ \ge \frac{1}{36} \int_{\mathcal{B}} \abs{\nabla_{A_\varepsilon'} u_\varepsilon -iu_\varepsilon Y_\varepsilon }^2
  + \pi \left( n\log{\frac{r}{n \varepsilon}} + \frac{\alpha}{4}(n' -n)\log{\frac{1}{\varepsilon}} \right) -Cn
 \end{multline*}
for $\varepsilon$ sufficiently small.
\end{lem}

\begin{proof}
Theorem \ref{energy_bound} provides the bound 
 \begin{equation*}
  \frac{1}{2} \int_{\mathcal{B}} \abs{\nabla_{A_\varepsilon'} u_\varepsilon}^2 + \frac{1}{2\varepsilon^2} (1- \abs{u_\varepsilon}^2)^2 + r^2(\curl{A_\varepsilon'})^2  \ge \frac{1}{18} \int_{\mathcal{B}} \abs{\nabla_{A_\varepsilon'} u_\varepsilon -iu_\varepsilon Y_\varepsilon }^2
  + \pi \left( n\log{\frac{r}{n \varepsilon}} -C\right).
 \end{equation*}
On the other hand, Lemma 9.1 of \cite{ss_book} provides the bound 
\begin{equation*}
  \frac{1}{2} \int_{\mathcal{B}} \abs{\nabla_{A_\varepsilon'} u_\varepsilon}^2 + \frac{1}{2\varepsilon^2} (1- \abs{u_\varepsilon}^2)^2 + r^2(\curl{A_\varepsilon'})^2  \ge 
  \pi \left( n\log{\frac{r}{n \varepsilon}} + \frac{\alpha}{2}(n' -n)\log{\frac{1}{\varepsilon}} \right) -Cn.
\end{equation*}
The result follows by averaging these two bounds.
\end{proof}

With this lemma in hand, we can prove the following Proposition, a variant of Proposition 9.3 from \cite{ss_book}.  It shows how the full energy, $G_\varepsilon$, can be split and bounded below by the free energy and various other terms.

\begin{prop}{(\cite{ss_book} Proposition 9.3 Redux)}\label{app_balls}
Suppose configurations $\{(u_\varepsilon,A_\varepsilon)\}$ satisfy the assumption (H1). 
Then there exist points $b_i \in
B_i$ such that, letting $\nu = \sum_i d_i \delta_{b_i}$, the
following estimates hold for $\varepsilon$ sufficiently small.
\begin{multline}\label{app_balls_0}
 G_{\varepsilon}(u_{\varepsilon},A_{\varepsilon}) \ge h_{ex}^2 J_0 + 2\pi h_{ex} \int \xi_0 d\nu + F_{\varepsilon}(u_{\varepsilon},A_{\varepsilon}') \\ - C(n'-n)r h_{ex} -C h_{ex} \varepsilon^{3\alpha / 2 -1} - C h_{ex}^2 \varepsilon^{\alpha}
\end{multline}
\begin{multline}\label{app_balls_1}
 F_{\varepsilon}(u_\varepsilon,A_\varepsilon') \ge \pi n \log{\frac{r}{n\varepsilon}} + F_{\varepsilon}(u_\varepsilon,A_\varepsilon',\Omega \backslash \mathcal{B}) + \frac{1-r^2}{2} \int_{\mathcal{B}}(\curl{A_\varepsilon'})^2 \\
 +\frac{1}{36} \int_{\mathcal{B}} \abs{\nabla_{A_\varepsilon'} u_\varepsilon -iu_\varepsilon Y_\varepsilon }^2 + \pi \frac{\alpha}{4}(n'-n)\abs{\log{\varepsilon}} -Cn.
\end{multline}
Here the vector field $Y_\varepsilon$ is the one from Theorem \ref{energy_bound}, and  $C$ is a
universal constant.
\end{prop}
\begin{proof}
 The proof is the same as the proof of Proposition 9.3 of \cite{ss_book}, except that we use our Lemma \ref{redux_92} in place of their Lemma 9.1 in order to recover the $Y_\varepsilon$ difference term.
\end{proof}

\subsection{Lower bounds in the annulus $B(p,\delta)\backslash B(p,K\ell)$}

We now show how to bound the energy contained in the annulus
around the point $p$ using a vector field $Y_\varepsilon$ similar to, but simpler,
than the one defined in the balls.  Denote the annulus by
$\mathcal{A} = B(p,\delta) \backslash B(p,K\ell)$ and note that $n \ll
h_{ex}$ implies that $K\ell \rightarrow 0$ as $\varepsilon
\rightarrow 0$, while $\delta$ stays fixed.  The only
difficulty in defining $Y_\varepsilon$ in $\mathcal{A}$ is that the annulus
can contain balls from $\mathcal{B}$ sprinkled throughout.  We get
around this by taking $Y_\varepsilon$ to vanish there (ultimately we view this
$Y_\varepsilon$ as extending the $Y_\varepsilon$ already defined in the balls).  Indeed,
define $Y_\varepsilon: \mathcal{A} \rightarrow \Rn{2}$ by
 \begin{equation}\label{ann_G_def}
  Y_\varepsilon(x) = \begin{cases}
          0, & \abs{x-p} \in T \\
      D(\abs{x-p}) \tau_{\partial B(p,\abs{x-p})}(x) \frac{1}{\abs{x-p}}, & \abs{x-p} \in (K\ell,\delta) \backslash T.
         \end{cases}
 \end{equation}
Here, as before, $\tau_{\partial B}$ is the counter-clockwise unit tangent vector field
to the boundary of a ball, $\partial B$, and $D(t)$ and $T$ are
defined by \eqref{app_d_def} and \eqref{app_t_def}.  Note that $Y_\varepsilon$ also depends on 
$\delta$ and $K$, but we do not write that in the notation.

Following Lemma 9.3 of \cite{ss_book}, we estimate from below the
energy contained in $\mathcal{A}$.

\begin{lem}{(\cite{ss_book} Lemma 9.3 Redux)}\label{93redux}
Suppose configurations $\{(u_\varepsilon,A_\varepsilon)\}$ satisfy assumption (H1).  Let $Y_\varepsilon$ be defined on $\mathcal{A}$ by \eqref{ann_G_def}.  Then for $\varepsilon$ sufficiently small,
\begin{multline}
 \frac{1}{2} \int_{\mathcal{A} \backslash \mathcal{B}} \abs{\nabla_{A_\varepsilon'} u_\varepsilon}^2 + \frac{1}{4} \int_{\mathcal{A}} (\curl{A_\varepsilon'})^2  \ge \frac{1}{36} \int_{\mathcal{A} \backslash \mathcal{B}} \abs{\nabla_{A_\varepsilon'} u_\varepsilon -iu_\varepsilon Y_\varepsilon }^2  \\
 +  \pi  \int_{K\ell}^{\delta} D^2(t)  \frac{dt}{t}   dt
 -\pi n^2 \delta^2 - 2\pi \frac{n^{3/2}}{K} - \pi n^2 \varepsilon^{\alpha/4} \log{\frac{\delta}{K\ell}}.
\end{multline}
\end{lem}
\begin{proof}
 We suppress the subscript $\varepsilon$ on $u_\varepsilon$, $Y_\varepsilon$
  and $A_\varepsilon'$.  The proof proceeds as in \cite{ss_book}
  except we use Lemma 3.2 of \cite{p1} with $\lambda = \frac{1}{2\delta}$ to recover the $Y$ term.  Indeed, for $t \notin T$ it yields
\begin{multline}
\frac{1}{2} \int_{\partial B(p,t)} \abs{\nabla_{A'} v}^2 + \frac{1}{4\delta} \int_{B(p,t)} (\curl{A'})^2
\ge \frac{1}{2} \int_{\partial B(p,t)} \abs{\nabla_{A'} v -ivY}^2 +
\pi D^2(t) \left(\frac{1}{t} - \delta \right).
\end{multline}
Then, using the fact  that $Y=0$ in $\mathcal{A} \cap \{ |x-p|\in T\}$, we have
\begin{equation}\label{93r_1}
\begin{split}
  & \frac{1}{2} \int_{\mathcal{A} \backslash \mathcal{B}} \abs{\nabla_{A'} v}^2 + \frac{1}{4} \int_{\mathcal{A}} (\curl{A'})^2 \ge \frac{1}{2} \int_{(K\ell,\delta) \backslash T} \int_{\partial B(p,t)} \abs{\nabla_{A'} v -ivY}^2 dt\\
  & + \frac{1}{2} \int_{(\mathcal{A} \backslash \mathcal{B})
  \cap \{ \abs{x-p} \in T\}} \abs{\nabla_{A'} v}^2 +
  \int_{(K\ell,\delta)\backslash T} \pi D^2(t) \left(\frac{1}{t} - \delta \right)dt \\
  & = \frac{1}{2} \int_{\mathcal{A} \backslash \mathcal{B}} \abs{\nabla_{A'} v -ivY}^2 + \int_{K\ell}^{\delta} \pi D^2(t) \left(\frac{1}{t} - \delta \right)dt
   - \int_T \pi D^2(t) \left( \frac{1}{t} -\delta \right).
\end{split}
\end{equation}
Now we bound
\begin{equation}
 \int_{K\ell}^\delta \pi D^2(t) \delta dt \le \pi n^2 \delta^2.
\end{equation}
The fact that $T\subset (K\ell,\delta)$ and $\abs{T} \le 2r$ implies that
\begin{equation}\label{93r_2}
 \int_T \pi D^2(t) \frac{dt}{t} \le \pi n^2 \int_{K\ell}^{K\ell + 2r}  \frac{dt}{t} = \pi n^2 \log \left(1+\frac{2r}{K\ell} \right) \le \pi n^2 \frac{2r}{K\ell} = 2\pi \frac{n^{3/2}}{K}.
\end{equation}
Then \eqref{93r_1} -- \eqref{93r_2} provide the bound
\begin{equation}\label{93r_3}
   \frac{1}{2} \int_{\mathcal{A} \backslash \mathcal{B}} \abs{\nabla_{A'} v}^2 + \frac{1}{4} \int_{\mathcal{A}} (\curl{A'})^2 \ge \frac{1}{2} \int_{\mathcal{A} \backslash \mathcal{B}} \abs{\nabla_{A'} v -ivY}^2  + \int_{K\ell}^{\delta} \pi D^2(t) \frac{dt}{t} - \pi n^2\delta^2  - 2\pi \frac{n^{3/2}}{K}.
\end{equation}

We now recall that $\abs{\nabla_{A'} u}^2 = \abs{\nabla \abs{u}}^2 +
\abs{u}^2 \abs{\nabla_{A'} v}^2$ and that $1-\varepsilon^{\alpha/4}
\le \abs{u} \le 1 + \varepsilon^{\alpha/4}$ on $\Omega \backslash
\mathcal{B}$.  This implies that $\abs{\nabla_{A'} u}^2 \ge  (1-2\varepsilon^{\alpha/4})\abs{\nabla_{A'} v}^2$.
To conclude, we multiply both sides of \eqref{93r_3} by $1-2\varepsilon^{\alpha/4}$ and use the fact that $1 \ge \frac{\abs{u}^2(1-2\varepsilon^{\alpha/4})}{(1+\varepsilon^{\alpha/4})^2} \ge
\frac{\abs{u}^2}{36}$ for $\varepsilon$ sufficiently small.
 The result then follows by noting that $
 \int_{K\ell}^{\delta} \pi D^2(t) \frac{dt}{t} \le \pi n^2 \log{\frac{\delta}{K\ell}}.
$

\end{proof}

 With this modification established we deduce a corresponding
 modification of Proposition 9.4 from \cite{ss_book}
 (the proof is exactly as in \cite{ss_book}):

\begin{prop}{(\cite{ss_book} Proposition 9.4 Redux)}\label{app_annuli}
Suppose configurations $\{(u_\varepsilon,A_\varepsilon)\}$ satisfy
assumption (H1).  Then there exist positive $K_0,
\delta_0$ such that if $K \ge K_0, \delta \le \delta_0$ and if $\ell$
and $\varepsilon$ are sufficiently small, letting $\nu = \sum_i
d_i \delta_{b_i}$, we have the estimate
\begin{multline}
  \frac{1}{2} \int_{\mathcal{A} \backslash \mathcal{B}} \abs{\nabla_{A_\varepsilon'} u_\varepsilon}^2 + \frac{1}{4} \int_{\mathcal{A}} (\curl{A_\varepsilon'})^2  + 2 \pi h_{ex} \int \xi_0 d\nu \ge \frac{1}{36} \int_{\mathcal{A} \backslash \mathcal{B}}
  \abs{\nabla_{A_\varepsilon'} u_\varepsilon -u_\varepsilon Y_\varepsilon }^2 \\
  + \pi n^2 \log{ \frac{\delta}{K\ell} } + 2 \pi n h_{ex} \underline{\xi_0}
  + 2 \pi h_{ex} \sum_{\substack{b_i \in B(p,K\ell) \\ d_i >0}} d_i(\xi_0(b_i) - \underline{\xi_0})
  - \pi n^2 \delta^2 + o(n^2).
  \end{multline}
Moreover, for any $t \in [K\ell,\delta]$ we have that
\begin{equation}\label{app_annuli_0}
 \abs{\frac{D(t) - n}{n}} \le C\ell^2 \left( \frac{1}{t^2} + 1\right).
\end{equation}

\end{prop}

\subsection{Lower bounds outside $B(p,\delta) \cup \mathcal{B}$}
\label{outside_energy}
In this section we find lower bounds in the region outside
the ball $B(p,\delta)$ and the collection of balls, $\mathcal{B}$.
These bounds are different from those found in \cite{ss_book} in
that we again use a completion of the square trick to find a novel
term in the lower bounds.  In this region, however, we use the
more natural function $G_p$ (see \eqref{Gp}) and its perpendicular gradient $\nabla^\bot G_p$ rather than the ad hoc $Y_\varepsilon$ vector fields used in previous sections.

We now state a result, which is part of Proposition 9.5 of \cite{ss_book}, that
provides information on the weak limits of $j' = (iu,\nabla_{A'} u)$,
$h' = \curl{A'}$, and $\mu' = \curl(j' +A')$ when suitably
normalized. This information on the weak limits will be used to
deal with the cross terms that arise when we ``complete the
square'' with Lemma \ref{j_square_completion}. In what follows it
will be useful to work with the function $f:\Rn{+} \rightarrow
\Rn{+}$ given by $f(x) = \vchi_{[0,1]}(x) +
x^{-1}\vchi_{[1,\infty)}(x)$.

\begin{lem}\label{leftover_lower_bound_pt1}
  Suppose configurations $\{(u_\varepsilon,A_\varepsilon)\}$
  satisfy (H1) -- (H4).
Then up to extraction as $\varepsilon \rightarrow 0$,
\begin{eqnarray*}
&  \frac{1}{n}  f(\abs{u_\varepsilon}) \vchi_{\Omega \backslash
\mathcal{B}} j_\varepsilon'  \rightharpoonup j_* & \text{weakly in } \ L^2_{loc}(\Omega \backslash \{p\})\\
&  \frac{1}{n} h_\varepsilon' \rightharpoonup h_* & \text{weakly
in } \ L^2(\Omega)\\ &
 \frac{1}{n}\mu_\varepsilon' \wstar 2\pi \delta_p & \text{weakly-$*$ in } \
 (C_c^{0,\gamma}(\Omega))^*  \ \text{for some }  \
\gamma \in (0,1).\end{eqnarray*}  The limits satisfy the relation 
 $\curl{j_*} + h_* = 2\pi \delta_p.$
Moreover, as $\delta \rightarrow 0,$
\begin{equation}\label{l_l_b_p1_1}
\int_{\Omega \backslash B(p,\delta)} G_p h_* - \nabla^\bot G_p \cdot j_* = \int_{\Omega \backslash B(p,\delta)} \abs{\nabla G_p}^2 +\abs{G_p}^2 + o_\delta(1).
\end{equation}
\end{lem}
\begin{proof}
 Everything except \eqref{l_l_b_p1_1} is proved directly in Proposition 9.5 of \cite{ss_book}.  For \eqref{l_l_b_p1_1} we make a minor modification of their argument.  They show in the proof, \cite[(9.80) -- (9.82)]{ss_book}, that
 \begin{equation}
  \int_{\Omega \backslash B(p,\delta)} G_p (h_* - G_p) - \nabla^\bot G_p \cdot (j_* + \nabla^\bot G_p) = o_\delta(1).
 \end{equation}
To conclude then, we simply write
\begin{equation}
 G_p h_* - \nabla^\bot G_p \cdot j_* = \abs{G_p}^2 + \abs{\nabla^\bot G_p}^2 + G_p (h_* - G_p) - \nabla^\bot G_p \cdot (j_* + \nabla^\bot G_p)
\end{equation}
and integrate.

\end{proof}

\begin{remark}
 This lemma guarantees that, up to extraction, the sequences $j_\varepsilon'/n$, $h_\varepsilon'/n$, and $\mu_\varepsilon'/n$ have weak limits.  Henceforth we assume that the extraction has already been performed so that these weak limits exist. 
\end{remark}

The following proposition provides the energy lower bounds in $\Omega \backslash (B(p,\delta) \cup
\mathcal{B})$.  The completion of the square trick is done with $-i u_\varepsilon n f(\abs{u_\varepsilon})\nabla^\bot G_p$.

\begin{prop}\label{exterior_est_w_g}
  Suppose configurations $\{(u_\varepsilon,A_\varepsilon)\}$
satisfy (H1) -- (H4). Define the set $W = \Omega \backslash (B(p,\delta) \cup
\mathcal{B})$.  Then for $\varepsilon$ sufficiently
small,
\begin{multline}\label{e_e_w_g_0}
\frac{1}{2} \int_W \abs{\nabla_{A_\varepsilon'} u_\varepsilon }^2 + \abs{\curl{A_\varepsilon'}}^2 \ge  \frac{1}{2} \int_W \abs{\nabla_{A_\varepsilon'} u_\varepsilon + iu_\varepsilon nf(\abs{u_\varepsilon}) \nabla^\bot G_p}^2 + \frac{1}{2}\int_W \abs{\curl{A_\varepsilon'}  - n G_p}^2 \\- \pi n^2 \log{\delta} +\pi n^2 S_\Omega(p,p) +o_\delta(n^2) +o(n^2).
\end{multline}
\end{prop}

\begin{proof}
 Suppress the subscript $\varepsilon$ for convenience.  To begin, we use $X = -n f(\abs{u})\nabla^\bot G_p$ in Lemma \ref{j_square_completion} to get the identity
\begin{equation}\label{e_e_w_g_1}
 \abs{\nabla_{A'} u}^2 = \abs{\nabla_{A'} u + i n u f(\abs{u}) \nabla^\bot G_p}^2  -2 n f(\abs{u}) \nabla^\bot G_p \cdot j' - n^2\abs{\nabla G_p}^2 f^2(\abs{u}) \abs{u}^2.
\end{equation}
We complement this with a ``completion of the square'' for $\curl{A'}$:
\begin{equation}\label{e_e_w_g_2}
 \abs{\curl{A'}}^2 = \abs{\curl{A'} - n G_p}^2 - n^2 \abs{G_p}^2 + 2 n G_p \curl{A'}.
\end{equation}
We add \eqref{e_e_w_g_1} and \eqref{e_e_w_g_2}, divide by $2$, and integrate over $W$ to arrive at the identity
\begin{multline}\label{e_e_w_g_3}
 \frac{1}{2}\int_W \abs{\nabla_{A'} u}^2 + \abs{\curl{A'}}^2 = \frac{1}{2} \int_W \abs{\nabla_{A'} u + i n u f(\abs{u}) \nabla^\bot G_p}^2 + \abs{\curl{A'} - n G_p}^2 \\ + n \int_W  G_p \curl{A'} - \nabla^\bot G_p \cdot j' f(\abs{u})  - \frac{n^2}{2} \int_W \abs{G_p}^2 + \abs{\nabla G_p }^2 f^2(\abs{u}) \abs{u}^2.
\end{multline}
We want to keep the first integral on the right, but we keep continue working with the second and third integrals.

The function $f$ satisfies the inequality $xf(x) \le 1$ for all $x\ge 0$, and hence $\abs{u}^2 f^2(\abs{u}) \le 1$.  This, when combined with the fact that $G_p \in H^1_{loc}(\Omega \backslash \{p\})$, provides an estimate for the third integral on the right side of \eqref{e_e_w_g_3}.  Indeed,
\begin{equation}\label{e_e_w_g_4}
\begin{split}
  -\frac{n^2}{2} \int_W \abs{G_p}^2 + \abs{\nabla G_p }^2 f^2(\abs{u}) \abs{u}^2 &\ge  -\frac{n^2}{2} \int_W \abs{G_p}^2 + \abs{\nabla G_p }^2 \\
&= -\frac{n^2}{2} \int_{\Omega \backslash B(p,\delta)} \abs{G_p}^2 + \abs{\nabla G_p }^2 + o(n^2),
\end{split}
\end{equation} as $\ep \to 0$,
where we have used the fact that $\abs{\mathcal{B}} = o(1)$ and
 $G_p \in H^1_{loc}(\Omega \backslash \{p\})$.

To estimate the second integral in \eqref{e_e_w_g_3} we  write
\begin{equation}\label{e_e_w_g_5}
\begin{split}
n \int_W  G_p \curl{A'} - \nabla^\bot G_p \cdot j' f(\abs{u}) &= n^2 \int_W  G_p \frac{1}{n}\curl{A'} - \nabla^\bot G_p \cdot \frac{1}{n}j' f(\abs{u}) \\
& = n^2 \int_{\Omega \backslash B(p,\delta)}  G_p \vchi_{\Omega \backslash \mathcal{B}} \frac{1}{n}\curl{A'} - \nabla^\bot G_p \cdot \frac{1}{n}j' f(\abs{u})\vchi_{\Omega \backslash \mathcal{B}}.
\end{split}
\end{equation}
Now the weak $L^2(\Omega)$ convergence $\curl{A'}/n
\rightharpoonup h_*$ of Lemma \ref{leftover_lower_bound_pt1}
implies that
\begin{equation}\label{e_e_w_g_6}
\int_{\Omega \backslash B(p,\delta)}  G_p \vchi_{\Omega \backslash \mathcal{B}} \frac{1}{n}\curl{A'} = o(1) + \int_{\Omega \backslash B(p,\delta)}  G_p h_*,
\end{equation}
while the weak $L^2_{loc}(\Omega \backslash \{p\})$ convergence $ \frac{1}{n}j' f(\abs{u})\vchi_{\Omega \backslash \mathcal{B}} \rightharpoonup j_*$ implies that
\begin{equation}\label{e_e_w_g_7}
\int_{\Omega \backslash B(p,\delta)}  \nabla^\bot G_p \cdot \frac{1}{n}j' f(\abs{u})\vchi_{\Omega \backslash \mathcal{B}} = o(1) + \int_{\Omega \backslash B(p,\delta)}  \nabla^\bot G_p \cdot j_*.
\end{equation}
Combining \eqref{e_e_w_g_5} -- \eqref{e_e_w_g_7} and applying \eqref{l_l_b_p1_1} allows us to conclude that
\begin{equation}\label{e_e_w_g_8}
n \int_W  G_p \curl{A'} - \nabla^\bot G_p \cdot j' f(\abs{u}) = n^2 \int_{\Omega \backslash B(p,\delta)} \abs{\nabla G_p}^2 +\abs{G_p}^2 + o_\delta(n^2) + o(n^2).
\end{equation}

From these calculations we see that the second and third integrals on the right of \eqref{e_e_w_g_3} can be written as an integral of $\abs{G_p}^2 +\abs{\nabla G_p}^2$ plus an error term.  Indeed, summing \eqref{e_e_w_g_4} and \eqref{e_e_w_g_8} yields the bound
\begin{multline}\label{e_e_w_g_9}
n \int_W  G_p \curl{A'} - \nabla^\bot G_p \cdot j' f(\abs{u})  - \frac{n^2}{2} \int_W \abs{G_p}^2 + \abs{\nabla G_p }^2 f^2(\abs{u}) \abs{u}^2 \\
\ge \frac{n^2}{2} \int_{\Omega \backslash B(p,\delta)} \abs{\nabla G_p}^2 +\abs{G_p}^2 + o_\delta(n^2) + o(n^2).
\end{multline}
To complete the proof we must estimate this integral of $G_p$ and its gradient.  For this we use the expansion $G_p(x) = -\log{\abs{x-p}} + S_\Omega(x,p)$.  Since $G_p$ vanishes on $\partial \Omega$, we may compute
\begin{equation}\label{e_e_w_g_10}
\begin{split}
 \frac{1}{2} \int_{\Omega \backslash B(p,\delta)}  \abs{\nabla G_p}^2 +\abs{G_p}^2 &= \frac{1}{2}\int_{\partial B(p,\delta)} G_p(x) \nabla G_p(x) \cdot \frac{p-x}{\delta} \\
& = \frac{1}{2} \int_{\partial B(p,\delta)} (-\log{\delta} +S_\Omega(p,p) + o_\delta(1)) \left(\frac{1}{\delta}-\frac{\partial S_\Omega}{\partial \nu}\right) \\
& = -\pi \log{\delta} +\pi S_\Omega(p,p) + o_\delta(1).
\end{split}
\end{equation}
We combine \eqref{e_e_w_g_3}, \eqref{e_e_w_g_9}, and \eqref{e_e_w_g_10} to prove \eqref{e_e_w_g_0}.

\end{proof}

\subsection{Lower bounds in $B(p,K\ell) \backslash \mathcal{B}$}\label{lb_blowup}
We now turn to finding bounds in $B(p,K\ell) \backslash \mathcal{B}$.
 In this case, we again use a completion of the
square trick, but the function we use is neither the natural
choice $-\nabla^\bot G_p$ nor a $Y_\varepsilon$ vector field as used before.
Instead, we use a function that arises as the weak-$*$ limit of a
renormalization and blow-up at scale $\ell$ of the superconducting current $j$.  This has the
disadvantage of being tied to
 the functions $u_\varepsilon, A_\varepsilon$
 and not just to the domain $\Omega$.

We define some notation related to $j$.  For $f(x) = \vchi_{[0,1]}(x) +
x^{-1}\vchi_{[1,\infty)}(x)$, we write $\hat{j} = j'f(\abs{u})$.  Define the blow-up of $\hat{j}$ as $\tilde{j}(x) = \ell \hat{j}(p+\ell x) \vchi_{\Omega \backslash \mathcal{B}}(p+\ell x)$.  Now we state a result, which is part of Proposition 9.5 of \cite{ss_book}, that defines which function we use in the completion of the square.

\begin{lem}\label{leftover_lower_bound}
 Suppose configurations $\{(u_\varepsilon,A_\varepsilon)\}$
 satisfy (H1) -- (H4), and that $1\ll n$ as well.
Let $\tilde{j}$ be the blow-up of $\hat{j}$ defined above, and define the blow-up measure on $\Rn{2}$ at scale
$\ell=\sqrt{n/h_{ex}}$ by
\begin{equation*}
\tilde{\mu}(u,A)(x) = \ell^2 \vchi_{\Omega}(\ell x+p)
\mu(u,A) (\ell x+p).
\end{equation*}
Then, up to extraction as $\varepsilon \rightarrow 0$,
\begin{equation}
 \frac{1}{n} \tilde{j} \rightharpoonup J_*
\end{equation}
weakly in $L^2_{loc}(\Rn{2})$, and
\begin{equation}
 \frac{\tilde{\mu}(u_\varepsilon,A_\varepsilon')}{2\pi n} \overset{*}{\rightharpoonup} \mu_*
\end{equation}
weakly-$*$ in the dual of $C_c^{0,\gamma}(\Rn{2})$ for some
$\gamma \in (0,1)$, where $\mu_*$ is a probability measure.
The limits satisfy $\curl{J_*} = 2\pi \mu_*.$  Moreover, as $K \rightarrow \infty$,
\begin{multline}\label{l_l_b_1}
 \frac{1}{2} \int_{B(0,K)} \abs{J_*}^2 = \frac{1}{2}\int_{B(0,K)} \abs{J_* - \nabla^\bot U_*}^2 - \pi\iint \log\abs{x-y}d\mu_*(x)d\mu_*(y) \\ + \pi \log{K} + o_K(1),
\end{multline}
where $U_*$ is the solution to $\Delta U_* = 2\pi \mu_*$ in $B(0,K)$ subject to the boundary condition $U_* =0$ on $\partial B(0,K)$.
\end{lem}

We now use the ``square completion" of Lemma \ref{j_square_completion} in conjunction with the weak convergence of the rescaled and blown-up currents to find energy bounds in $B(p,K\ell) \backslash \mathcal{B}$.

\begin{lem}\label{j_trick}
 Suppose configurations $\{(u_\varepsilon,A_\varepsilon)\}$
 satisfy (H1) -- (H4), and that $1\ll n$ as well.
 Let $J_*$ and $f$ be as defined above, and define the blow-down of $J_*$ by
\begin{equation*}
\bar{J}_*(x) = \frac{1}{\ell} J_*\left(\frac{x-p}{\ell}\right).
\end{equation*}
Then
 \begin{equation}
  \frac{1}{2} \int_{B(p,K\ell)\backslash \mathcal{B}} \abs{\nabla_{A'} u}^2 = \frac{1}{2}\int_{B(p,K\ell)\backslash \mathcal{B}} \abs{\nabla_{A'} u -i u  f(\abs{u}) n \bar{J}_*}^2 + \frac{n^2}{2} \int_{B(0,K)} \abs{J_*}^2 + o(n^2).
 \end{equation}
\end{lem}

\begin{proof}
 Set $X = n f(\abs{u})\bar{J}_*$ in Lemma \ref{j_square_completion} and integrate over $B(p,K\ell)\backslash \mathcal{B}$ to arrive at
 \begin{multline}\label{j_t_1}
  \frac{1}{2} \int_{B(p,K\ell)\backslash \mathcal{B}} \abs{\nabla_{A'} u}^2 = \frac{1}{2}\int_{B(p,K\ell)\backslash \mathcal{B}} \abs{\nabla_{A'} u -i u  f(\abs{u}) n \bar{J}_*}^2 \\ + n^2\int_{B(p,K\ell)} \bar{J}_* \cdot \frac{\hat{j}}{n} \vchi_{\Omega \backslash \mathcal{B}} - \frac{n^2}{2} \int_{B(p,K\ell)\backslash \mathcal{B}} \abs{\bar{J}_*}^2\abs{u f(\abs{u})}^2.
 \end{multline}
We must only deal with the last two integrals.  By the weak convergence $\tilde{j}/n \rightharpoonup J_*$, we have, blowing up via a change of variables $x \mapsto p+\ell x$,
\begin{equation}
 \int_{B(p,K\ell)} \bar{J}_* \cdot \frac{\hat{j}}{n} \vchi_{\Omega \backslash \mathcal{B}} = \int_{B(0,K)} J_* \frac{\tilde{j}}{n} = \int_{B(0,K)} \abs{J_*}^2 + o(1).
\end{equation}
Using the fact that $\abs{u} = 1 +o(1)$ outside of $\mathcal{B}$, and making the same change of variables, we also have that
\begin{equation}\label{j_t_2}
 \int_{B(p,K\ell)\backslash \mathcal{B}} \abs{\bar{J}_*}^2\abs{u f(\abs{u})}^2 = \int_{B(p,K\ell)} \abs{\bar{J}_*}^2 + o(1) = \int_{B(0,K)} \abs{J_*}^2 + o(1).
\end{equation}
Combining \eqref{j_t_1} -- \eqref{j_t_2} yields the result.
\end{proof}

This result can be combined with the properties of $J_*$ to arrive at a more useful estimate.  Indeed, the previous lemma and \eqref{l_l_b_1} immediately yield the following.

\begin{prop}\label{bkl_lower_bound}
Let $\{(u_\ep,A_\ep)\} $ satisfy (H1) -- (H4), suppose $1 \ll n$, and  let $J_*$ and
$\mu_*$ be as above. Then
 \begin{multline}
  \frac{1}{2} \int_{B(p,K\ell)\backslash \mathcal{B}} \abs{\nabla_{A'} u}^2 = \frac{1}{2}\int_{B(p,K\ell)\backslash \mathcal{B}} \abs{\nabla_{A'} u -i u  f(\abs{u}) n \bar{J}_*}^2 + \frac{n^2}{2}\int_{B(0,K)} \abs{J_* - \nabla^\bot U_*}^2 \\
  + \pi n^2 \log{K} -\pi n^2 \iint\log{\abs{x-y}} d\mu_*(x)d\mu_*(y) + o_K(n^2).
 \end{multline}
\end{prop}

\subsection{Synthesis: lower bounds on all of $\Omega$}

Finally, we combine all of the lower bounds of the previous
sections with the energy-splitting result, Proposition
\ref{app_balls}, to find a lower bound for $G_\varepsilon$.

We introduce the vector field $X_\varepsilon$ to be a single field that consists of all of the different fields we have completed the square with.  Indeed, define
\begin{equation}\label{X_def}
 X_\ep := \begin{cases}
   \frac{1}{n}Y_\varepsilon & \text{in }\mathcal{B} \cup \mathcal{A} \\
   f(\abs{u_\varepsilon}) \overline{\nabla^\bot U_*} & \text{in }B(p,K\ell) \backslash \mathcal{B} \\
     - f(\abs{u_\varepsilon})\nabla^\bot G_p & \text{in }\Omega \backslash ( B(p,\delta) \cup \mathcal{B}).
 \end{cases}
\end{equation}
Here we have used the notation $\overline{\nabla^\bot U_*}$ for the blow-down at
scale $\ell$ of $\nabla^\bot U_*$, where $U_*$ is defined in Lemma \ref{leftover_lower_bound}.  The function $G_p$ is defined by \eqref{Gp}, and $Y_\varepsilon$ is the field initially defined in $\mathcal{B}$ by Theorem \ref{energy_bound} and then extended to $\mathcal{A} \backslash \mathcal{B}$ according to \eqref{ann_G_def}.  Although $X_\ep$ also depends on $ \delta$ and
$K$, we do not build the dependence into the notation. We also
recall the definitions of $f_\ep$ in \eqref{fep_def} and $I$ in
\eqref{I}.

Our main result of this section implies the first assertion
\eqref{asser1}  of Theorem \ref{3}:
\begin{thm}\label{squeeze_lower}
  Suppose configurations $\{(u_\varepsilon,A_\varepsilon)\}$
  satisfy (H1) -- (H4), and that $1 \ll n$.
Let $X_\ep$ be the vector field defined by \eqref{X_def}.  Then
for $\varepsilon$ sufficiently small, and as $K\rightarrow \infty$, $\delta \rightarrow 0$,
\begin{multline}
  G_\varepsilon(u_\varepsilon,A_\varepsilon) \ge f_\varepsilon(n) + n^2 I(\mu_*) + o_{K,\delta}(n^2) +o(n^2) 
   + \frac{n^2}{36} \int_{\Omega}
  \abs{\frac{1}{n}\nabla_{A_\varepsilon'} u_\varepsilon - i u_\varepsilon X}^2\\
 + \frac{1}{5} \int_{B(p,\delta) \cup \mathcal{B}} \abs{\curl{A_\varepsilon'}}^2
    + \frac{1}{2}\int_{\Omega\backslash(B(p,\delta)\cup \mathcal{B})} \abs{\curl{A_\varepsilon'}  - n G_p}^2.
\end{multline}
\end{thm}

\begin{proof}
For convenience we drop the subscript $\varepsilon$ for the rest of the proof. Applying both bounds of Proposition \ref{app_balls} provides the initial lower bound
\begin{equation}\label{s_l_1}
 G_\varepsilon(u,A)  \ge I + II + III,
\end{equation}
where
\begin{equation}
 I:= h_{ex}^2 J_0 + \pi n \log{\frac{r}{n\varepsilon}} + \frac{1}{36} \int_\mathcal{B} \abs{\nabla_{A'} u-iuY}^2,
\end{equation}
\begin{equation}
 II: = F_\varepsilon(u,A',\Omega \backslash \mathcal{B}) + 2\pi h_{ex} \int \xi_0 d\nu + \frac{1-r^2}{2} \int_\mathcal{B} \abs{\curl{A'}}^2,
\end{equation}
and
\begin{equation}
 III: = \frac{\pi \alpha}{4} (n'-n)\abs{\log\varepsilon} -Crh_{ex} (n'-n) - Ch_{ex}\varepsilon^{3\alpha/2 -1} - Ch_{ex}^2 \varepsilon^\alpha  -Cn.
\end{equation}

The term $III$ is easiest to deal with, so we dispatch it first.  The hypotheses $h_{ex} \le C\varepsilon^{-\beta}$, $2/3 < \alpha <1$, and $0 < \beta < 3\alpha/2-1$ imply that
\begin{equation}\label{s_l_2}
 \begin{split}
  &\varepsilon^{3\alpha/2 -1} h_{ex} = o(1)  \\
  &h_{ex}^2 \varepsilon^{\alpha} = o(1).
 \end{split}
\end{equation}
Since $r = 1/\sqrt{h_{ex}}$, we have that $r h_{ex} = \sqrt{h_{ex}}$, and hence the hypothesis 
(H3) implies that
\begin{equation}\label{s_l_2_0}
 \frac{\pi \alpha}{4} (n'-n)\abs{\log\varepsilon} -Crh_{ex} (n'-n) = (n'-n)\left( \frac{\pi \alpha}{4} \abs{\log\varepsilon} - C\sqrt{h_{ex}} \right) \ge 0
\end{equation}
for $\varepsilon$ sufficiently small.  These bounds and that $1 \ll n$ then imply that
\begin{equation}\label{s_l_3}
 III \ge o(1) - o(n^2).
\end{equation}

The terms in $II$ essentially constitute the energy content of the exterior of the balls; in bounding $II$ we will employ all of the estimates of the previous sections.  We begin by splitting $F_\varepsilon(u,A',\Omega \backslash \mathcal{B})$ into parts corresponding to the different regions considered in the previous sections.  Indeed, we write
\begin{equation}
 II = IV + V + VI,
\end{equation}
where
\begin{equation}
 IV := F_\varepsilon(u,A',\mathcal{A} \backslash \mathcal{B}) + 2\pi h_{ex} \int \xi_0 d\nu + \frac{1-r^2}{2} \int_\mathcal{B} \abs{\curl{A'}}^2
\end{equation}
is the energy content of the annulus around $p$,
\begin{equation}
 V := F_\varepsilon(u,A',\Omega\backslash(B(p,\delta)\cup \mathcal{B}))
\end{equation}
is the energy content outside of the ball $B(p,\delta)$, and
\begin{equation}
 VI := F_\varepsilon(u,A',B(p,K\ell)\backslash \mathcal{B})
\end{equation}
is the energy content in the ball $B(p,K\ell)$.  In the annulus around $p$, Proposition \ref{app_annuli} shows that
\begin{equation}\label{s_l_4}
\begin{split}
 IV & \ge   \frac{1}{36} \int_{\mathcal{A} \backslash \mathcal{B}}
  \abs{\nabla_{A'} u -iuY}^2 +  \frac{1}{4} \int_{\mathcal{A} \backslash \mathcal{B}} \abs{\curl{A'}}^2 + \left(\frac{1}{4} - \frac{r^2}{2} \right)\int_{\mathcal{B}} \abs{\curl{A'}}^2 \\
  &+ \pi n^2 \log{ \frac{\delta}{K\ell} } + 2 \pi n h_{ex} \underline{\xi_0}
  + 2 \pi h_{ex} \sum_{\substack{b_i \in B(p,K\ell) \\ d_i >0}} d_i(\xi_0(b_i)
   - \underline{\xi_0})
  - \pi n^2 \delta^2 + o(n^2).
\end{split}
\end{equation}
In $B(p,\delta)$ we use Proposition \ref{exterior_est_w_g} to estimate
\begin{equation}\label{s_l_5}
\begin{split}
 V &\ge  \frac{1}{2} \int_{\Omega\backslash(B(p,\delta)\cup \mathcal{B})} \abs{\nabla_{A'} u + iu n f(\abs{u})\nabla^\bot G_p}^2 + \frac{1}{2}\int_{\Omega\backslash(B(p,\delta)\cup \mathcal{B})} \abs{\curl{A'}  - n G_p}^2 \\&- \pi n^2 \log{\delta} +\pi n^2 S_\Omega(p,p) +o_\delta(n^2) + o(n^2).
\end{split}
\end{equation}
Finally, in $B(p,K\ell)$ we utilize Proposition \ref{bkl_lower_bound} to get
\begin{equation}\label{s_l_6}
\begin{split}
 VI & \ge  \frac{1}{2}\int_{B(p,K\ell)\backslash \mathcal{B}} \abs{\nabla_{A'} u -i u  f(\abs{u}) n \bar{J}_*}^2 + \frac{n^2}{2}\int_{B(0,K)} \abs{J_* - \nabla^\bot U_*}^2 + \frac{1}{2}\int_{B(p,K\ell)\backslash \mathcal{B}} \abs{\curl{A'}}^2\\
  &+ \pi n^2 \log{K} -\pi n^2 \iint\log{\abs{x-y}} d\mu_*(x)d\mu_*(y) + o_K(n^2).
\end{split}
\end{equation}
By changing variables to blow-down at scale $\ell$, we have that
\begin{equation}\label{s_l_7_0}
\begin{split}
\frac{n^2}{2}\int_{B(0,K)} \abs{J_* - \nabla^\bot U_*}^2 &= \frac{n^2}{2}\int_{B(p,K\ell)} \abs{\bar{J_*} - \overline{\nabla^\bot U_*}}^2 \\
&\ge \frac{1}{2}\int_{B(p,K\ell) \backslash \mathcal{B}} \abs{iu f(\abs{u})  n\bar{J_*} - iu f(\abs{u}) n\overline{\nabla^\bot U_*}}^2 +o(n^2),
\end{split}
\end{equation}
where we have used the fact that $\abs{u} = 1 +o(1)$ in $B(p,K\ell) \backslash \mathcal{B}$.  Hence,
\begin{multline}\label{s_l_7_1}
\frac{1}{2}\int_{B(p,K\ell)\backslash \mathcal{B}} \abs{\nabla_{A'} u -i u  f(\abs{u}) n \bar{J}_*}^2 + \frac{n^2}{2}\int_{B(0,K)} \abs{J_* - \nabla^\bot U_*}^2 \\
\ge \frac{n^2}{4} \int_{B(p,K\ell)\backslash \mathcal{B}} \abs{\frac{1}{n}\nabla_{A'} u -i u  f(\abs{u})  \overline{\nabla^\bot U}_*}^2
\end{multline}
We now note that the proof of Proposition 9.1 of \cite{ss_book} shows the inequality
\begin{equation}\label{s_l_7}
 2 \pi h_{ex} \sum_{\substack{b_i \in B(p,K\ell) \\ d_i >0}} d_i(\xi_0(b_i) - \underline{\xi_0}) \ge
 \pi n^2 \int Q(x) d\mu_*(x) +o_K(1),
\end{equation}
where $Q$ is the quadratic form of the Hessian of $\xi_0$ at $p$.  Also, the assumption that $10 \le h_{ex}$ implies that $1/4 - r^2/2 \ge 1/5$.  Then, summing \eqref{s_l_4} -- \eqref{s_l_6} and employing \eqref{s_l_7_0} -- \eqref{s_l_7}, we arrive at the bound
\begin{equation} \label{s_l_8}
\begin{split}
 II &\ge  \pi n^2 \log{ \frac{1}{\ell} } + 2 \pi n h_{ex} \underline{\xi_0} +\pi n^2 S_\Omega(p,p) + n^2 I(\mu_*) + o_\delta(n^2) + o_K(n^2) + o(n^2) \\
& +\frac{n^2}{36} \int_{\mathcal{A} \backslash \mathcal{B}}
  \abs{\frac{1}{n}\nabla_{A'} u - i u\frac{1}{n}Y }^2 +\frac{n^2}{4} \int_{B(p,K\ell)\backslash \mathcal{B}} \abs{\frac{1}{n}\nabla_{A'} u -i u  f(\abs{u})  \overline{\nabla^\bot U}_*}^2 \\
& + \frac{n^2}{2} \int_{\Omega\backslash(B(p,\delta)\cup \mathcal{B})} \abs{\frac{1}{n}\nabla_{A'} u + iu f(\abs{u})\nabla^\bot G_p}^2
+ \frac{1}{5} \int_{B(p,\delta) \cup \mathcal{B}} \abs{\curl{A'}}^2 \\
&+  \frac{1}{2}\int_{\Omega\backslash(B(p,\delta)\cup \mathcal{B})} \abs{\curl{A'}  - n G_p}^2 .
\end{split}
\end{equation}

In order to conclude we turn back to $I$.  Note that since $r = 1/\sqrt{h_{ex}}$ and $\ell = \sqrt{n/h_{ex}}$, we have that
\begin{equation}\label{s_l_9}
 n \log{\frac{r}{n\varepsilon}} = n \left( \log{\frac{\ell}{\varepsilon}} - \frac{3}{2}\log{n} \right)
 = n \log{\frac{\ell}{\varepsilon}} -\frac{3n}{2}\log{n} =n \log{\frac{\ell}{\varepsilon}} -o(n^2) .
\end{equation}
We use this expansion in $I$; the result follows in view of the
definition of $X_\ep$ from summing
 $I$, \eqref{s_l_3}, and \eqref{s_l_8}.

\end{proof}

\begin{remark} The condition $10 \le h_{ex}$ can be relaxed at the cost of
a different constant in front of the $\int_{B(p,\delta)\cup
\mathcal{B}} \abs{\curl{A'}}^2$ term.  Indeed, for any
$\gamma<1/4$, the relaxed condition $h_{ex}\ge 2/(1-4\gamma)$ puts
a term $\gamma$ in front of the $\curl{}$ integral.  If we are
willing to drop the $\curl{}$ integral altogether in the lower
bound, we may set $\gamma=0$ and assume only that $h_{ex} \ge 2$,
the minimum requirement for $1/4 - r^2/2$ to be nonnegative.
\end{remark}

\section{The case of a priori upper bounds on $G_\varepsilon$}\label{G_ep_apriori}
We now utilize the result of Theorem \ref{squeeze_lower}
 in conjunction with various a priori upper bounds
on the full energy $G_\varepsilon$.  The result is the following
more  versatile extension of Proposition \ref{degree_control},
which together with Theorem \ref{squeeze_lower} implies Theorem
\ref{3}.  The main assertion is that  $\pqnorm{\nabla_{A'}
u}{2}{\infty}$ is of the same order as $n$. As mentioned in the
introduction, this is significant in that the quantity $n$ is
determined by the ball construction, and is thus not intrinsically
defined, whereas  $\|\nab_{A'} u\|_{\lti}$ is.

\begin{thm}\label{application_1}
  Suppose configurations $\{(u_\varepsilon,A_\varepsilon)\}$
  satisfy (H1) -- (H4) and that $1 \ll n$.

1. It always holds that, as $\ep \to 0$,
\begin{equation}\label{app_1_02}
\pnormspace{\frac{1}{n}\nabla_{A_\varepsilon'} u_\varepsilon -
iu_\varepsilon X_\varepsilon}{2}{\Omega} = O(1),
\end{equation}
\begin{equation}\label{app_1_01}
\frac{1}{n}\pqnormspace{\nabla_{A_\varepsilon'}
u_\varepsilon}{2}{\infty}{\Omega} \le \pqnormspace{\nabla
G_p}{2}{\infty}{\Omega} + C(C_0+1),
\end{equation}
where $C$ depends on $\Omega$ and $C_0$ is the constant from the bound (H4),  and
\begin{equation}\label{app_1_03}
\pnormspace{\frac{1}{n}\curl{A_\varepsilon'}}{2}{\Omega} = O(1).
\end{equation}


 2. If we
assume that the a priori bound
\begin{equation}\label{app_1_000}
G_\varepsilon(u_\varepsilon,A_\varepsilon) \le f_\varepsilon(n) +
n^2 I(\mu_*)  + o(n^2)
\end{equation}
holds, then, as $\varepsilon \rightarrow 0$, 
\begin{equation}\label{app_1_0002}
\pnormspace{\frac{1}{n}\nabla_{A_\varepsilon'} u_\varepsilon -
iu_\varepsilon X_\varepsilon}{2}{\Omega} =  o(1) + o_{K,\delta}(1),
\end{equation}
and
\begin{equation}\label{app_1_0001}
\pqnormspace{\nabla G_p}{2}{\infty}{\Omega} - o(1) \le \frac{1}{n}
\pqnormspace{\nabla_{A_\varepsilon'}
u_\varepsilon}{2}{\infty}{\Omega} \le \pqnormspace{\nabla
G_p}{2}{\infty}{\Omega} + C  + o(1),
\end{equation}
where $C$ is a universal constant.
Also,
\begin{equation}\label{app_1_0003}
\frac{1}{n} \curl{A_\varepsilon'} \rightarrow G_p
\end{equation}
strongly in $L^2(\Omega)$.

\end{thm}
 The  reason for the lack of strong convergence  in $\lti$ of
$\frac{1}{n} \nab_{A'} u$, or rather  of $X_\ep$, to $-\np G_p$
was given in the introduction, Section \ref{oc}.  This explains the presence of the constant $C$ in \eqref{app_1_0001}. We start with two lemmas that establish that we can, however, estimate $\|X_\ep\|_{\lti(\om)} $ in
terms of $\|\nab G_p\|_{\lti(\om)}$, but without equality. Note
that the second lemma is the place where we crucially apply  the
$\lti$ control in the balls proved in the previous paper
\cite{p1}.
\begin{lem}\label{annular_Gp_swap}
 Let $X_\varepsilon$ be the vector field defined by \eqref{X_def}.  Define the set
\begin{equation*}
W_\varepsilon =  \supp(X_\varepsilon) \cap (\Omega\backslash (B(p,K\ell) \cup \mathcal{B})).
\end{equation*}
Then
\begin{equation}\label{a_g_s_0}
 \pqnormspace{X_\varepsilon+ \nabla^\bot G_p}{2}{\infty}{W_\varepsilon} = o(1) +o_{K,\delta}(1).
\end{equation}
\end{lem}

\begin{proof}
Again suppress the subscripts $\varepsilon$.  To begin we decompose $W$ into two components, one part inside the annulus $\mathcal{A}$ and the other outside.  Let $W_1 =\Omega \backslash (B(p,\delta) \cup \mathcal{B})$ and
 \begin{equation*}
  \mathcal{A_{T}} = \{x \in \mathcal{A} \;\vert\; \abs{x-p} \notin T  \},
 \end{equation*}
where $T$ is defined in \eqref{app_t_def}.  Then $W = W_1 \cup \mathcal{A}_T$ is a disjoint union, and we may trivially bound
\begin{equation}\label{a_g_s_1}
 \pqnormspace{X+ \nabla^\bot G_p}{2}{\infty}{W} \le \pqnormspace{X+ \nabla^\bot G_p}{2}{\infty}{W_1} + \pqnormspace{X+ \nabla^\bot G_p}{2}{\infty}{\mathcal{A}_T}.
\end{equation}

On $W_1$ we have that $X = -f(\abs{u}) \nabla^\bot G_p$, and hence $X + \nabla^\bot G_p = (1-f(\abs{u}))\nabla^\bot G_p$ there.  By construction, $\abs{u} = 1+o(1)$ on $W_1$, so $f(\abs{u}) = 1+o(1)$, and we may bound
\begin{equation}\label{a_g_s_2}
 \pqnormspace{X+ \nabla^\bot G_p}{2}{\infty}{W_1} \le o(1) \pqnormspace{\nabla^\bot G_p}{2}{\infty}{W_1} \le o(1) \pqnormspace{\nabla^\bot G_p}{2}{\infty}{\Omega} =o(1).
\end{equation}
On $\mathcal{A}_T$ we have that
\begin{equation*}
 X(x) = \frac{D(t)}{n} \frac{\tau_p}{\abs{x-p}},
\end{equation*}
where $D(t)$ is defined by \eqref{app_d_def}.  The decomposition $G_p(x) = -\log{\abs{x-p}} + S_\Omega(x,p)$ yields $\nabla^\bot G_p(x) = -\frac{\tau_p}{\abs{x-p}} + \nabla^\bot S_\Omega(x,p)$.  Hence, on $\mathcal{A}_T$ we have
 \begin{equation*}
  \frac{D(t)}{n} \frac{\tau_p}{\abs{x-p}} + \nabla^\bot G_p = \left(\frac{D(t) - n}{n} \right)\frac{\tau_p}{\abs{x-p}} + \nabla^\bot S_\Omega(x,p).
 \end{equation*}
We may use \eqref{app_annuli_0} to estimate
\begin{equation*}
 \abs{\frac{D(t) -n}{n}} \le C\ell^2 \left(\frac{1}{t^2}+1\right) \le C\left(\frac{1}{K^2}  + \ell^2  \right)
\end{equation*}
since $t \in (K\ell,\delta)$.  Using this, we see that
\begin{equation}\label{a_g_s_3}
 \begin{split}
  \pqnormspace{X + \nabla^\bot G_p}{2}{\infty}{\mathcal{A}_T}
  & \le C\left(\frac{1}{K^2}  + \ell^2\right) \pqnorm{\frac{1}{\abs{\cdot-p}}}{2}{\infty} + \pqnormspace{\nabla^\bot S_\Omega(\cdot,p)}{2}{\infty}{\mathcal{A}_T}\\
  &=2\sqrt{\pi}(o(1) + o_K(1) ) + o_\delta(1).
 \end{split}
\end{equation}
Here we have used $\pqnorm{1/\abs{x}}{2}{\infty} = 2\sqrt{\pi}$; the fact that $\nabla^\bot S_\Omega(\cdot,p)$ is continuous and $\abs{\mathcal{A}_T} = o_\delta(1)$ allows us to write
$\pqnormspace{\nabla^\bot
S_\Omega(\cdot,p)}{2}{\infty}{\mathcal{A}_T} = o_\delta(1)$.  The
result \eqref{a_g_s_0} follows from \eqref{a_g_s_1} --
\eqref{a_g_s_3}.
\end{proof}

Now we show that $\pqnorm{u_\varepsilon X_\varepsilon}{2}{\infty}$ can be estimated above and below by $\pqnorm{\nabla G_p}{2}{\infty}$.

\begin{lem}\label{X_control}
Let $X_\varepsilon$ be the vector field defined by \eqref{X_def}.  We have that, for $\varepsilon$ sufficiently small, 
\begin{equation}\label{X_c_0}
\pqnormspace{\nabla G_p}{2}{\infty}{\Omega\backslash (B(p,\delta) \cup \mathcal{B})} - o(1) \le  \pqnormspace{u_\varepsilon X_\varepsilon}{2}{\infty}{\Omega} \le \pqnormspace{\nabla G_p}{2}{\infty}{\Omega} + C(C_0+1)  +o_{K,\delta}(1),
\end{equation}
where $C_0$ is the constant from (H4) and $C$ depends on $\Omega$. 
\end{lem}
\begin{proof}
 Recall that, by construction, $\abs{u} = 1+o(1)$ on $\Omega \backslash \mathcal{B}$.  In the balls, $\mathcal{B}$, the construction of the vector field $Y$ is such that $\abs{u} \le 3/2$ on $\supp(Y) \cap \mathcal{B}$ (see Proposition 4.3 of \cite{p1}).  Since $X= -f(\abs{u})\nabla^\bot G_p$ on $\Omega \backslash(B(p,\delta)\cup \mathcal{B})$, the lower bound follows from the pointwise inequality
 \begin{equation}
  (1-o(1)) \abs{\nabla G_p(x)}  \le \abs{u(x)} \abs{X(x)} \text{ for } x\in \Omega \backslash(B(p,\delta)\cup \mathcal{B}).
 \end{equation}

For the upper bound we define the sets
\begin{equation*}
 \begin{split}
  & \Omega_1 = \Omega \backslash (B(p,K\ell)\cup \mathcal{B}) \\
  & \Omega_2 = B(p,K\ell) \backslash \mathcal{B} \\
  & \Omega_3 = \mathcal{B},
 \end{split}
\end{equation*}
where, as usual we abuse notation by writing $\mathcal{B}$ for $\cup_{B\in\mathcal{B}} B$.  We then have that
\begin{equation*}
 \pqnormspace{uX}{2}{\infty}{\Omega} \le \pqnormspace{uX}{2}{\infty}{\Omega_1}+\pqnormspace{uX}{2}{\infty}{\Omega_2}+\pqnormspace{uX}{2}{\infty}{\Omega_3},
\end{equation*}
and we estimate each term separately.  On $\Omega_1$ we apply Lemma \ref{annular_Gp_swap} to see  that
\begin{equation}
 \begin{split}
 \pqnormspace{uX}{2}{\infty}{\Omega_1} & = \pqnormspace{uX}{2}{\infty}{\Omega_1 \cap \supp(X)}\\
  & \le (1+o(1)) \pqnormspace{X}{2}{\infty}{\Omega_1 \cap \supp(X)}  \\
  & \le (1+o(1)) (\pqnormspace{ \nabla G_p}{2}{\infty}{\Omega_1} + \pqnormspace{X+\nabla^\bot G_p}{2}{\infty}{\Omega_1 \cap \supp(X)})
   \\ &\le \pqnormspace{\nabla G_p}{2}{\infty}{\Omega} + o(1) + o_{K,\delta}(1).
 \end{split}
\end{equation}

For $\Omega_3$ we employ Proposition 6.4 of \cite{p1} along with the bound $\abs{u} \le 3/2$ to bound
\begin{equation}
 \pqnormspace{uX}{2}{\infty}{\Omega_3}^2 \le \frac{9}{4n^2} \pqnormspace{Y}{2}{\infty}{\Omega_3}^2 \le  \frac{C}{n^2}\left(F_\varepsilon^r(u, A',\mathcal{B}) - \pi n \left(\log{\frac{r}{\varepsilon n}} \right) + n^2 \right),
\end{equation}
where $C$ is a universal constant and 
\begin{equation}
F_\varepsilon^r(u,A',\mathcal{B}) = \frac{1}{2} \int_{\mathcal{B}} \abs{\nabla_{A'} u}^2 + r^2\abs{\curl{A'}}^2 + \frac{1}{2\varepsilon^2}(1-\abs{u}^2)^2.
\end{equation} 
We claim that for $\varepsilon$ sufficiently small, 
\begin{equation}
 F^r_\varepsilon(u, A', \mathcal{B}) - \pi n \log{\frac{r}{n\varepsilon}} \le C_0 n^2 + Cn^2,
\end{equation}
where $C_0$ is the constant from (H4) and $C$ depends on $\Omega$.  This immediately implies that 
\begin{equation}
 \pqnormspace{uX}{2}{\infty}{\Omega_3} \le C(C_0 +1),
\end{equation}
where $C$ depends on $\Omega$.  To prove this claim we must use a modification of lower bounds argument of Theorem \ref{squeeze_lower}, and then compare with the matching upper bound (H4).  We argue as we did in \eqref{s_l_1}, except now we wish to retain the term $F^r_\varepsilon(u, A', \mathcal{B})$ rather than bounding it from below by Lemma \ref{redux_92}.  This yields a lower bound of the form  $G_\varepsilon(u,A) \ge I + II$, where $II$ is identical to the $II$ used in \eqref{s_l_1}, and 
\begin{equation}
 I := h_{ex}^2 J_0 + F^r_\varepsilon(u, A', \mathcal{B}) - C\sqrt{h_{ex}} (n'-n) - Ch_{ex}\varepsilon^{3\alpha/2 -1} - Ch_{ex}^2 \varepsilon^\alpha. 
\end{equation}
The term $II$, which corresponds to the free energy outside the balls, we bound in exactly the same way, yielding \eqref{s_l_8}.  For the purposes of the claim we can disregard all of the integrals in the second and third lines of \eqref{s_l_8} and retain only the first line.  For $I$ we employ \eqref{s_l_2} to write its last two terms as $o(1)$.  Combining this with the bounds on $II$ and comparing to the upper bound (H4), we have that, for $\varepsilon$ sufficiently small, 
\begin{equation}\label{X_c_1}
 \pi n \log{\frac{\ell}{\varepsilon}} + C_0 n^2 \ge F^r_\varepsilon(u, A', \mathcal{B}) + n^2 I(\mu_*) - C\sqrt{h_{ex}} (n'-n) + o_{K,\delta}(n^2) +o(n^2). 
\end{equation}
Letting $K\rightarrow \infty$ and $\delta \rightarrow 0$, we absorb the $o_{K,\delta}(n^2)$ term into the $o(n^2)$ term.  The functional $I(\cdot)$, defined over probability measures, has a unique minimizer $\mu_0$ (see \cite{saff_totik}); we bound $I(\mu_*)\ge I(\mu_0)$, a constant that depends only on $\Omega$.  We now borrow half of $F^r_\varepsilon(u,A',\mathcal{B})$ and use Lemma 9.1 of \cite{ss_book} to bound
\begin{equation}\label{X_c_2}
 \frac{1}{2} F^r_\varepsilon(u, A', \mathcal{B}) \ge \frac{\pi n}{2}\log{\frac{r}{n\varepsilon}} + \frac{\pi \alpha}{4}(n'-n)\log\frac{1}{\varepsilon} - Cn.
\end{equation}
Putting \eqref{X_c_2} into \eqref{X_c_1} and employing \eqref{s_l_2_0} to deal with the $n'-n$ terms and \eqref{s_l_8} to rewrite the  $\log{\ell/\varepsilon}$ term, we find that
\begin{equation}\label{X_c_3}
\frac{3\pi n}{2}\log{n} + C_0 n^2 - n^2 I(\mu_0) + Cn -o(n^2) \ge \frac{1}{2} \left(F^r_\varepsilon(u, A', \mathcal{B}) - \pi n \log{\frac{r}{\varepsilon n}} \right).
\end{equation}
The claim follows.

In $\Omega_2$ we note that $\abs{uX} = f(\abs{u})\abs{u} \abs{\overline{\nabla U_*}} \le \abs{\overline{\nabla U_*}}$.  This and a blow-up at scale $\ell$ imply that
\begin{equation}
 \pqnormspace{uX}{2}{\infty}{\Omega_2} \le \pqnormspace{\overline{\nabla U_*}}{2}{\infty}{\Omega_2} \le \pqnormspace{\nabla U_*}{2}{\infty}{B(0,K)}.
\end{equation}
Since $\Delta U_* = 2\pi \mu_*$ in $B(0,K)$ with vanishing Dirichlet boundary condition, we may write
\begin{equation}
 U_*(x) = \int_{B(0,K)} H_K(x,y) d\mu_*(y),
\end{equation}
 where
\begin{equation*}
 H_K(x,y) = \log{\abs{x-y}} - \log{\abs{ x\frac{K}{\abs{x}} - y \frac{\abs{x}}{K}  }}
\end{equation*}
is the Green's kernel on $B(0,K)$.  Since $H_K(x,y) = H_1(x/K,y/K)$, and the gradient of the $H_1$ kernel can have at worst a singularity like $1/\abs{x-y}$, we have that
\begin{equation}
 \sup_{y\in B(0,K)} \pqnormspace{\nabla H_K(\cdot,y)}{2}{\infty}{B(0,K)} = \sup_{y\in B(0,1)} \pqnormspace{\nabla H_1(\cdot,y)}{2}{\infty}{B(0,1)} = C < \infty
\end{equation}
for some universal constant $C$ that does not depend on $K$.  Now, for any set $E\subseteq B(0,K)$, we have that
\begin{equation}
\begin{split}
 \int_E \abs{\nabla U_*(x)} dx & \le \int_E \int_{B(0,K)}  \abs{\nabla H_K(x,y)} d\mu_*(y) dx \\
 & = \int_{B(0,K)} \int_E \abs{\nabla H_K(x,y)} dx d\mu_*(y) \\
 & \le \int_{B(0,K)} \abs{E}^{1/2} \pqnormspace{\nabla H_K(\cdot,y)}{2}{\infty}{B(0,K)} d\mu_*(y) \\
 & = \abs{E}^{1/2}  \sup_{y\in B(0,K)} \pqnormspace{\nabla H_K(\cdot,y)}{2}{\infty}{B(0,K)},
\end{split}
 \end{equation}
where we have utilized the fact that $\mu_*$ is a probability measure.  Hence
\begin{equation}
 \pqnormspace{\nabla U_*}{2}{\infty}{B(0,K)} \le \sup_{y\in B(0,K)} \pqnormspace{\nabla H_K(\cdot,y)}{2}{\infty}{B(0,K)} = C < \infty.
\end{equation}
The result follows.

\end{proof}

\begin{remark}\label{constant_shuffle}
The dependence of the term $C(C_0+1)$ on the domain is only through the dependence of $I(\mu_0)$ on the quadratic form of $D^2 \xi_0$.  If the stronger a priori upper bound $G_\varepsilon(u_\varepsilon,A_\varepsilon) \le f_\varepsilon(n) + n^2 I(\mu_*)  + o(n^2)$ holds, then we may replace the $C_0 n^2 - n^2 I(\mu_0)$ term on the left side of \eqref{X_c_3} with $o(n^2)$.  In \eqref{X_c_0}, this allows us to replace the term  $C(C_0+1)$ with a universal constant $C$. 
\end{remark}

We can now conclude the
\begin{proof}[Proof of Theorem \ref{application_1}]
For convenience we drop the subscript $\varepsilon$ on
$u_\varepsilon$, $A_\varepsilon$, and $X_\varepsilon$.  The hypotheses allow us to
employ Theorem \ref{squeeze_lower} for a lower bound on
$G_\varepsilon(u,A)$.  Comparing this with the upper bound
\eqref{H4} and dividing by $n^2$, we find that
\begin{equation}
\pnormspace{\frac{1}{n}\nabla_{A'} u - iuX}{2}{\Omega}^2 \le  C_0 + o(1) + o_{K,\delta}(1) .
\end{equation}
This implies \eqref{app_1_02}.  From this and the bound $\pqnorm{\cdot}{2}{\infty} \le \pnorm{\cdot}{2}$, we have that
\begin{equation}
\abs{\frac{1}{n}\pqnormspace{\nabla_{A'} u}{2}{\infty}{\Omega} - \pqnormspace{uX}{2}{\infty}{\Omega}} \le C_0 + o(1) + o_{K,\delta}(1).
\end{equation}
Moreover, using Lemma \ref{X_control} and letting $\delta \rightarrow 0$, $K \rightarrow \infty$ implies that
\begin{equation}
\frac{1}{n} \pqnormspace{\nabla_{A'} u}{2}{\infty}{\Omega} \le   \pqnormspace{\nabla G_p}{2}{\infty}{\Omega} + C(C_0+1) + o(1),
\end{equation}
where $C$ depends on $\Omega$ and $C_0$ is from the bound (H4).  This is \eqref{app_1_01}.  A similar argument, using the extra terms in the lower bounds of Theorem \ref{squeeze_lower}, proves \eqref{app_1_03}.

Suppose now that the bound \eqref{app_1_000} holds.  Then, again comparing with the bound from Theorem \ref{squeeze_lower}, we find that
\begin{equation}
\pnormspace{\frac{1}{n}\nabla_{A'} u - iuX}{2}{\Omega}^2 \le  o(1) + o_{K,\delta}(1).
\end{equation}
This is \eqref{app_1_0002}.  We use Lemma \ref{X_control} and Remark \ref{constant_shuffle} and let $\delta \rightarrow 0$ and $K \rightarrow \infty$ to arrive at the bounds
\begin{equation}
\begin{split}
 \pqnormspace{\nabla G_p}{2}{\infty}{\Omega}  - o(1) &\le \frac{1}{n} \pqnormspace{\nabla_{A'} u}{2}{\infty}{\Omega} \\
&\le \pqnormspace{\nabla G_p}{2}{\infty}{\Omega} + C  + o(1),
\end{split}
\end{equation}
where $C$ is a universal constant.  This is \eqref{app_1_0001}.  A similar argument proves that 
\begin{equation}
\pnormspace{\frac{1}{n} \curl{A'} - G_p\vchi_{\Omega\backslash (B(p,K\ell) \cup \mathcal{B})}  }{2}{\Omega} = o(1) + o_{K,\delta}(1).
\end{equation} 
Then
\begin{equation}
\begin{split}
\pnormspace{\frac{1}{n} \curl{A'} - G_p}{2}{\Omega} & \le \pnormspace{\frac{1}{n} \curl{A'} - G_p\vchi_{\Omega\backslash (B(p,K\ell) \cup \mathcal{B})}  }{2}{\Omega} + \pnormspace{G_p}{2}{B(p,K\ell) \cup \mathcal{B}} \\
& \le o(1) + o_{K,\delta}(1) + \pnormspace{G_p}{2}{\mathcal{B}} + \pnormspace{G_p}{2}{B(p,K\ell)}.
\end{split}
\end{equation}
Let $\varepsilon \rightarrow 0$ and then send $\delta \rightarrow 0$ and $K \rightarrow \infty$.  Then the right hand side tends to zero and \eqref{app_1_0003} follows.

\end{proof}

\section{Convergence results}\label{convergence_results}
This section provides several applications of Theorem \ref{application_1}.
Throughout we will assume that (H1) -- (H4) hold, and that $1 \ll n$.

\subsection{Compactness of Jacobians}\label{cpt_jac}

In this section we will use the results of Theorem
\ref{application_1} to prove compactness of the gauge-invariant
Jacobians (defined by \eqref{jacdef}) in a function space
based on Lorentz spaces, which we call $\mathcal{X}(\Omega)$.  We recall
(see \eqref{jacest}) that the  best estimates and compactness
results for Jacobians in the literature are in the dual of the
H\"older spaces $C^{0,\gamma}_c$ (their limit being generally
bounded Radon measures).

 Before defining the space $\mathcal{X}$, we recall
the main problem that leads to considering it. In two space
dimensions the exponent $p=2$ is critical in the sense that its
Sobolev conjugate $2^* = \frac{4}{2-2} = \infty$. This leads to
embeddings $H^1_0(\Omega) \hookrightarrow L^p(\Omega)$ for each $1
\le p < \infty$, but the embedding into $L^\infty(\Omega)$ fails.
Indeed, it is possible to construct functions in $H^1_0(\Omega)$
that are unbounded in a neighborhood of every point of $\Omega$
(see Section 5.6 of \cite{stein_sing} for details).  However, for
any $p>2$, we get embeddings $W_0^{1,p}(\Omega) \hookrightarrow
C^{0,\gamma}(\om)$ with $\gamma = 1 - 2/p$.  We thus see a sharp
transition from $p=2$, where we can find very poorly behaved
functions, to $p>2$ where we have gained enough control so that
the functions are H\"older continuous. This suggests that it might
be possible to find an intermediate space $\mathcal{X}$,
\begin{equation*}
 W_0^{1,p}(\Omega) \hookrightarrow \mathcal{X}(\Omega) \hookrightarrow H_0^1(\Omega)  \text{ for all } p>2,
\end{equation*}
such that $\mathcal{X}(\Omega) \hookrightarrow C^0(\Omega)$.

Since the Sobolev spaces consist of functions whose weak
derivatives are in some $L^p$ space, it is natural to look to the
Lorentz spaces, $L^{p,q}$, which are generalizations of $L^p$
spaces, in order to define $\mathcal{X}(\Omega)$.  Though we will only use
two of the Lorentz spaces, $L^{2,1}$ and $L^{2,\infty}$, we give
the definition for all $1\le p \le \infty,$ $1\le q \le \infty$.
See Chapter 5 of \cite{stein_intro} or Chapter 1 of
\cite{grafakos} for a more thorough treatment.

Recall that we define the decreasing rearrangement $f^*:\Rn{+}
\rightarrow \Rn{+} $ by
\begin{equation*}
 f^{*}(t) = \inf \{ s>0  \;|\;  \lambda_f(s) \le t \},
\end{equation*}
where
\begin{equation*}
 \lambda_f(s) = \abs{\{x \in \Omega \;|\; \abs{f(x)} > s   \}}.
\end{equation*}
For $1\le p,q \le \infty$ we define the Lorentz space
\begin{equation*}
L^{p,q}(\Omega) = \{ f \text{ measurable } \;|\;
\pqnormspace{f}{p}{q}{\Omega} < \infty \},
\end{equation*}
where
\begin{equation}\label{lorentz_def_1}
 \pqnormspace{f}{p}{q}{\Omega} =
 \left\{
 \begin{alignedat}{2}
    &\left( \int_0^{\infty} \left(t^{\frac{1}{p}} f^*(t)\right)^q \frac{dt}{t}   \right)^{\frac{1}{q}}
     &&\text{ for }q< \infty,\\
     &\sup_{t>0} t^{\frac{1}{p}} f^*(t) &&\text{ for }q= \infty.
 \end{alignedat}
 \right.
\end{equation}
Below we summarize some useful properties of Lorentz spaces;
proofs can be found in \cite{grafakos}.

\begin{lem}\label{lorentz_properties} The following hold.
\begin{enumerate}
 \item  For $1\le p,q \le \infty$, the spaces $L^{p,q}$ are quasi-Banach spaces, i.e. complete with respect to the quasi-norm \eqref{lorentz_def_1}.
 \item  The space $L^{p,p}$ coincides with the Lebesgue space $L^p$, and the space $L^{p,\infty}$ coincides with weak-$L^p$.
 \item  For $1\le p,q\le \infty$ the topology of $L^{p,q}$ generated by the quasi-norms is metrizable, and for $1<p<\infty,$ $1 \le q \le \infty$ also normable (see \eqref{ltinorm1} for the $p=2$, $q=\infty$ norm).
 \item For $1 \le p \le \infty$, $1\le q < r \le \infty$ there are constants $c_{p,q,r}$ such that
 \begin{equation}
  \pqnorm{f}{p}{r} \le c_{p,q,r} \pqnorm{f}{p}{q}.
 \end{equation}
This shows that the Lorentz spaces embed as the second index increases.
 \item  For $1 < p,q < \infty$, $(L^{p,q})^* = L^{p',q'}$ where $p'$ and $q'$ are the conjugate exponents of $p$ and $q$ respectively.  The duality is achieved via integration:
\begin{equation}
 \abs{\int_{\Omega} f(x) g(x) dx } \le C \pqnormspace{f}{p}{q}{\Omega} \pqnormspace{g}{p'}{q'}{\Omega}.
\end{equation}
\end{enumerate}
\end{lem}

Note in particular the embeddings $L^{2,q}(\Omega) \hookrightarrow
L^{2,2}(\Omega) = L^{2}(\Omega)$ for $1 \le q < 2$.  This suggests
defining our intermediate space  as follows: for any open set $V
\csubset \Rn{2}$ with $C^1$ boundary, we set
\begin{equation*}
\mathcal{X}(V) = \{ f \in H^1_0(V) \;|\; \nabla f \in L^{2,1}(V) \},
\end{equation*}
and endow it with the norm $\norm{f}_{\mathcal{X}(V)} = \pqnormspace{\nabla
f}{2}{1}{V}$, which makes $\mathcal{X}(V)$ into a Banach space.  Here we
abuse notation by writing $\pqnormspace{\nabla f}{2}{1}{V}$ for
the norm on $L^{2,1}(\Omega)$, which exists by item 4 of Lemma
\ref{lorentz_properties}, not the quasi-norm defined by
\eqref{lorentz_def_1}.  We write $\mathcal{X}^*(V) := (\mathcal{X}(V))^*$ for the dual
of $\mathcal{X}(V)$ and define the space
\begin{equation*}
 \mathcal{X}^*_{loc} = \mathcal{X}^*_{loc}(\Rn{2}) = \{ f \;|\; f \in \mathcal{X}^*(B(0,R))\;\; \forall R > 0    \}.
\end{equation*}
We say that a sequence $\{ f_n \} \subset \mathcal{X}^*_{loc}$ converges
locally-weak-$*$ in $\mathcal{X}^*_{loc}$ to $f$ if for every $V\csubset
\Rn{2}$, $f_n \overset{*}{\rightharpoonup} f$ in $\mathcal{X}^*(V)$.

It turns out that $\mathcal{X}(\Omega)$ has exactly the properties we
sought. Indeed, we have the following lemma, the proof of which
can be found in Theorem 3.3.4 of \cite{helein}.

\begin{lem}
 Let $V$ be an open subset of $\Rn{2}$ with $C^1$ boundary. Then $\mathcal{X}(V) \hookrightarrow C_c^0(V)$, and
 \begin{equation}
  \pnormspace{f}{\infty}{V} \le C \pqnormspace{\nabla f}{2}{1}{V}.
 \end{equation}
\end{lem}

In addition to being bounded and continuous, the functions in
$\mathcal{X}(V)$ are also differentiable almost everywhere.  See
\cite{stein_diff} for a proof of this fact.  However, it is not
possible to find a modulus of continuity for the functions in
$\mathcal{X}(V)$.

\begin{lem}
 There is no embedding $\mathcal{X}(V) \hookrightarrow C^{0,\omega}(V)$
  for any modulus of continuity $\omega$.
\end{lem}
\begin{proof}
 For simplicity we suppose $V=B(0,1)$.  Suppose there exists a modulus of continuity $\omega$ such that $\mathcal{X}(V) \hookrightarrow C^{0,\omega}(V)$, i.e.
 \begin{equation}
  \sup_{x\neq y} \frac{\abs{f(x)-f(y)}}{\omega(\abs{x-y})} \le C \pqnorm{\nabla
  f}{2}{1},
 \end{equation} with $\omega(s) \to 0 $ as $s \to 0$.
Then any bounded set in $\mathcal{X}(V)$ is equicontinuous, and hence, by
Arzela-Ascoli, pre-compact in $C^0(V)$.  It is easy to check that
$\mathcal{X}(V)$ is scale-invariant.  That is, $\pqnorm{\nabla f(\lambda
\cdot)}{2}{1} = \pqnorm{\nabla f}{2}{1}$ for all $\lambda>0$.  For
any function $f$ such that $f(0)\neq 0$, we consider
$\{f(n\cdot)\}_{n\in \mathbb{N}}$, which is pre-compact in
$C^0(V)$.  Since the support of $f(n\cdot)$ is contained in
$B(0,1/n)$, any convergent subsequences must converge uniformly to
$0$.  However, $f(n0) = f(0)\neq 0$, which contradicts the uniform
convergence to $0$.
\end{proof}

With these facts about the space $\mathcal{X}(\Omega)$ in hand we can prove a
compactness result for the Jacobians.

\begin{prop}\label{app_jac}
 Suppose that configurations
 $\{(u_\varepsilon,A_\varepsilon)\}$ satisfy (H1) -- (H4), and that $1\ll n$.
Further assume that $\pnorm{u_\varepsilon}{\infty} \le 1$.  Then up to extraction,
\begin{equation}\label{app_jac_0}
 \frac{\mu(u_\varepsilon,A_\varepsilon)}{2\pi n} \overset{*}{\rightharpoonup} \delta_p \;\text{ and} \;\;\;
 \frac{\mu(u_\varepsilon,A'_\varepsilon)}{2 \pi n} \overset{*}{\rightharpoonup} \delta_p
\end{equation}
weakly-$*$ in $\mathcal{X}^*(\Omega)$.  Define the blow-up measures on
$\Rn{2}$ at scale $\ell =\sqrt{n/h_{ex}}$ by
\begin{equation*}
 \tilde{\mu}(u,A)(x) = \ell^2 \vchi_{\Omega}(\ell x+p)  \mu(u,A) (\ell x+p).
\end{equation*}
Then the blow-up measures also converge up to extraction:
\begin{equation}\label{app_jac_00}
 \frac{\tilde{\mu}(u_\varepsilon,A_\varepsilon)}{2 \pi n} \overset{*}{\rightharpoonup} \mu_* \; \text{ and}\;\;\; \frac{\tilde{\mu}(u_\varepsilon,A'_\varepsilon)}{2 \pi n} \overset{*}{\rightharpoonup} \mu_*
\end{equation}
locally-weak-$*$ in  $\mathcal{X}^*_{loc}$, where $\mu_*$ is a probability measure on $\Rn{2}$.
\end{prop}

\begin{proof}
 We again suppress the subscript $\varepsilon$ on $u$, $A$, and $A'$ in calculations.  The pointwise bound $\abs{j_{A'}(u)} = \abs{(iu,\nabla_{A'} u) } \le \abs{u} \abs{\nabla_{A'} u} \le \abs{\nabla_{A'} u}$, together with the result of Theorem \ref{application_1}, shows that
\begin{equation}\label{app_jac_2}
  \pqnormspace{j'}{2}{\infty}{\Omega}
   \le  \pqnormspace{\nabla_{A'} u}{2}{\infty}{\Omega} \le C n.
\end{equation}
Invoking the $L^{2,1}$ -- $L^{2,\infty}$ duality and using
\eqref{app_jac_2} then proves that for $f \in \mathcal{X}(\Omega)$
\begin{equation}
\begin{split}
 \abs{\int_{\Omega} f \curl{j'}} &= \abs{\int_{\Omega} \nabla^{\bot}
 f \cdot j'}  \le \pqnormspace{\nabla f}{2}{1}{\Omega}
 \pqnormspace{j'}{2}{\infty}{\Omega} \\
 &\le C n \pqnormspace{\nabla f}{2}{1}{\Omega}.
\end{split}
\end{equation}
Theorem \ref{application_1} also showed that
the bound
\begin{equation}
 \pnormspace{\curl{A'}}{2}{\Omega} \le C n
\end{equation}
also holds.  This fact, combined with the Cauchy-Schwarz and
Poincar\'{e} inequalities, allows us to deduce the bound
\begin{equation}
 \abs{\int_{\Omega} f \curl{A'}} \le \pnormspace{f}{2}{\Omega}  \pnormspace{\curl{A'}}{2}{\Omega} \le C n \pqnormspace{\nabla f}{2}{1}{\Omega}.
\end{equation}

Thus, for any function $f \in \mathcal{X}(\Omega)$,
\begin{equation}\label{app_jac_5}
 \abs{\int_{\Omega} f \mu(u,A')} = \abs{\int_{\Omega} f \curl{(j' + A')}}
  \le  C n \pqnormspace{\nabla f}{2}{1}{\Omega}.
\end{equation}
This proves that the collection $ \left \{
\frac{1}{n}\mu(u,A') \right \}$
 is bounded in $\mathcal{X}^*(\Omega)$, the dual of $\mathcal{X}(\Omega)$.  Since $\mathcal{X}(\Omega)$ is separable, there exists a weak-$*$ sequential limit, and up to extraction
 \begin{equation}
  \frac{\mu(u,A')}{n} \overset{*}{\rightharpoonup} \nu
 \end{equation}
in $\mathcal{X}^*(\Omega)$.  Now, Proposition 9.5 of \cite{ss_book} shows
that up to extraction
\begin{equation*}
  \frac{\mu(u,A')}{n} \overset{*}{\rightharpoonup} 2 \pi \delta_p
\end{equation*}
weakly-$*$ in $(C_c^{0,\gamma}(\Omega))^*$ for $\gamma > 2/3$, and
hence we may conclude that $\nu = 2\pi \delta_p$.

Recall that $A' = A - h_{ex} \nabla^{\bot} \xi_0,$ which implies
that
\begin{equation}\label{app_jac_3}
j' + A'= j + A  - h_{ex}(1-\abs{u}^2)\nabla^{\bot} \xi_0.
\end{equation}
For any $f \in \mathcal{X}(\Omega)$ we estimate:
\begin{equation}\label{app_jac_4}
\begin{split}
& \abs{ h_{ex}\int_{\Omega} f \curl{((1-\abs{u}^2)\nabla^{\bot} \xi_0)} }
 = \abs{h_{ex} \int_{\Omega} \nabla^{\bot} f \cdot (1-\abs{u}^2)\nabla^{\bot} \xi_0} \\
& \le h_{ex} \varepsilon \left(\int_{\Omega} \frac{(1-\abs{u}^2)^2}{\varepsilon^2} \right)^{\frac{1}{2}}  \pnormspace{\nabla \xi_0}{\infty}{\Omega} \pnormspace{\nabla f}{2}{\Omega} \\
& \le h_{ex} \varepsilon \left( F_\varepsilon(\abs{u},\Omega)
\right)^{\frac{1}{2}} C \pqnormspace{\nabla f}{2}{1}{\Omega} \le C
\varepsilon^{\frac{1+\alpha}{2} -\beta} \pqnormspace{\nabla
f}{2}{1}{\Omega},
\end{split}
\end{equation}
where the first inequality follows from the embedding $L^{2,1}
\hookrightarrow L^2$ and the last follows from the assumptions
(H1).
%
Note also that we have absorbed $\pnormspace{\nabla
\xi_0}{2}{\Omega}$ into the constant since $\xi_0$ depends only on
the geometry of $\Omega$.  We conclude from \eqref{app_jac_3} and
\eqref{app_jac_4} that
\begin{equation}
 \norm{\frac{\mu(u,A)}{n}  - \frac{\mu(u,A')}{n} }_{\mathcal{X}^*(\Omega)} \rightarrow 0  \text{ as } \varepsilon \rightarrow 0.
\end{equation}
This proves \eqref{app_jac_0}.

It remains to prove \eqref{app_jac_00}.  First note that since
$\mu(u,A) \in L^1(\Omega)$, the blow-up, $\tilde{\mu}$, is an
element of $\mathcal{X}^*_{loc}$.  Fix $R>0$, and consider $\mathcal{X}(B(0,R))$.  We
will show that, up to extraction,
\begin{equation*}
 \frac{\tilde{\mu}(u,A')}{2 \pi n} \overset{*}{\rightharpoonup} \mu_*
\end{equation*}
weakly-$*$ in $\mathcal{X}^*(B(0,R))$, where $\mu_*$ is a probability measure.

Fix a function $f \in \mathcal{X}(B(0,R))$.  Recall that the blow-up of
$\mu$ is given by $\tilde{\mu}(x) = \ell^2 \mu(\ell x+p)
\vchi_{\Omega}(\ell x+p)$. By changing variables, we have
\begin{equation}
 \int_{B(0,R)} f(x) \tilde{\mu}(u,A')(x) dx = \int_{\Omega} f\left( \frac{x-p}{\ell}\right) \mu(u,A')(x)dx.
\end{equation}
For $\ell$ sufficiently small, i.e. for $\varepsilon$ sufficiently
small, the blow-down of the support of $f$ is contained in
$\Omega$, and hence $f((\cdot - p)/\ell) \in \mathcal{X}(\Omega)$.  This allows
us to apply the bound \eqref{app_jac_5} to conclude that
\begin{equation}
\begin{split}
 \abs{\int_{B(0,R)} f(x) \tilde{\mu}(u,A')(x) dx} & \le C n \pqnormspace{\nabla f((\cdot-p)/\ell)}{2}{1}{\Omega} \\ &= Cn \pqnormspace{\nabla f}{2}{1}{B(0,R)}.
\end{split}
\end{equation}
Then, as above, we conclude that up to extraction
\begin{equation}
 \frac{\tilde{\mu}(u,A')}{2 \pi n} \overset{*}{\rightharpoonup} \nu
\end{equation}
weakly-$*$ in $\mathcal{X}^*(B(0,R))$. Proposition 9.5 of \cite{ss_book}
shows that up to extraction
\begin{equation*}
  \frac{\tilde{\mu}(u,A')}{2\pi n}
  \overset{*}{\rightharpoonup} \mu_*
\end{equation*}
weakly-$*$ in $(C_c^{0,\gamma}(B(0,R)))^*$ for $\gamma > 2/3$,
where $\mu_*$ is a probability measure on $\Rn{2}$. This
proves that
\begin{equation}
 \frac{\tilde{\mu}(u,A')}{2\pi n}
\overset{*}{\rightharpoonup} \mu_*
\end{equation}
weakly-$*$ in $\mathcal{X}^*(B(0,R))$.  Applying \eqref{app_jac_4} to the
blow-up, we conclude that
\begin{equation}
 \frac{\tilde{\mu}(u,A)}{2\pi n} \overset{*}{\rightharpoonup} \mu_*
\end{equation}
weakly-$*$ in $\mathcal{X}^*(B(0,R))$ as well.  Since the above analysis
works for any choice of $R$ we conclude that the blow-up
convergence is locally-weak-$*$ convergence in $\mathcal{X}^*_{loc}$, i.e.
\eqref{app_jac_00} holds.

\end{proof}

\subsection{$L^{2,\infty}$ weak-$*$ convergence of $j'/n$}\label{current_convergence}

This section employs the stronger a priori bound
\begin{equation}
G_\varepsilon(u_\varepsilon,A_\varepsilon) \le f_\varepsilon(n) + n^2 I(\mu_*) + o(n^2)
\end{equation}
in addition to (H1) -- (H4) and $1 \ll n$.  These assumptions
allow us to employ item (2) of Theorem \ref{application_1} to find
more convergence results.  In particular, we will establish the
$L^{2,\infty}$ weak-$*$ convergence of the superconducting
currents.

We first show that $f(\abs{u})\abs{u}^2 X$ and $-\nabla^\bot G_p$ are close
in the weak-$*$ topology.

\begin{lem}\label{X_weakstar_conv}
 Let $X_\ep$ be the vector field defined by \eqref{X_def}.  Fix a vector field $H\in L^{2,1}(\Omega)$.  Then
 \begin{equation}
\abs{\int_\Omega (f(\abs{u})\abs{u}^2 X_\ep + \nabla^\bot
G_p)\cdot H} \le o(1) (\pqnormspace{H}{2}{1}{\Omega} +1) +
o_\delta(1).
 \end{equation}
\end{lem}

\begin{proof}
 Define the subsets $\Omega_1 = \Omega \backslash(B(p,\delta) \cup \mathcal{B})$
 and $\Omega_2 = B(p,\delta) \cup \mathcal{B}$, and drop the $\ep$ subscripts.
   Note that on the set
  $\Omega_1$ we have that $X = - f(\abs{u})\nabla^\bot G_p$, and so
 \begin{equation*}
 f(\abs{u}) \abs{u}^2 X + \nabla^\bot G_p = (1-f^2(\abs{u}) \abs{u}^2)\nabla^\bot G_p.
 \end{equation*}
Hence
 \begin{equation}
   \pqnormspace{f(\abs{u})\abs{u}^2 X + \nabla^\bot G_p}{2}
   {\infty}{\Omega_1} \le \pnormspace{1-f^2(\abs{u})\abs{u}^2}{\infty}{\Omega_1}
    \pqnormspace{\nabla G_p}{2}{\infty}{\Omega}  \le
    o(1).
 \end{equation}
Since $\abs{\Omega_2} \le \abs{\mathcal{B}} + \abs{B(p,\delta)} =  o(1)+ o_\delta(1),$ we have that
\begin{equation}
\pqnormspace{H}{2}{1}{\Omega_2} \le \int_0^{\abs{\Omega_2}} H^*(t) \frac{dt}{t^{1/2}} = o(1)+ o_\delta(1),
\end{equation}
where the equality follows from the absolute continuity of the integral.
Using these two bounds and Lemma \ref{X_control}, we then have that
\begin{equation}
\begin{split}
 & \abs{\int_\Omega (f(\abs{u})\abs{u}^2 X + \nabla^\bot G_p)\cdot H}  \\ &\le \abs{\int_{\Omega_1} (f(\abs{u})\abs{u}^2 X + \nabla^\bot G_p)\cdot H}
 + \abs{\int_{\Omega_2} (f(\abs{u})\abs{u}^2 X + \nabla^\bot G_p)\cdot H} \\
 & \le \pqnormspace{f(\abs{u})\abs{u}^2 X + \nabla^\bot G_p}{2}{\infty}{\Omega_1}
 \pqnormspace{H}{2}{1}{\Omega_1}  \\&+
 \pqnormspace{f(\abs{u})\abs{u}^2 X   + \nabla^\bot G_p}{2}{\infty}{\Omega_2}
  \pqnormspace{H}{2}{1}{\Omega_2} \\
 &\le o(1)\pqnormspace{H}{2}{1}{\Omega} + O(1)(o(1) + o_\delta(1))  = o(1)\pqnormspace{H}{2}{1}{\Omega} + o_\delta(1) + o(1).
\end{split}
\end{equation}

\end{proof}

\begin{remark} The above lemma also holds with the $f(\abs{u})$ terms
removed everywhere.  They are present in the lemma for ease of
use in what follows.  The reason the term is harmless is because
the field $X$ is only nonzero where $\abs{u} = 1+o(1)$, so adding
the $f$ term only modifies the powers of $\varepsilon$ that show
up in the $o(1)$ terms.
\end{remark}

This lemma allows us to deduce the convergence of the currents.  This is the content of the following proposition.

\begin{prop}\label{j_weakstar_conv}
Suppose (H1) -- (H4) hold, that $1\ll n$, and the a priori bound
\begin{equation*}
 G_\varepsilon(u_\varepsilon,A_\varepsilon) \le f_\varepsilon(n) + n^2 I(\mu_*) + o(n^2)
\end{equation*}
also holds.   Then for $j_\varepsilon' = (iu_\varepsilon, \nabla_{A_\varepsilon'} u_\varepsilon)$, we have that
 \begin{equation}
  \frac{f(\abs{u_\varepsilon})}{n} j_\varepsilon' \wstar -\nabla^\bot G_p
 \end{equation}
weakly-$*$ in $L^{2,\infty}(\Omega)$.  In particular, this implies that under the additional assumption that $\pnormspace{u_\varepsilon}{\infty}{\Omega} \le 1$, we have
\begin{equation}
 \frac{1}{n} j_\varepsilon' \wstar -\nabla^\bot G_p
\end{equation}
weakly-$*$ in $L^{2,\infty}(\Omega)$.
\end{prop}
\begin{proof}
 Again we suppress the subscript $\varepsilon$.  We have the pointwise bound
 \begin{equation}
 \begin{split}
   \abs{f(\abs{u})\frac{j'}{n} - f(\abs{u})\abs{u}^2 X } & = \abs{f(\abs{u})} \abs{(iu,\frac{\nabla_{A'} u}{n} - iuX)} \\
   & \le \abs{f(\abs{u}) \abs{u}} \abs{\frac{\nabla_{A'} u}{n} - iuX} \\
   & \le \abs{\frac{\nabla_{A'} u}{n} - iuX}
 \end{split}
 \end{equation}
since $xf(x) \le 1$.  Then the strong $L^{2,\infty}$ convergence $ \frac{\nabla_{A'} u}{n} - iuX \rightarrow 0$, given by Theorem \ref{application_1}, implies the strong $L^{2,\infty}$ convergence   $f(\abs{u}) j'/n - f(\abs{u})\abs{u}^2 X \rightarrow 0$.

We now prove the weak-$*$ convergence.  Let $H\in L^{2,1}(\Omega)$.  Then
\begin{multline}\label{j_w_c_1}
 \abs{\int_\Omega \left( \frac{f(\abs{u})j'}{n} + \nabla^\bot G_p  \right)\cdot H} \le
 \abs{\int_\Omega \left( \frac{f(\abs{u})j'}{n} - f(\abs{u})\abs{u}^2 X  \right)\cdot H} \\+
 \abs{\int_\Omega \left( f(\abs{u})\abs{u}^2 X + \nabla^\bot G_p  \right)\cdot H}
\end{multline}
From the above analysis, we know that
\begin{equation}\label{j_w_c_2}
 \abs{\int_\Omega \left( \frac{f(\abs{u})j'}{n} - f(\abs{u})\abs{u}^2 X
  \right)\cdot H}
  = o(1)\pqnormspace{H}{2}{1}{\Omega}.
\end{equation}
We then combine \eqref{j_w_c_1}, \eqref{j_w_c_2}, and Lemma \ref{X_weakstar_conv} to conclude that
\begin{equation}
 \abs{\int_\Omega \left( \frac{f(\abs{u})j'}{n} + \nabla^\bot G_p  \right)\cdot H} \le
 \pqnormspace{H}{2}{1}{\Omega} o(1) + o_\delta(1) + o(1).
\end{equation}
Let $\delta \rightarrow 0$; we conclude that $f(\abs{u}) j'/n  \wstar -\nabla^\bot G_p$ weakly-$*$ in $L^{2,\infty}(\Omega)$.  The second result follows by noting that $\abs{u} \le 1$ implies that $f(\abs{u}) = 1$.
\end{proof}

Together, Propositions \ref{app_jac} and \ref{j_weakstar_conv} demonstrate the convergence  $j'/n \wstar -\nabla^\bot G_p$ and $\mu'/n \wstar 2\pi \delta_p$ weakly-$*$ in $L^{2,\infty}(\Omega)$ and $\mathcal{X}^*(\Omega)$ respectively.  We also know from \eqref{app_1_0003} of Theorem \ref{application_1} that for $h' = \curl{A'}$, we have $h'/n \rightarrow G_p$ strongly in $L^2(\Omega)$.  We thus see the consistency between the relations
\begin{equation*}
 \curl{j'} + h' = \mu' \text{ and } 
\end{equation*}
and
\begin{equation*}
 \curl(-\nabla^\bot G_p) + G_p = 2\pi \delta_p.
\end{equation*}

We summarize below all the convergence results that hold in
addition to those of Proposition \ref{app_jac}.
\begin{cor}\label{london_convergence}
 Assume (H1) -- (H4)
 hold in addition to the assumptions $1\ll n$ and
 \begin{equation*}
  G_\varepsilon(u_\varepsilon,A_\varepsilon) \le f_\varepsilon(n) + n^2 I(\mu_*) + o(n^2).
 \end{equation*}
 Suppose that $\pnorm{u_\varepsilon}{\infty}\le 1$.  Then
\begin{equation}
 \left\{
 \begin{alignedat}{2}
  &\frac{1}{n}j_\varepsilon' \wstar -\nabla^\bot G_p && \text{ weakly-$*$} \text{ in }  L^{2,\infty}(\Omega) \\
 & \frac{1}{n} \mu_\varepsilon' \wstar 2\pi \delta_p  && \text{ weakly-$*$} \text{ in }  \mathcal{X}^*(\Omega) \\
 & \frac{1}{n}h_\varepsilon' \rightarrow G_p     && \text{ strongly} \text{ in }  L^2(\Omega).
 \end{alignedat}
 \right.
\end{equation}
\end{cor}

\section{Results in Lorentz-Zygmund spaces}\label{results_in_lz}

\subsection{Motivation}

We now return to the setting of Section \ref{current_convergence} in order to improve the convergence results from weak-$*$ $L^{2,\infty}$ to strong convergence in a slightly larger space.  In particular we assume that $1 \ll n$ and that the a priori upper bound
\begin{equation*}
 G_\varepsilon(u_\varepsilon,A_\varepsilon) \le f_\varepsilon(n) + n^2 I(\mu_*) +o(n^2)
\end{equation*}
holds.  To motivate what follows, we recall the mechanism that allowed us to prove the $L^{2,\infty}$ weak-$*$ convergence $j'/n \wstar -\nabla^\bot G_p$ in Proposition \ref{j_weakstar_conv}.  For simplicity we temporarily assume $\pnorm{u}{\infty} \le 1$ so that the normalization by $f(\abs{u})$ is not needed.  For $H\in L^{2,1}(\Omega)$, we bounded the integral
\begin{equation*}
 \abs{\int_\Omega \left( \frac{j'}{n} + \nabla^\bot G_p  \right)\cdot H} \le
 \abs{\int_\Omega \left( \frac{j'}{n} - \abs{u}^2 X  \right)\cdot H} +
 \abs{\int_\Omega \left( \abs{u}^2 X + \nabla^\bot G_p  \right)\cdot H}.
\end{equation*}
The first of these terms was bounded using the duality between $L^{2,1}$ and $L^{2,\infty}$ and the $L^2$ estimate $\pnorm{j'/n - \abs{u}^2 X}{2} = o(1).$  So, the first term is actually no obstacle to strong $L^{2,\infty}$ convergence.  On the other hand, we fail to control the second term by $o(1)\pqnorm{H}{2}{1}$, which would immediately give the strong convergence if it held.  Instead, we have to use the estimates of Lemma \ref{X_weakstar_conv}, which bound the second term by $o(1) \pqnorm{H}{2}{1} +o(1) + o_\delta(1)$.  These residual terms $o(1) + o_\delta(1)$ that do not multiply the $L^{2,1}$ norm of $H$ come from the product
\begin{equation}
 \pqnormspace{j'/n - \abs{u}^2 X}{2}{\infty}{B(p,K\ell)\cup \mathcal{B}} \pqnormspace{H}{2}{1}{B(p,K\ell)\cup \mathcal{B}}.
\end{equation}
The first of these quantities is bounded by a universal constant and the second vanishes because the measure  $\abs{B(p,K\ell)\cup \mathcal{B}} = o(1) + o_\delta(1)$.  It is therefore clear that the obstruction to strong $L^{2,\infty}$ convergence comes from the fact that we cannot prove estimates
$\pqnormspace{j'/n - \abs{u}^2 X}{2}{\infty}{B(p,K\ell)\cup \mathcal{B}}=o(1)$ on the small sets $B(p,K\ell)\cup \mathcal{B}$.

The problem is that the $L^{2,\infty}$ norm, like the $L^{\infty}$
norm, does not necessarily shrink to zero when it is calculated
over sets of vanishing measure. As a simple example consider the
function $f(x) = 1/\abs{x}$ in $\Rn{2}$.  A simple calculation
shows that
\begin{equation*}
\pqnorm{\vchi_{B(0,R)}f}{2}{\infty} = \sqrt{2} \pi.
\end{equation*}
for all $R>0$.  This points to a natural solution: we seek a space
that is slightly larger than $L^{2,\infty}(\Omega)$ with the
property that if $E_n$ is a sequence of sets with measure going to
zero, then $f\vchi_{E_n} \rightarrow 0$ strongly in the larger
space's norm for every $f\in L^{2,\infty}$.  Such spaces are found
in the Lorentz-Zygmund spaces.

\subsection{Lorentz-Zygmund spaces: definitions and properties}
The Lorentz-Zygmund spaces constitute a natural generalization of the Lorentz spaces $L^{p,q}$ and the Zygmund spaces $L^p\log^\alpha L = \{ f \;\vert\; \int(\abs{f}\log^\alpha(1+\abs{f}))^p< \infty\}$.   They are constructed by introducing the Zygmund space logarithmic weight with index $\alpha$ to the Lorentz spaces.  That is, for $1\le p,q \le \infty$, $\alpha \in \Rn{}$, we define the Lorentz-Zygmund space $L^{p,q}\log^\alpha L (\Omega)$ to be the collection of all measurable functions $f$ defined on $\Omega$ such that the quasi-norm
\begin{equation}\label{l_z_def}
 \lzspace{f}{p}{q}{\alpha}{\Omega} =
 \left\{
 \begin{alignedat}{2}
    &\left( \int_0^{\infty} \left(t^{\frac{1}{p}}\log^\alpha\left(e+\frac{1}{t}\right) f^*(t)\right)^q \frac{dt}{t}   \right)^{\frac{1}{q}}
     &&\text{ for }q< \infty,\\
     &\sup_{t>0} t^{\frac{1}{p}} \log^\alpha\left(e+\frac{1}{t}\right) f^*(t) &&\text{ for }q= \infty
 \end{alignedat}
 \right.
\end{equation}
is finite.  Here $f^*$ denotes the decreasing rearrangement of $f$; see the discussion of Lorentz spaces in Section \ref{cpt_jac} for the definition.

We summarize the crucial properties of these spaces in the following Lemma.  See \cite{ben_rud} for a thorough treatment of the spaces and the proofs.  Note that we use a slightly different logarithmic weight than is used in \cite{ben_rud}, but it makes no difference in the results.

\begin{lem}\label{l_z_properties}
 The following hold.
 \begin{enumerate}
\item For $1\le p,q\le \infty$, $\alpha \in \Rn{}$, the space $L^{p,q}\log^\alpha L (\Omega)$ is a quasi-Banach space, i.e. complete with respect to the quasi-norm \eqref{l_z_def}.
\item For $1 < p < \infty$, $1\le q \le \infty$, $\alpha \in \Rn{}$, the space $L^{p,q}\log^\alpha L (\Omega)$ is normable, i.e. there is a norm equivalent to the quasi-norm \eqref{l_z_def} that generates the same topology.
\item The space $L^{p,q}\log^0 L (\Omega)$ coincides with the Lorentz spaces $L^{p,q}(\Omega)$, and the space $L^{p,p}\log^\alpha L(\Omega)$ coincides with the Zygmund space $L^p \log^\alpha L(\Omega)$.
\item For $1\le p \le \infty$, $1\le q_0,q_1 \le \infty$, $\alpha,\beta \in \Rn{}$, we have the embedding
\begin{equation*}
 L^{p,q_1}\log^{\alpha}L (\Omega) \hookrightarrow  L^{p,q_0}\log^{\beta}L (\Omega)
\end{equation*}
when either $
 q_1 \le q_0 $  and $\alpha \ge \beta$
 or $  q_1 > q_0  $  and  $ \alpha + \frac{1}{q_1} > \beta +
 \frac{1}{q_0}$.
\item For $1 < p < \infty$, $1\le q \le \infty$, $\alpha\in \Rn{}$, we have that
\begin{equation*}
 \left(  L^{p,q}\log^{\alpha}L (\Omega)   \right)^* = L^{p',q'}\log^{-\alpha}L (\Omega),
\end{equation*}
where $p'$ and $q'$ are the conjugate exponents of $p$ and $q$ respectively.  The duality is achieved via integration:
\begin{equation*}
 \abs{\int_\Omega f g} \le \lzspace{f}{p}{q}{\alpha}{\Omega} \lzspace{g}{p'}{q'}{-\alpha}{\Omega}.
\end{equation*}
\end{enumerate}
\end{lem}

Now we must determine which Lorentz-Zygmund spaces have the property that we sought at the end of the last section.  The following lemma points the way.

\begin{lem}\label{desired_properties}
 Suppose a sequence of functions $f_n:\Omega \rightarrow \Rn{k}$, $k\ge 1$ is uniformly bounded in $L^{p,\infty}(\Omega)$, i.e. $\sup\limits_{n}\pqnormspace{f_n}{p}{\infty}{\Omega} \le C < \infty$.  Let $E_n \subset \Omega$ be a sequence of subsets so that $\abs{E_n} \rightarrow 0$.  If $1 \le q \le \infty$ and $\alpha < -1/q$, we have that
\begin{equation}
 \lzspace{\vchi_{E_n} f_n}{p}{q}{\alpha}{\Omega} \rightarrow 0.
\end{equation}
\end{lem}
\begin{proof}
The heart of the proof is the simple inequality
\begin{equation}\label{des_p_1}
 (\vchi_{E_n} f_n)^*(t) \le f_n^*(t) \vchi_{[0,\abs{E_n}]}(t).
\end{equation}
Then for $q=\infty$ we have
\begin{equation}
\begin{split}
 \lzspace{\vchi_{E_n} f_n}{p}{\infty}{\alpha}{\Omega} &= \sup_{t>0}  (\vchi_{E_n} f_n)^*(t) t^{1/p} \log^\alpha(e+1/t) \\
& \le \sup_{t>0}  f_n^*(t) \vchi_{[0,\abs{E_n}]}(t) t^{1/p} \log^\alpha(e+1/t) \\
& \le \sup_{t>0}  t^{1/p} f_n^*(t) \;\sup_{t>0}  \vchi_{[0,\abs{E_n}]}(t) \log^\alpha(e+1/t) \\
& \le \log^{\alpha}\left(e+ \frac{1}{\abs{E_n}}\right)  \sup_n \pqnormspace{f_n}{p}{\infty}{\Omega}  \\
& \le C \log^{\alpha}\left(e+ \frac{1}{\abs{E_n}}\right).
\end{split}
\end{equation}
Since $\alpha <0$ and $\abs{E_n} \rightarrow 0$, the conclusion follows.

Suppose now that $1\le q<\infty$.  Then, again using \eqref{des_p_1}, we have
\begin{equation}
\begin{split}
 \lzspace{\vchi_{E_n} f_n}{p}{q}{\alpha}{\Omega}^q  & = \int_0^{\infty} \left(t^{1/p}\log^\alpha\left(e+\frac{1}{t}\right) (\vchi_{E_n}f_n)^*(t)\right)^q \frac{dt}{t} \\
 & \le \left(\sup_{t>0} t^{1/p}f_n^*(t) \right)^q  \int_0^{\abs{E_n}}  \log^{q\alpha}\left(e+\frac{1}{t}\right) \frac{dt}{t} \\
 & \le C^q \int_{\log(e+1/\abs{E_n})}^\infty s^{q\alpha} \frac{e^s}{e^s-e}ds \\
 & \le C^q \frac{e}{e-1} \int_{\log(e+1/\abs{E_n})}^\infty s^{q\alpha}ds.
\end{split}
\end{equation}
Here we have used the change of variables $s = \log(e+1/t)$ for the second inequality.  Then, since $\alpha < -1/q$, we have that
\begin{equation}
 \int_{\log(e+1/\abs{E_n})}^\infty s^{q\alpha}ds = \frac{1}{-\alpha q-1} \log^{\alpha q +1}\left(e+\frac{1}{\abs{E_n}}\right) \rightarrow 0,
\end{equation}
from which the result follows.

\end{proof}

We must have $p=2$, so the above lemma tells us that our candidate spaces are $L^{2,q}\log^\alpha L(\Omega)$ with $\alpha < -1/q$.  We also see from this lemma that $L^{2,\infty}$ embeds into each of these spaces and that the function $x \mapsto 1/\abs{x}$ has finite norm in all of them.  However, item $4$ of Lemma \ref{l_z_properties} guarantees that $L^{2,\infty} \log^{\alpha} L(\Omega) \hookrightarrow L^{2,q} \log^{\beta} L(\Omega)$ for $\beta + 1/q < \alpha < 0$.  So, if we can prove the convergence results in $L^{2,\infty} \log^{\alpha} L(\Omega)$ for all $\alpha <0$, then this proves convergence in every one of the possible candidates.  In a sense, this says that the scale of spaces $L^{2,\infty}\log^{\alpha} L(\Omega)$, $\alpha <0$, is the smallest extension of $L^{2,\infty}(\Omega)$ in the Lorentz-Zygmund scale with the desired properties.

\subsection{Convergence results in Lorentz-Zygmund spaces}

With the motivation and definitions in place, we proceed to proving the convergence result.
For $\alpha \in \Rn{}$ we define the space
\begin{equation}
\mathcal{X}_\alpha(\Omega) = \{ f\in H^1_0(\Omega) \;\vert\; \nabla f \in L^{2,1}\log^\alpha L(\Omega) \},
\end{equation}
and endow it with the norm $\norm{f}_{\mathcal{X}_\alpha(\Omega)} = \lzspace{\nabla f}{2}{1}{\alpha}{\Omega}$, which makes $\mathcal{X}_\alpha(\Omega)$ into a Banach space.  For the purpose of calculations we will work with the quasi-norms that define the Lorentz-Zygmund spaces, but in defining the norm on $\mathcal{X}_\alpha(\Omega)$ we use the equivalent norm.  Note that $\mathcal{X}_0(\Omega) = \mathcal{X}(\Omega)$ as defined in Section \ref{cpt_jac}, and that $\mathcal{X}_\alpha(\Omega) \hookrightarrow \mathcal{X}_\beta(\Omega)$ for $\alpha \ge \beta$.  Write $\mathcal{X}_\alpha^*(\Omega) = (\mathcal{X}_{-\alpha}(\Omega))^*$ for the dual; we use the notation with the negative sign so that the scales match in the natural embedding $L^{2,\infty}\log^\alpha L(\Omega) \hookrightarrow \mathcal{X}^*_\alpha(\Omega)$.  Also define the space
\begin{equation*}
 \mathcal{X}^*_{\alpha, loc} = \mathcal{X}^*_{\alpha,loc}(\Rn{2}) = \{ f \;|\; f \in \mathcal{X}^*_{\alpha}(B(0,R))\;\; \forall R > 0    \}.
\end{equation*}
 
 We now prove that for any $\alpha <0$,  $j'/n \rightarrow -\nabla^\bot G_p$ strongly in $L^{2,\infty}\log^\alpha L(\Omega)$, and  $\mu'/n \rightarrow 2\pi \delta_p$ strongly in $\mathcal{X}^*_\alpha(\Omega)$, as well as the corresponding results for the blown-up vorticity. Recall that we define the function $f(x) = \vchi_{[0,1]}(x) + x^{-1}\vchi_{[1,\infty]}(x),$ and that we define the blow-up measures on $\Rn{2}$ at scale $\ell =\sqrt{n/h_{ex}}$ by
\begin{equation*}
 \tilde{\mu}(x) = \ell^2 \vchi_{\Omega}(\ell x+p)  \mu (\ell x+p).
\end{equation*}

\begin{prop}\label{l_z_convergence}
  Suppose configurations $\{(u_\varepsilon,A_\varepsilon)\}$
  satisfy the same  assumptions as in Proposition
  \ref{j_weakstar_conv}.
Then for any $\gamma <0$,
\begin{equation}\label{l_z_c_01}
\frac{1}{n}f(\abs{u_\varepsilon})j_\varepsilon' \rightarrow -\nabla^\bot G_p
\end{equation}
strongly in $L^{2,\infty}\log^\gamma L(\Omega)$.  If we further assume that $\pnormspace{u_\varepsilon}{\infty}{\Omega} \le 1$, then for any $\gamma <0$,
\begin{equation}\label{l_z_c_02}
\frac{1}{n}j_\varepsilon' \rightarrow -\nabla^\bot G_p
\end{equation}
strongly in $L^{2,\infty}\log^\gamma L(\Omega)$, and
\begin{equation}\label{l_z_c_03}
\frac{1}{n}\mu_\varepsilon' \rightarrow 2\pi \delta_p \;\text{ and } \;\frac{1}{n}\mu_\varepsilon \rightarrow 2\pi \delta_p
\end{equation}
strongly in $\mathcal{X}_\gamma^*(\Omega)$, where $\mu_\varepsilon =
\curl{j_\varepsilon + A_\varepsilon}$.  The blow-up measures also converge up to extraction:
\begin{equation}\label{l_z_c_04}
 \frac{\tilde{\mu}_\varepsilon'}{2 \pi n} \rightarrow \mu_* \; \text{ and}\;\;\; \frac{\tilde{\mu}_\varepsilon}{2 \pi n} \rightarrow \mu_*
\end{equation}
strongly in $\mathcal{X}^*_{\alpha,loc}(\Rn{2})$, where $\mu_*$ is a probability measure.

\end{prop}

\begin{proof}
We suppress the subscript $\varepsilon$.  Define the sets $\Omega_1 = \Omega \backslash (B(p,K\ell) \cup \mathcal{B})$, $\Omega_2 = B(p,K\ell) \cup \mathcal{B} $.  We trivially bound
\begin{multline}\label{l_z_c_0}
\lzspace{f(\abs{u})\frac{j'}{n} + \nabla^\bot G_p}{2}{\infty}{\gamma}{\Omega} \le \lzspace{f(\abs{u})\frac{j'}{n} + \nabla^\bot G_p}{2}{\infty}{\gamma}{\Omega_1}\\ +
\lzspace{f(\abs{u})\frac{j'}{n} + \nabla^\bot G_p}{2}{\infty}{\gamma}{\Omega_2}
\end{multline}
in order to treat each piece separately.  We deal with $\Omega_1 $ first.  With the vector field $X$ defined by \eqref{X_def}, we write
\begin{equation}\label{l_z_c_1}
f(\abs{u})\frac{j'}{n} + \nabla^\bot G_p = f(\abs{u})(iu,\frac{1}{n}\nabla_{A'} u - iuX) + f(\abs{u})\abs{u}^2(\nabla^\bot G_p+X) + \nabla^\bot G_p (1-\abs{u}^2f(\abs{u})).
\end{equation}
Since $xf(x) \le 1$ and $ 1-\varepsilon^{\alpha/4} \le  \abs{u} \le 1+\varepsilon^{\alpha/4}$ in $\Omega_1$, we may estimate
\begin{multline}\label{l_z_c_2}
\lzspace{f(\abs{u})\frac{j'}{n} + \nabla^\bot G_p}{2}{\infty}{\gamma}{\Omega_1} \le \lzspace{\frac{1}{n}\nabla_{A'} u - iuX}{2}{\infty}{\gamma}{\Omega_1}
\\ + (1+o(1))\lzspace{X+\nabla^\bot G_p}{2}{\infty}{\gamma}{\Omega_1}
+ o(1) \lzspace{\nabla^\bot G_p}{2}{\infty}{\gamma}{\Omega_1}.
\end{multline}
The first of these terms is $o(1)+o_{K,\delta}(1)$ by \eqref{app_1_0002} of Theorem \ref{application_1}, and the third is $o(1)$ since $\nabla^\bot G_p$ is in $L^{2,\infty}$.
The second term is $o(1) + o_\delta(1)$; to see this we employ Lemmas \ref{annular_Gp_swap} and \ref{desired_properties} and the fact that $\abs{\Omega_1 \backslash \supp(X)} =o(1)$ to estimate
\begin{equation}\label{l_z_c_10}
\begin{split}
\lzspace{X + \nabla^\bot G_p}{2}{\infty}{\gamma}{\Omega_1} & \le \pqnormspace{X + \nabla^\bot G_p}{2}{\infty}{\Omega_1 \cap \supp(X)}  \\& + \lzspace{\nabla^\bot G_p \vchi_{\Omega_1 \backslash \supp(X)}}{2}{\infty}{\gamma}{\Omega_1} \\
&= o(1)+o_\delta(1) + o(1) \pqnormspace{\nabla G_p}{2}{\infty}{\Omega_1}  
 = o(1) +o_\delta(1).
\end{split}
\end{equation}
Hence
\begin{equation}\label{l_z_c_5}
\lzspace{f(\abs{u})\frac{j'}{n} + \nabla^\bot G_p}{2}{\infty}{\gamma}{\Omega_1} = o(1) + o_{K,\delta}(1).
\end{equation}

We now turn to the $\Omega_2$ term.  Here we again utilize the crucial properties of the space $L^{2,\infty}\log^\gamma L(\Omega)$ that we proved in Lemma \ref{desired_properties}.  Indeed,   \eqref{app_1_0001} of Theorem \ref{application_1} and the fact that  $xf(x) \le 1$ guarantee that
\begin{equation}
\pqnormspace{f(\abs{u})\frac{j'}{n}}{2}{\infty}{\Omega} \le \pqnormspace{\frac{1}{n} \nabla_{A'} u}{2}{\infty}{\Omega} \le C_\Omega,
\end{equation}
where $C_\Omega$ is a constant that depends only on $\Omega$.  Applying Lemma \ref{desired_properties} then shows that
\begin{equation}\label{l_z_c_6}
\begin{split}
\lzspace{f(\abs{u})\frac{j'}{n} + \nabla^\bot G_p}{2}{\infty}{\gamma}{\Omega_2} &\le
\lzspace{\vchi_{\Omega_2} f(\abs{u})\frac{j'}{n}}{2}{\infty}{\gamma}{\Omega} +
\lzspace{ \vchi_{\Omega_2} \nabla^\bot G_p}{2}{\infty}{\gamma}{\Omega} \\
& \le (C_\Omega + \pqnormspace{\nabla G_p}{2}{\infty}{\Omega}) \log^\gamma\left( e+ \frac{1}{\abs{B(p,K\ell) \cup \mathcal{B}}}\right).
\end{split}
\end{equation}

Now, \eqref{l_z_c_0}, \eqref{l_z_c_5}, and \eqref{l_z_c_6} show that
\begin{equation}\label{l_z_c_7}
\lzspace{f(\abs{u})\frac{j'}{n} + \nabla^\bot G_p}{2}{\infty}{\gamma}{\Omega} \le
o(1) + o_{K,\delta}(1) + C\log^\gamma\left( e+ \frac{1}{\abs{B(p,K\ell) \cup \mathcal{B}}}\right),
\end{equation}
where $C = C_\Omega + \pqnormspace{\nabla G_p}{2}{\infty}{\Omega}$.  Let $\varepsilon \rightarrow 0$ and then let $K \rightarrow \infty$ and $\delta \rightarrow 0$.  Then, since $\gamma <0$, the right hand side of \eqref{l_z_c_7} goes to zero, and the strong convergence \eqref{l_z_c_01} is proved.

Suppose now that $\pnorm{u}{\infty} \le 1$.  Then $f(\abs{u}) =1$ everywhere, and so
\eqref{l_z_c_02} follows directly from \eqref{l_z_c_01}.  Let $g \in X_{-\gamma}(\Omega)$.  Then, since $-\Delta G_p + G_p = 2\pi \delta_p$, we have that
\begin{equation}\label{l_z_c_8}
\begin{split}
\abs{\int_\Omega g \frac{\mu'}{n} - 2\pi g(p)}  &= \abs{ \int_\Omega -\nabla^\bot g \cdot \left(\frac{j'}{n} +\nabla^\bot G_p \right)  + \int_\Omega g \left(\frac{\curl{A'}}{n} - G_p  \right)} \\
& \le \lzspace{\nabla g}{2}{1}{-\gamma}{\Omega} \lzspace{\frac{j'}{n} +\nabla^\bot G_p}{2}{\infty}{\gamma}{\Omega} \\
& + \pnormspace{g}{2}{\Omega} \pnormspace{\frac{\curl{A'}}{n} - G_p}{2}{\Omega} \\
& \le \norm{g}_{\mathcal{X}_{-\gamma}(\Omega)} \left(\lzspace{\frac{j'}{n} +\nabla^\bot G_p}{2}{\infty}{\gamma}{\Omega} +  \pnormspace{\frac{\curl{A'}}{n} - G_p}{2}{\Omega}  \right).
\end{split}
\end{equation}
Hence, by \eqref{l_z_c_02} and \eqref{app_1_0003}, we have that
\begin{equation}
\norm{\frac{\mu'}{n} - 2\pi \delta_p}_{\mathcal{X}^*_\gamma(\Omega)}  \rightarrow 0.
\end{equation}
An obvious modification of \eqref{app_jac_3} and \eqref{app_jac_4} shows that
\begin{equation}
\norm{\frac{\mu'}{n} - \frac{\mu}{n}}_{\mathcal{X}^*_\gamma(\Omega)} \rightarrow 0,
\end{equation}
 so we may conclude \eqref{l_z_c_03}.

To prove \eqref{l_z_c_04} we must use the set $B(p,K\ell)\backslash \mathcal{B}$ and blow up at scale $\ell$.  Indeed, from \eqref{app_1_0002} and the definition of the vector field $X$ there, we have that 
\begin{equation}
 \pnormspace{\frac{1}{n}\nabla_{A'}u - iuf(\abs{u}) \overline{\nabla^\bot U_*} }{2}{B(p,K\ell)\backslash \mathcal{B}}= o(1) + o_{K,\delta}(1),
\end{equation}
where $U_*$ solves $\Delta U_* = \mu_*$ in $B(0,K)$ and vanishes on $\partial B(0,K)$, and $\overline{\nabla^\bot U_*} $ is the blow-down at scale $\ell$ of $\nabla^\bot U_*$.  Arguing as above and writing $\tilde{j'}$ for the blow up of $j'$, we find that
\begin{equation}
 \lzspace{\frac{\tilde{j}'}{n} - \nabla^\bot U_*}{2}{\infty}{\gamma}{B(0,K)} \le
o(1) + o_{K,\delta}(1), 
\end{equation}
from which we deduce that $\tilde{\mu'}/n \rightarrow \mu_*$ in $\mathcal{X}^*_{\gamma,loc}(\Rn{2})$.
\end{proof}

\section{Results for solutions with $n$ bounded}\label{res_bnded}
Recall that in most of the results in Sections \ref{more_lower_bounds} -- \ref{results_in_lz} we have assumed that the vorticity mass diverges, i.e. $1 \ll n$.  This condition was needed to show the existence of weak limits after blow-up in Lemma \ref{leftover_lower_bound}, and these limits and their properties were crucial in proving most of the results in these sections.  Moreover, terms of the form $C/n$ and $(\log{n})/n$ were often written as $o(1)$, which certainly required the condition $1\ll n$ to hold.

In this section we examine the case of $n$ bounded.  The difficulties are two-fold.  First, without knowing the weak-limits after blow-up, it is not entirely clear what the correct vector field is to complete the square with in the region near $p$.  Second, to achieve lower bounds that match up to $o(1)$ the upper bounds for locally minimizing solutions with $n$ bounded, we need finer control on the lower bounds in the vortex balls.  In particular, we would need something like $F_\varepsilon(u,A',\mathcal{B}) \ge \pi n \log(r/\varepsilon) + n \gamma + o(1)$, where $\gamma$ is a specific constant related to the energy of a radial, degree-one vortex profile (see \eqref{u0_prop} below).  While it is possible to find such lower bounds by comparing the energy of a configuration to that of a local minimizer, there appears to be some difficulty in adapting our completion of the square technique to this setting.

We thus restrict our attention to the case of configurations $(u_\varepsilon,A_\varepsilon)$ that are solutions to the Ginzburg-Landau equations with $n$ independent of $\varepsilon$.  In particular, we will assume the solutions are of the type that we will consider in Section \ref{solns_apps}.  That is, $\{(u_\varepsilon,A_\varepsilon) \}$ satisfy the following assumptions.
\begin{enumerate}
\item[(J1)] $\{(u_\varepsilon,A_\varepsilon)\}$ are solutions satisfying $F_\varepsilon(u_\varepsilon,A_\varepsilon') \le C\log(1/\varepsilon)$ and $h_{ex} \le \varepsilon^{-\beta}$ for some $0<\beta <1$.
\item[(J2)] There exists an $R_0 >0$ and points $(a_1(\varepsilon),\dotsc,a_n(\varepsilon))\in \Omega^n$ such that $\abs{a_i(\varepsilon) - a_j(\varepsilon)} \gg \varepsilon$ for $i\neq j$, $d(a_i(\varepsilon),\partial \Omega) \gg \varepsilon$ for each $i$, and $\{\abs{u_\varepsilon} \le 1/2\} \subset \cup_i B(a_i(\varepsilon),R_0 \varepsilon)$.
\item[(J3)]  $\deg(u_\varepsilon,\partial B(a_i,R_0\varepsilon)) = 1$, and $u_\varepsilon$ has exactly one zero in each $B(a_i(\varepsilon),R_0 \varepsilon)$.
\item[(J4)] We have fixed the Coulomb gauge so that $\diverge{A_\varepsilon}=0$ in $\Omega$ and $ A_\varepsilon \cdot \nu = 0$ on $\partial \Omega$. 
\item[(J5)] The configurations satisfy the bounds 
\begin{equation*}
 \begin{split}
 &G_\varepsilon(u_\varepsilon,A_\varepsilon) \le f_\varepsilon(n) + B_0 n^2 \\
 &\abs{F_\varepsilon(u_\varepsilon,A_\varepsilon') - f^0_\varepsilon(n) } \le B_1 n^2,
\end{split}
\end{equation*}
where $B_0$ and $B_1$ are fixed positive constants that depend on $\Omega$, and $f_\varepsilon^0(n) = f_\varepsilon(n) - 2\pi n h_{ex} \underline{\xi_0} - h_{ex}^2 J_0$.

\end{enumerate}


We define the function $\Phi_\varepsilon$ to be the solution to
\begin{equation}\label{phi_def}
 \begin{cases}
  -\Delta \Phi_\varepsilon + \Phi_\varepsilon = 2\pi \sum_i \delta_{a_i} & \text{ in } \Omega \\
  \Phi_\varepsilon = 0 & \text{ on } \partial \Omega.
 \end{cases}
\end{equation}
We can write $\Phi_\varepsilon$ in a way similar to how we wrote $G_p$:
\begin{equation}\label{phi_exp}
 \Phi(x) = \sum_{i=1}^n (-\log{\abs{x-a_i}} + S_\Omega(x,a_i)), 
\end{equation}
where $S_\Omega(\cdot,\cdot) \in C^1(\Omega \times \Omega)$.  We now present a lemma that bounds the cross terms that appear when we complete the square with $\nabla^\bot \Phi_\varepsilon$ outside of the balls $\cup B(a_i(\varepsilon),R\varepsilon)$.  This result is essentially proved, up to a sign error, in Proposition 10.2 of \cite{ss_book}.  For clarity we present the proof here.

\begin{lem}\label{phi_cross_terms}
Suppose $\{(u_\varepsilon,A_\varepsilon)\}$ satisfy (J1) -- (J5).  Fix $R\ge R_0$, let $\Phi_\varepsilon$ be the function defined above, and define the set $\tilde{\Omega}:= \Omega \backslash (\cup_i B(a_i(\varepsilon),R\varepsilon))$.  Then 
\begin{multline} 
 \int_{\tilde{\Omega}}  \Phi_\varepsilon \curl{A_\varepsilon'} - \nabla^\bot \Phi_\varepsilon \cdot j_\varepsilon'  \ge \int_{\tilde{\Omega}} \abs{\nabla \Phi_\varepsilon}^2+ \abs{\Phi_\varepsilon}^2 \\ 
 -  \frac{C}{R} \left( \int_{\tilde{\Omega}} \frac{(1-\abs{u_\varepsilon}^2)^2}{\varepsilon^2} +\left( \int_{\tilde{\Omega}} \frac{(1-\abs{u_\varepsilon}^2)^2}{\varepsilon^2}\right)^{1/2} + \int_{\tilde{\Omega}}\abs{\nabla_{A_\varepsilon'} u_\varepsilon + i u_\varepsilon \nabla^\bot \Phi_\varepsilon}^2 \right) + o(1).
\end{multline}
\end{lem}
\begin{proof}
As usual, we suppress the subscript $\varepsilon$.  We begin by stating three bounds for the function $\Phi$ and its derivatives.  Using the expansion \eqref{phi_exp}, it can be shown that 
\begin{equation}\label{phi_bounds}
\begin{split}
& \abs{\Phi(x)} \le C \abs{\log{\abs{x-a_i}}} \text{ in } B(a_i,R\varepsilon) \\
& \pnormspace{\nabla \Phi}{\infty}{\tilde{\Omega}} \le C/(R\varepsilon) \\
& \pnormspace{\nabla \Phi}{4}{\tilde{\Omega}}^4 \le C/(R^2\varepsilon^2). 
\end{split}
\end{equation}
See Section 10.1 of \cite{ss_book}, for instance, for a proof of this fact. 
 
 We rewrite
\begin{equation}
\int_{\tilde{\Omega}}  \Phi \curl{A'} - \nabla^\bot \Phi \cdot j' = \int_{\tilde{\Omega}} \abs{\nabla \Phi}^2+ \abs{\Phi}^2 + \int_{\tilde{\Omega}}  \Phi (\curl{A'} - \Phi) - \nabla^\bot \Phi \cdot (j'+\nabla^\bot \Phi).
\end{equation}
 Writing $u = \rho e^{i\varphi}$, we will show that we can essentially set $\rho=1$ in $j' = \rho^2(\nabla \varphi - A')$.  Indeed,
\begin{equation}
 \int_{\tilde{\Omega}}  \Phi (\curl{A'} - \Phi) - \nabla^\bot \Phi \cdot (j'+\nabla^\bot \Phi) = I + II,
\end{equation}
where 
\begin{equation}
 I:= \int_{\tilde{\Omega}}  \Phi (\curl{A'} - \Phi) - \nabla^\bot \Phi \cdot (\nabla \varphi - A' +\nabla^\bot \Phi)
\end{equation}
and 
\begin{equation}
 II:= \int_{\tilde{\Omega}} (1-\rho^2) \nabla^\bot \Phi \cdot (\nabla \varphi - A' +\nabla^\bot \Phi) +(\rho^2-1)\abs{\nabla \Phi}^2.
\end{equation}

An application of Cauchy-Schwarz shows that 
\begin{multline}
  \abs{II} \le \varepsilon \pnormspace{\nabla \Phi}{\infty}{\tilde{\Omega}} \left( \int_{\tilde{\Omega}} \frac{(1-\rho^2)^2}{\varepsilon^2} + \int_{\tilde{\Omega}}\abs{\nabla \varphi - A' +\nabla^\bot \Phi}^2 \right) \\+ \varepsilon \pnormspace{\nabla \Phi}{4}{\tilde{\Omega}}^2 \left( \int_{\tilde{\Omega}} \frac{(1-\rho^2)^2}{\varepsilon^2}\right)^{1/2}. 
\end{multline}
Since $\rho \ge 1/2$ on $\tilde{\Omega}$, we have that 
\begin{equation*}
\abs{\nabla \varphi - A' +\nabla^\bot \Phi}^2 \le 4 \abs{\nabla_{A'} u + i u \nabla^\bot \Phi}^2
\end{equation*}
This bound and the bounds \eqref{phi_bounds} imply that
\begin{equation}\label{cpt_3}
\abs{II} \le  \frac{C}{R} \left( \int_{\tilde{\Omega}} \frac{(1-\rho^2)^2}{\varepsilon^2} +\left( \int_{\tilde{\Omega}} \frac{(1-\rho^2)^2}{\varepsilon^2}\right)^{1/2} + \int_{\tilde{\Omega}}\abs{\nabla_{A'} u + i u \nabla^\bot \Phi}^2 \right).
 \end{equation}

To handle $I$, we integrate by parts and use the fact that $-\Delta \Phi + \Phi = 2\pi \sum_i \delta_{a_i}$ to get
\begin{multline}\label{cpt_1}
   I =  \int_{\tilde{\Omega}}  \Phi (\curl{A'} - \Phi + \curl(\nabla \varphi - A') + \Delta \Phi )  - \int_{\partial \tilde{\Omega}}  \Phi(\nabla \varphi - A'  + \nabla^\bot \Phi)\cdot \tau \\ 
  =  - \int_{\partial \tilde{\Omega}}  \Phi(\nabla \varphi - A' + \nabla^\bot \Phi)\cdot \tau.
\end{multline}
Since $\Phi = 0$ on $\partial \Omega$, only the boundaries of the balls are important in $\partial \tilde{\Omega}$.  Then, writing $\bar{\Phi}_i$ for the average of $\Phi$ on $\partial B(a_i,R\varepsilon)$, we rewrite
\begin{multline}
 \int_{\partial B(a_i,R\varepsilon)}  \Phi(\nabla \varphi - A' + \nabla^\bot \Phi)\cdot \tau 
 = \int_{\partial B(a_i,R\varepsilon)}  (\Phi-\bar{\Phi}_i)(\nabla \varphi - A' + \nabla^\bot \Phi)\cdot \tau \\ + \bar{\Phi}_i \int_{\partial B(a_i,R\varepsilon)}  (\nabla \varphi - A' + \nabla^\bot \Phi)\cdot \tau.
\end{multline} 
We now argue as in (10.25) of \cite{ss_book} to bound $\abs{\tau \cdot (\nabla \varphi - A' + \nabla^\bot \Phi)} \le C/(R\varepsilon)$, and hence that 
\begin{equation}
 \abs{\int_{\partial B(a_i,R\varepsilon)}  (\Phi-\bar{\Phi}_i)(\nabla \varphi - A' + \nabla^\bot \Phi)\cdot \tau} \le \int_{\partial B(a_i,R\varepsilon)}  o(1) \frac{C}{R\varepsilon} = o(1).
\end{equation}
Also, 
\begin{equation}
\begin{split}
 \abs{\int_{\partial B(a_i,R\varepsilon)}  (\nabla \varphi - A' + \nabla^\bot \Phi)\cdot \tau} & =  \abs{2\pi + \int_{B(a_i,R\varepsilon)}\curl{A'} + \int_{\partial B(a_i,R\varepsilon)}\frac{\partial \Phi}{\partial \nu} }\\
 &= \abs{2\pi  - 2\pi +\int_{B(a_i,R\varepsilon)}\curl{A'}  + \Phi } \\
 & \le CR\varepsilon F_\varepsilon(u,A')^{1/2} + CR^2 \varepsilon^2 \abs{\log{R\varepsilon}},
\end{split}
\end{equation}
where the inequality follows from Cauchy-Schwarz for the $\curl{A'}$ term and \eqref{phi_bounds} for the $\Phi$ term.  Since $F_\varepsilon(u,A') \le C \abs{\log{\varepsilon}}$, and \eqref{phi_bounds} implies that $\abs{\bar{\Phi}_i} \le C \abs{\log{R \varepsilon}}$, we then have that
\begin{equation}\label{cpt_2}
 \abs{\bar{\Phi}_i \int_{\partial B(a_i,R\varepsilon)}  (\nabla \varphi - A' + \nabla^\bot \Phi)\cdot \tau} \le C \abs{\log{R\varepsilon}} (R\varepsilon \abs{\log{\varepsilon}}^{1/2} + R^2 \varepsilon^2 \abs{\log{R \varepsilon}}) =o(1).
\end{equation}
Then \eqref{cpt_1} -- \eqref{cpt_2} show that $I = o(1)$, which, along with \eqref{cpt_3} proves the result.
\end{proof}

Since we are dealing with solutions and each vortex ball has degree one, the natural candidate for the vector field to use in place of the weak limit of the blown-up currents is the perpendicular gradient of the unique radial, degree-one vortex solution in $\Rn{2}$.  We thus define the function $u_0:\Rn{2} \rightarrow \mathbb{C}$ to be the (unique radial) solution of $-\Delta u_0 = u_0(1-\abs{u_0}^2)$ in $\Rn{2}$.  Existence and uniqueness of a solution of the form $u_0 = f(r) e^{i\theta}$ in polar coordinates are established in \cite{hh}.  The fact that this $u_0$ is the unique degree-one solution was established in \cite{mironescu}.  It is known (see \cite{bbh} or Proposition 3.11 of \cite{ss_book}) that there exists a constant $\gamma>0$ such that
\begin{equation}\label{u0_prop}
 \frac{1}{2}\int_{B(0,R)} \abs{\nabla u_0}^2 + \frac{1}{2} (1-\abs{u_0}^2)^2 = \pi \log{R} + \gamma + o_R(1),
\end{equation}

We now present a result that provides a lower bound for the free energy of these solutions.  It is essentially the analogue of the results in Sections \ref{lb_split} -- \ref{lb_blowup}, which bounded the free energy in the case $1\ll n$.  In this case, the analogue of the vector field $X_\varepsilon$ is the vector field $Z_{\varepsilon,R}$, which is defined as follows.  For $R\ge R_0$, we define 
\begin{equation}\label{Z_def}
 Z_{\varepsilon,R} = \begin{cases}
      \varepsilon^{-1} \nabla u_0\left(\frac{\cdot-a_i}{\varepsilon}\right) & \text{in } B(a_i(\varepsilon),R\varepsilon), i=1,\dotsc,n \\
      -iu_\varepsilon \nabla^\bot \Phi_\varepsilon & \text{in } \Omega \backslash (\cup_i B(a_i(\varepsilon),R\varepsilon)),
     \end{cases}
\end{equation}
where $u_0$ is the radial one-vortex solution in $\Rn{2}$.

\begin{prop}\label{solns_fin_vort}
Assume configurations $\{(u_\varepsilon,A_\varepsilon)\}$ satisfy (J1) -- (J5).  Let $R\ge R_0$, and let $Z_{\varepsilon,R}$ be the vector field defined by \eqref{Z_def}.  Then for $\varepsilon$ sufficiently small and $R$ sufficiently large, 
 \begin{multline}\label{sfv_0}
 F_\varepsilon(u_\varepsilon,A_\varepsilon') \ge \frac{1}{4} \int_\Omega \abs{\nabla_{A_\varepsilon'} u_\varepsilon - Z_{\varepsilon,R}}^2 +\frac{1}{2} \int_\Omega \abs{\curl{A_\varepsilon'} - \Phi_\varepsilon  }^2 + \pi n \log{\frac{1}{\varepsilon}} + n \gamma 
 \\- \pi \sum_{i \neq j} \log{\abs{a_i - a_j}} + \pi \sum_{i,j} S_\Omega(a_i,a_j) +o(1) +o_R(1),
\end{multline}
where $o_R(1)$ vanishes as $R \rightarrow \infty$.
\end{prop}
\begin{proof}
We split the energy into into two components: that on $\tilde{\Omega}:= \Omega \backslash (\cup_i B(a_i,R\varepsilon))$, and that in the balls $\cup_i B(a_i,R\varepsilon)$.  We begin with the balls. In the ball $B(a_i,R\varepsilon)$ we complete the square with the blow-down of $u_0$ given by $\nabla v_i(x) := \varepsilon^{-1}\nabla u_0((x -a_i)/\varepsilon)$:
\begin{equation}
 \abs{\nabla_{A'} u }^2 = \abs{\nabla_{A'} u  - \nabla v_i}^2 + 2\Re(\nabla_{A'} u \cdot \nabla v_i) -\abs{\nabla v_i}^2.
\end{equation}
Now, from Proposition 3.12 of \cite{ss_book}, we know that, up to extraction, the blow-ups at scale $\varepsilon$ of $(u,A')$ converge to $(u_0,0)$ in $C^1_{loc}(\Rn{2})$.  We note that the vanishing limiting magnetic potential is a direct consequence of the Coulomb gauge condition (J4): the Coulomb gauge allows estimates of the $H^2$ norm of $A$, which in turn gives $L^\infty$ estimates.  So, making a blow-up change of variables, we have that 
\begin{equation}
 \int_{B(a_i,R\varepsilon)} 2\Re(\nabla_{A'} u \cdot \nabla v_i) -\abs{\nabla v_i}^2 = \int_{B(0,R)} \abs{\nabla u_0}^2 + o(1), 
\end{equation}
and 
\begin{equation}
  \int_{B(a_i,R\varepsilon)} \frac{1}{2\varepsilon^2} (1-\abs{u}^2)^2 = \int_{B(0,R)} \frac{1}{2} (1-\abs{u_0}^2)^2 +o(1).
\end{equation}
Combining these and summing over the $i$, we get that
\begin{multline}
 \frac{1}{2} \int_{\cup_i B(a_i,R\varepsilon)}\abs{\nabla_{A'} u }^2 +\frac{1}{2\varepsilon^2} (1-\abs{u}^2)^2 = \sum_i \frac{1}{2} \int_{ B(a_i,R\varepsilon)} \abs{\nabla_{A'} u  - \nabla v_i}^2 \\ + \frac{n}{2} \int_{B(0,R)} \abs{\nabla u_0}^2 + \frac{1}{2} (1-\abs{u_0}^2)^2 + o(1).
\end{multline}
Employing the property of $u_0$ given by \eqref{u0_prop}, we then have that 
\begin{multline}\label{sfv_3}
F_\varepsilon(u,A',\cup_i B(a_i,R\varepsilon)) = \sum_i \frac{1}{2} \int_{ B(a_i,R\varepsilon)} \abs{\nabla_{A'} u  - \nabla v_i}^2  +\abs{\curl{A'}}^2\\
+ \pi n \log{R} + n\gamma + o(1) +o_R(1).
\end{multline}

Outside of the balls, in $\tilde{\Omega}$, we complete the square with $iu\nabla^\bot \Phi$ as in Lemma \ref{j_square_completion} to get 
\begin{equation}
 \abs{\nabla_{A'} u }^2 = \abs{\nabla_{A'}u + iu \nabla^\bot \Phi}^2 - 2 \nabla^\bot \Phi \cdot j' - \abs{\nabla^\bot \Phi}^2\abs{u}^2. 
\end{equation}
We also complete the square with the $\curl{A'}$ term to get
\begin{equation}
 \abs{\curl{A'}}^2 = \abs{\curl{A'} - \Phi}^2 - \abs{\Phi}^2 + 2  \Phi \curl{A'}.
\end{equation}
From this we see that
\begin{multline}\label{sfv_1}
 \frac{1}{2}\int_{\tilde{\Omega}} \abs{\nabla_{A'} u}^2 + \abs{\curl{A'}}^2 = \frac{1}{2} \int_{\tilde{\Omega}} \abs{\nabla_{A'} u + i u \nabla^\bot \Phi}^2 + \abs{\curl{A'} - \Phi}^2 \\ +  \int_{\tilde{\Omega}}  \Phi \curl{A'} - \nabla^\bot \Phi \cdot j'   - \frac{1}{2} \int_{\tilde{\Omega}} \abs{\Phi}^2 + \abs{\nabla \Phi }^2 \abs{u}^2.
\end{multline} 
Combining \eqref{sfv_1} with Lemma \ref{phi_cross_terms} and using that $\abs{u} \le 1$,  we find that 
\begin{multline}\label{sfv_2}
 F_\varepsilon(u,A',\tilde{\Omega}) \ge \left(\frac{1}{2}-\frac{C}{R}\right) \int_{\tilde{\Omega}} \abs{\nabla_{A'} u + i u \nabla^\bot \Phi}^2 + \frac{1}{2} \int_{\tilde{\Omega}} \abs{\curl{A'} - \Phi}^2 + \int_{\tilde{\Omega}} \frac{1}{4\varepsilon^2}(1-\abs{u}^2)^2 \\
 -  \frac{C}{R} \left( \int_{\tilde{\Omega}} \frac{(1-\abs{u}^2)^2}{4 \varepsilon^2} +\left( \int_{\tilde{\Omega}} \frac{(1-\abs{u}^2)^2}{4 \varepsilon^2}\right)^{1/2} \right) 
 +  \frac{1}{2} \int_{\tilde{\Omega}}  \abs{\nabla \Phi }^2 + \abs{\Phi}^2 +o(1).
\end{multline} 

Now, for $x\ge 0$, we have that the minimum of $x - \frac{C(x+\sqrt{x})}{R}$ is $-C^2/(4R(R-C)) = o_R(1)$ as $R\rightarrow \infty$. We use this with $x = \int_{\tilde{\Omega}} \frac{1}{4\varepsilon^2}(1-\abs{u}^2)^2$ to replace the $(1-\abs{u}^2)^2$ integrals in \eqref{sfv_2} by $o_R(1)$.  For $R$ large enough we can also bound $1/2 -C/R \ge 1/4$.  To deal with the $\abs{\nabla \Phi }^2 + \abs{\Phi}^2$ term, we use an argument from Proposition 10.2 of \cite{ss_book} that uses the expansion \eqref{phi_exp} to show that
\begin{equation}\label{sfv_4}
 \frac{1}{2} \int_{\tilde{\Omega}} \abs{\nabla \Phi }^2 + \abs{\Phi}^2 = \pi n \log\frac{1}{R\varepsilon}- \pi \sum_{i \neq j} \log{\abs{a_i - a_j}} + \pi \sum_{i,j} S_\Omega(a_i,a_j) +o(1).
\end{equation}
Thus, for $R$ sufficiently large, 
\begin{multline}\label{sfv_5}
 F_\varepsilon(u,A',\tilde{\Omega}) \ge \frac{1}{4} \int_{\tilde{\Omega}} \abs{\nabla_{A'} u + i u \nabla^\bot \Phi}^2 + \frac{1}{2} \int_{\tilde{\Omega}} \abs{\curl{A'} - \Phi}^2  + \pi n \log\frac{1}{R\varepsilon}\\ - \pi \sum_{i \neq j} \log{\abs{a_i - a_j}} + \pi \sum_{i,j} S_\Omega(a_i,a_j) +o(1) +o_R(1).
\end{multline} 
Adding \eqref{sfv_5} to \eqref{sfv_3} and using the fact that $\int_{\cup_i B(a_i,R\varepsilon)} \abs{\Phi} = o(1)$ then proves \eqref{sfv_0}.

\end{proof}

\begin{remark}
The condition (J2) guarantees that making $R$ large does not lead to the possibility of $B(a_i(\varepsilon),R\varepsilon) \cap B(a_j(\varepsilon),R\varepsilon) \neq \varnothing$ for $i\neq j$.  Thus we are free to eventually let $R \rightarrow \infty$.
\end{remark}

Eventually we will need an estimate of $\pqnorm{Z_{\varepsilon,R}}{2}{\infty}$.  We prove this now in the analogue of Lemma \ref{X_control}.

\begin{lem}\label{Z_control}
Assume configurations $\{(u_\varepsilon,A_\varepsilon)\}$ satisfy (J1) -- (J5).  There is a constant $C_\Omega >0$, depending only on $\Omega$, such that if $h_{ex} \ge C_\Omega n^2$, then
for $\varepsilon$ sufficiently small, 
\begin{equation}\label{zc_0}
 C_0 n - o(1) -o_R(1) \le \pqnormspace{Z_{\varepsilon,R}}{2}{\infty}{\Omega} \le C_1 n +o(1) +o_R(1),
\end{equation}
where $C_0$ is a positive universal constant and $C_1$ is a positive constant that depends only on $\Omega$.
\end{lem}
\begin{proof}
 Arguing as in Lemma \ref{X_control} and employing the bounds $1 \ge \abs{u} \ge 1/2$ in $\cup_i B(a_i,R\varepsilon)$, we see that 
\begin{multline}\label{zc_01}
 \frac{1}{2}\pqnormspace{\nabla \Phi_\varepsilon}{2}{\infty}{\Omega} - o(1) -o_R(1) \le \pqnormspace{Z_{\varepsilon,R}}{2}{\infty}{\Omega} \\ 
 \le \pqnormspace{\nabla \Phi_\varepsilon}{2}{\infty}{\Omega} + n \pqnormspace{\nabla u_0}{2}{\infty}{\Rn{2}} + o(1) +o_R(1).
\end{multline}
We claim that if $h_{ex}\ge C_\Omega n^2$, then $C_0 n \le \pqnormspace{\nabla \Phi_\varepsilon}{2}{\infty}{\Omega} \le C_1 n$, where $C_0$ is a positive universal constant and $C_1$ depends only on $\Omega$.  Once the claim is established, \eqref{zc_0} follows immediately from \eqref{zc_01}.

We begin the proof of the claim by showing that the energy bounds (J5) imply that the points $a_i$ are a distance from $p$ controlled by $1/\sqrt{h_{ex}}$.  Using a modification of the energy splitting result, Proposition \ref{app_balls}, found in (11.32) -- (11.33) of \cite{ss_book}, we have that
\begin{equation}
G_\varepsilon(u,A) = h_{ex}^2 J_0 + F_\varepsilon(u,A') + 2\pi h_{ex}\sum_i \xi_0(a_i) +o(1).
\end{equation}   
Plugging in the upper bound of $G_\varepsilon$ and the lower bound of $F_\varepsilon$ given by (J5) and recalling that $\underline{\xi_0} = \xi_0(p)$, we find  
\begin{equation}
B_0 n^2 \ge 2\pi h_{ex} \sum_i (\xi_0(a_i) - \xi_0(p)) - B_1 n^2 +o(1). 
\end{equation}
From this we conclude that for $\varepsilon$ sufficiently small, 
\begin{equation}\label{zc_1}
\sup_i  (\xi_0(a_i) - \xi_0(p)) \le \frac{(B_0+B_1)n^2 }{4\pi h_{ex}}. 
\end{equation}
Now an analysis of the function $\xi_0$ will allow us to pass from the bound of \eqref{zc_1} to a bound on $\sup_i \abs{a_i -p}$. Recall that we have assumed that $\xi_0$ achieves its unique minimum at $p$ and that $D^2 \xi_0(p)$ is positive definite. It can be shown (see \cite{ss_glmin}) that the set of critical points of $\xi_0$ is finite.  Using these facts with Taylor's theorem, we may conclude that if $h_{ex} \ge C_{\Omega} n^2$, then 
\begin{equation}
\sup_i \abs{a_i -p} \le C_{\Omega} \frac{n}{\sqrt{h_{ex}}} = :R_{h_{ex}},
\end{equation} 
where $C_{\Omega}$ is a constant that depends on $\Omega$ (via dependence on the smaller eigenvalue of $D^2 \xi_0(p)$, $\pnorm{D^3 \xi_0}{\infty}$, etc).



This concentration of vortices inside $B(p,R_{h_{ex}})$ allows us to obtain a lower bound on $\pqnorm{\nabla \Phi_\varepsilon}{2}{\infty}$ of order $n$.  To see this, assume that $h_{ex}$ is sufficiently large so that $R_{h_{ex}}\le \dist(p,\partial\Omega)/2$, and fix any $r\in (R_{h_{ex}},2 R_{h_{ex}})$.  Note that the bounds on $R_{h_{ex}}$ imply that $B(p,2R_{h_{ex}}) \subseteq \Omega$.  Then, since each $a_i \in B(p,R_{h_{ex}})$, we have that
\begin{equation}\label{zc_40}
2 \pi n = \int_{B(p,R_{h_{ex}})} -\Delta \Phi_\varepsilon + \Phi_\varepsilon = \int_{B(p,r)} -\Delta \Phi_\varepsilon + \Phi_\varepsilon \le  \int_{B(p,r)} \abs{\Phi_\varepsilon} +  \int_{\partial B(p,r)} \abs{\nabla \Phi_\varepsilon} .
\end{equation}
Recalling the definition of the $L^{2,\infty}$ norm defined in \eqref{ltinorm1}, we may bound
\begin{equation}
 \int_{B(p,r)} \abs{\Phi_\varepsilon}  \le \sqrt{\pi} r \pqnorm{\Phi_\varepsilon}{2}{\infty}.
\end{equation}
Plugging this into \eqref{zc_40} and integrating over $(R_{h_{ex}},2 R_{h_{ex}})$ we see that
\begin{equation}
 2 \pi n R_{h_{ex}} \le \frac{\sqrt{\pi}}{2} 3R_{h_{ex}}^2 \pqnorm{\Phi_\varepsilon}{2}{\infty} + \int_{B(p,2R_{h_{ex}})\backslash B(p,R_{h_{ex}})} \abs{\nabla \Phi_\varepsilon}
\end{equation}
Dividing by the square roof of the area of the annulus $B(p,2R_{h_{ex}})\backslash B(p,R_{h_{ex}})$ and again using \eqref{ltinorm1}, we get
\begin{equation}\label{zc_4}
\frac{2\sqrt{\pi} n}{\sqrt{3}} \le \frac{\sqrt{3} R_{h_{ex}}}{2}\pqnorm{\Phi_\varepsilon}{2}{\infty} + \pqnorm{\nabla \Phi_\varepsilon}{2}{\infty}.
\end{equation}
To estimate $\pqnorm{\Phi_\varepsilon}{2}{\infty}$, we use the expansion \eqref{phi_exp} and the constant $B_\Omega$ defined by
\begin{equation}
B_\Omega : = \sup_{y\in \Omega} \pqnormspace{-\log\abs{\cdot - y} + S_\Omega(\cdot,y)}{2}{\infty}{\Omega} < \infty.
\end{equation}  
We find that
\begin{equation}
\pqnorm{\Phi_\varepsilon}{2}{\infty} \le \sum_{i=1}^n \pqnorm{-\log\abs{\cdot - a_i} + S_\Omega(\cdot,a_i)}{2}{\infty} \le n B_\Omega.
\end{equation}
Note that if $h_{ex} \ge C_\Omega n^2$, for some constant depending only on $\Omega$, then \begin{equation*}
\frac{\sqrt{3}B_\Omega R_{h_{ex}}}{2} \le \sqrt{\frac{\pi}{3}},     
\end{equation*}
and we conclude the lower bound 
\begin{equation}\label{zc_5}
\frac{\sqrt{\pi} n}{\sqrt{3}} \le \pqnorm{\nabla \Phi_\varepsilon}{2}{\infty}.
\end{equation}

The upper bound is far easier to prove.  Indeed, we use the expansion \eqref{phi_exp} to calculate
\begin{equation}
\pqnorm{\nabla \Phi_\varepsilon}{2}{\infty} \le \sum_{i=1}^n \pqnorm{-(\cdot-a_i)/\abs{\cdot - a_i}^2 + \nabla S_\Omega(\cdot,a_i)}{2}{\infty}  \le C_1 n,
\end{equation}
for some positive $C_1$ that depends on $\Omega$.

\end{proof}

With these results in hand, we can prove the analogue of Theorems \ref{squeeze_lower} and \ref{application_1} for the case of $n$ independent of $\varepsilon$ and $h_{ex}$ either bounded or divergent.  First we introduce a bit of notation. In the case $h_{ex} = O(1)$ we define the renormalized energy $R_{n,h_{ex}}: \Omega^n \rightarrow \Rn{}$ by 
\begin{equation}\label{ren_en_1}
 R_{n,h_{ex}}(x_1,\dotsc,x_n) = -\pi \sum_{i\neq j} \log{\abs{x_i-x_j}} + \pi \sum_{i,j} S_{\Omega}(x_i,x_j) + 2\pi h_{ex}\sum_i \xi_0(x_i).
\end{equation} 
In the case $1 \ll h_{ex}$ we define the renormalized energy $w_n:(\Rn{2})^n \rightarrow \Rn{}$ by
\begin{equation}\label{ren_en_2}
 w_n(x_1,\dotsc,x_n) = -\pi \sum_{i\neq j} \log{\abs{x_i - x_j}} + \pi n \sum_i Q(x_i),
\end{equation}
where $Q$ is the quadratic form of $D^2\xi_0(p)$.  Notice that $w_n$ is defined on $(\Rn{2})^n$ and not on $\Omega^n$; this is because $w_n$ is applied after blowing up at scale $\ell = \sqrt{n/h_{ex}}$, which is $o(1)$ when $1 \ll h_{ex}$.  In particular it is applied to the points $\tilde{a}_i = (a_i -p)/\ell$.

\begin{prop}\label{bounded_n_thm}
 Assume configurations $\{(u_\varepsilon,A_\varepsilon)\}$ satisfy (J1) -- (J5).

1. Suppose $h_{ex} = O(1)$.  Then for $\varepsilon$ sufficiently small and $R$ sufficiently large, 
\begin{multline}\label{bnt_0}
G_\varepsilon(u_\varepsilon,A_\varepsilon) \ge \frac{1}{4} \int_\Omega \abs{\nabla_{A_\varepsilon'} u_\varepsilon - Z_{\varepsilon,R}}^2 +\frac{1}{2} \int_\Omega \abs{\curl{A_\varepsilon'} - \Phi_\varepsilon }^2 \\
+ h_{ex}^2 J_0 + \pi n \log{\frac{1}{\varepsilon}} + \min_{\Omega^n} R_{n,h_{ex}}  + n \gamma + o(1) + o_R(1).
\end{multline} 
We always have that the bounds 
\begin{equation}\label{bnt_11}
 C_0\sqrt{n} \le  \pqnormspace{\nabla_{A_\varepsilon'} u_\varepsilon}{2}{\infty}{\Omega} \le C_1 n 
\end{equation}
hold, where $C_0$ and $C_1$ are positive constants.
Moreover, if the solutions satisfy the upper bound  
\begin{equation}\label{bnt_1}
G_\varepsilon(u_\varepsilon,A_\varepsilon) \le 
h_{ex}^2 J_0 + \pi n \log{\frac{1}{\varepsilon}} + \min_{\Omega^n} R_{n,h_{ex}}  + n \gamma + o(1),
\end{equation}
and $h_{ex} \ge C_\Omega n^2$ (the constant from Lemma \ref{Z_control}), then 
\begin{equation}\label{bnt_2}
 C_0 n \le \pqnormspace{\nabla_{A_\varepsilon'} u_\varepsilon}{2}{\infty}{\Omega} \le C_1n, 
\end{equation}
where $C_0$ is a universal positive constant and $C_1$ depends only on $\Omega$.

2. Suppose $1 \ll h_{ex}$.  Then for $\varepsilon$ sufficiently small and $R$ sufficiently large, 
\begin{multline}\label{bnt_3}
G_\varepsilon(u_\varepsilon,A_\varepsilon) \ge \frac{1}{4} \int_\Omega \abs{\nabla_{A_\varepsilon'} u_\varepsilon - Z_{\varepsilon,R}}^2 +\frac{1}{2} \int_\Omega \abs{\curl{A_\varepsilon'} - \Phi_\varepsilon  }^2 \\
+f_\varepsilon(n) + \min_{(\Rn{2})^n} w_n  + n \gamma + o(1) + o_R(1).
\end{multline} 
Moreover, if the solutions satisfy the upper bound  
\begin{equation}\label{bnt_4}
G_\varepsilon(u_\varepsilon,A_\varepsilon) \le 
f_\varepsilon(n) + \min_{(\Rn{2})^n} w_n  + n \gamma + o(1),
\end{equation}
then 
\begin{equation}\label{bnt_5}
 C_0 n \le  \pqnormspace{\nabla_{A_\varepsilon'} u_\varepsilon}{2}{\infty}{\Omega} \le C_1 n, 
\end{equation}
where $C_0$ is a universal positive constant and $C_1$ depends only on $\Omega$.
\end{prop}

\begin{proof}
The inequalities \eqref{bnt_11} follow from (J1) -- (J5) and Proposition \ref{degree_control}.
The proof of the rest is very similar to the proofs of Theorems \ref{squeeze_lower} and \ref{application_1}; here we use a slightly different form of the energy splitting lemma and we use the free energy bounds of Proposition \ref{solns_fin_vort}.  Indeed, Lemma 7.3  and bounds  (11.32) -- (11.33) of \cite{ss_book} show that 
\begin{equation}
G_\varepsilon(u,A) = h_{ex}^2 J_0 + F_\varepsilon(u,A',\Omega) + 2\pi h_{ex}\sum_i \xi_0(a_i) +o(1).
\end{equation}   
In the case $h_{ex} = O(1)$ we then insert Proposition \ref{solns_fin_vort} into this to get \eqref{bnt_0}.   In the case $1 \ll h_{ex}$, the above and a blow-up at scale $\ell$ (see the arguments following (11.33) in \cite{ss_book}) show \eqref{bnt_3}.
In either case, we may compare the matching upper and lower bounds to show that 
\begin{equation}
 \pnormspace{\nabla_{A'} u - Z_{\varepsilon,R}}{2}{\Omega} = o(1) +o_R(1).
\end{equation}
Then \eqref{bnt_2} and \eqref{bnt_5} follow from this and Lemma \ref{Z_control}.

\end{proof}

\section{Stable solutions}\label{solns_apps}
In this section we apply Theorem \ref{application_1} and Proposition \ref{bounded_n_thm} to the branches of stable solutions constructed in Chapter 11 of \cite{ss_book} and to energy 
minimizers in the regime $1 \ll n(\varepsilon) \ll h_{ex}(\varepsilon) \le C \lep$ constructed there in
Chapter 9.

The stable solutions are constructed to have a prescribed number of vortices $n$.  As before, the typical inter-vortex distance scale is given by $\ell = \sqrt{n/h_{ex}}$.  We say that $h_{ex}(\varepsilon)$ and $n(\varepsilon)$ are admissible if they satisfy the following two conditions.\\
1. There exists $\beta_0< 1/2$ such that $h_{ex} < \varepsilon^{-\beta_0}$.\\
2. If $n \neq 0$, then $n^2 \le \eta h_{ex}$, and $n^2
\log{\frac{1}{\ell}} \le \eta \log{\frac{\ell}{\varepsilon}}$.  Here
$\eta$ is a small parameter depending on $\Omega$ and $\beta_0$.

The behavior of the quantity $\ell$ is separated into three distinct
cases, and each case produces solutions with different
asymptotics.  The first case assumes $\ell$ does not tend to zero,
and the admissibility conditions then ensure that $n$ and $h_{ex}$
are both bounded.  Up to extraction we then assume that $n$ is
independent of $\varepsilon$.  The second case lets $\ell$ go to zero
but assumes that $n$ stays bounded.  In the third case $\ell$ goes to
zero and $n$ diverges to infinity.

The following theorem includes  the new $\lti$ bounds in the
results of \cite{ss_book}.

\begin{thm}{(\cite{ss_book} Theorem 11.1 Redux)}\label{stable_solns}

Given $\beta_0 \in (0,1/2)$, taking $\eta =
\eta(\Omega,\beta_0)$ sufficiently small, and given
$n(\varepsilon)$ and $h_{ex}(\varepsilon)$ admissible, there
exists $\varepsilon_0 >0 $ such that for $0 < \varepsilon <
\varepsilon_0$ there exists a configuration
$(u_\varepsilon,A_\varepsilon)$ with the following properties.

The configuration $(u_\varepsilon,A_\varepsilon)$ is a locally
minimizing critical point of $G_\varepsilon$ and hence a stable
solution of the Ginzburg-Landau equations.  The function
$u_\varepsilon$ has exactly $n$ zeroes, located at points
$a_1(\varepsilon),\dotsc,a_n(\varepsilon) \in \Omega$, and there
exists $R>0$ such that $\abs{u_\varepsilon} \ge 1/2$ on the set
$\Omega \backslash \cup_i B(a_i(\varepsilon),R\varepsilon)$ and
$\deg{(u_\varepsilon,\partial B(a_i(\varepsilon),R\varepsilon))} =
1$.  Finally, depending on which of the three cases described
above holds, we have one of the following.

1. (Case 1) If $\ell$ does not tend to zero so that $n$ is
independent of $\varepsilon$ and $h_{ex}$ is bounded independently
of $\varepsilon$, then up to extraction the n-tuple
$(a_1(\varepsilon),\dotsc,a_n(\varepsilon))$ converges as
$\varepsilon \rightarrow 0$ to a minimizer of $R_{n,h_{ex}}$,
which was defined in \eqref{ren_en_1}.
The energy of these solutions as $\varepsilon
\rightarrow 0$ is given asymptotically by
\begin{equation*}
 G_\varepsilon(u_\varepsilon,A_\varepsilon) = h_{ex}^2 J_0 + \pi n \abs{\log{\varepsilon}} + \min_{\Omega^n}{R_{n,h_{ex}}} + n \gamma + o(1),
\end{equation*}
where $\gamma>0$ is the constant from \eqref{u0_prop}.  We always have that the bounds 
\begin{equation}
 C_0\sqrt{n} \le  \pqnormspace{\nabla_{A_\varepsilon'} u_\varepsilon}{2}{\infty}{\Omega} \le C_1 n 
\end{equation}
hold, where $C_0$ and $C_1$ are positive constants.  Moreover, if $h_{ex} \ge C_\Omega n^2$ (the constant from Lemma \ref{Z_control}), then 
\begin{equation}
 C_0 n \le \pqnormspace{\nabla_{A_\varepsilon'} u_\varepsilon}{2}{\infty}{\Omega} \le C_1n, 
\end{equation}
where $C_0$ is a positive universal constant and $C_1$ depends only on $\Omega$.

2. (Case 2) If $n$ is independent of $\varepsilon$ and $h_{ex}
\rightarrow \infty$ then, up to extraction the rescaled n-tuple
$(\tilde{a}_1(\varepsilon),\dotsc,\tilde{a}_n(\varepsilon))$, where
$\tilde{a}_i(\varepsilon) = (a_i(\varepsilon) -p)/\ell$, converges as
$\varepsilon \rightarrow 0 $ to a minimizer of $w_n$, which was defined in \eqref{ren_en_2}.
The energy of these solutions as $\varepsilon
\rightarrow 0$ is given asymptotically by
\begin{equation*}
 G_\varepsilon(u_\varepsilon,A_\varepsilon) = f_\varepsilon(n) + \min_{(\Rn{2})^n}{w_n} + n\gamma + o(1).
\end{equation*}
Finally, we have
\begin{equation}\label{stable_solns_12}
C_0 \le \frac{1}{n} \pqnormspace{\nabla_{A_\varepsilon'} u_\varepsilon}{2}{\infty}{\Omega} \le C_1,
\end{equation}
where $C_0$ is a positive universal constant and $C_1$ depends only on $\Omega$.

3. (Case 3) If $n,h_{ex} \rightarrow \infty$, then up to
extraction
\begin{equation*}
 \frac{1}{n} \sum_{i=1}^n \delta_{\tilde{a}_i(\varepsilon)} \rightharpoonup \mu_0
\end{equation*}
in the narrow sense of measures, and $\mu_0$ is the unique
probability measure minimizing $I$, as defined by \eqref{I}.
The energy of these solutions as $\varepsilon \rightarrow 0$ is
given asymptotically by
\begin{equation*}
 G_\varepsilon(u_\varepsilon,A_\varepsilon) = f_\varepsilon(n) + n^2 I(\mu_0) + o(n^2).
\end{equation*}
We have
\begin{equation}\label{stable_solns_13}
 \pqnormspace{\nabla G_p}{2}{\infty}{\Omega} - o(1)
\le \frac{1}{n} \pqnormspace{\nabla_{A_\varepsilon'}
u_\varepsilon}{2}{\infty}{\Omega}
 \le \pqnormspace{\nabla G_p}{2}{\infty}{\Omega} + C + o(1),
\end{equation}
where $C$ is a universal constant. Finally, the convergence results of Corollary \ref{london_convergence} and Proposition \ref{l_z_convergence} hold.

\end{thm}

\begin{proof}
 The construction of the solutions and the proof of the energy asymptotics are done in Theorem 11.1 of \cite{ss_book}.  It remains to prove that in each case the result comparing $\pqnormspace{\nabla_{A'_\varepsilon} u_\varepsilon}{2}{\infty}{\Omega} $ to $n$ holds.  In the first and second case, the construction of the solutions is such that assumptions (J1) -- (J5) hold, and so we may apply Proposition \ref{bounded_n_thm}.  In the third case, (H1) -- (H4) are satisfied, and the asymptotics of $G_\varepsilon$ allow us to apply the second part of Theorem \ref{application_1} directly to conclude \eqref{stable_solns_13}.  The convergence results follow  directly from Corollary \ref{london_convergence} and Proposition \ref{l_z_convergence}.

\end{proof}

\begin{proof}[Proof of Theorem \ref{4}]
The result for branches of solutions follows immediately from the
above.
 
 We now discuss the case of energy-minimizers.
 In the regime $\llep\ll \he - \hci \ll \lep$,
 Proposition 9.1 of \cite{ss_book} establishes that minimizers
  $(u_\varepsilon,A_\varepsilon)$ of $G_\varepsilon$ satisfy $1\ll
  n \ll \he$,
\begin{equation*}
G_\varepsilon(u_\varepsilon,A_\varepsilon) = f_\varepsilon(n) +
n^2 I(\mu_0) + o(n^2),
\end{equation*}
and
\begin{multline*}
 F_\varepsilon(u_\varepsilon,A_\varepsilon') = \pi n \log{\frac{\ell}{\varepsilon}} + \pi S_\Omega(p,p)n^2 + \pi n^2 \log{\frac{1}{\ell}} \\
 - \pi n^2 \iint \log{\abs{x-y}} d\mu_0(x) d\mu_0(y) + o(n^2),
\end{multline*}
where $\mu_0$ is the minimizer of $I$ defined in \eqref{I}. The assumptions on $\alpha$ then
guarantee that  $F_\varepsilon(u_\varepsilon,A_\varepsilon') \le
\varepsilon^{\alpha-1}$, and so Theorem \ref{application_1} is
applicable.  The result follows.

For the case $\he - \hci \le O (\llep)$, it is proved in Theorem
12.1 of \cite{ss_book} that minimizers have $n$ vortices with $n =
O(1)$ and that they are among the solutions found in case 2 of
Theorem \ref{stable_solns}. Thus from that theorem, the
result holds in this case as well.

 For higher $\he$, in regime 4 when
$ \he \le \lep$, there is nothing new to
prove: the  upper bound follows  from bounds on $\|\nab_{A}
u\|_{L^2}$ and the fact that $n$ and $\he$ are of the same order.
The lower bound follows from the weak convergence of $j/\he$ to $-
\np h_*$ with $h_*$ nonconstant.

\end{proof}




\end{document}